\documentclass[preprint,12pt,authoryear]{elsarticle}

\usepackage{amsthm}
\usepackage{amsmath,amssymb,amsfonts}
\usepackage{mathtools}
\usepackage{xcolor}
\renewcommand{\textcolor}[2]{#2}

\usepackage{relsize} %
\usepackage{cancel}
\usepackage{algorithm}
\usepackage{algpseudocode}

\usepackage[utf8]{inputenc}
\usepackage[T1]{fontenc} 

\usepackage{graphicx} 
\usepackage{appendix}
\usepackage{bbm}
\usepackage[font=footnotesize,skip=0pt]{caption}
\usepackage{float}

\usepackage{subcaption}
\usepackage{comment}
\usepackage{booktabs}

\theoremstyle{definition}
\newtheorem{remark}{Remark}%
\renewcommand{\theremark}{\thesection.\arabic{remark}}
\usepackage{hyperref}
\hypersetup{
    colorlinks=true,       
    linkcolor=blue,        
    citecolor=ForestGreen, 
    filecolor=magenta,     
    urlcolor=cyan          
}

\usepackage{svg}
\usepackage{bbold}
\usepackage{multirow,makecell}
\usepackage{array}
\usepackage{comment}

\newtheorem{observation}{Observation}

\makeatletter
\def\ps@pprintTitle{%
   \let\@oddhead\@empty
   \let\@evenhead\@empty
   \def\@oddfoot{\hfill\thepage} 
   \let\@evenfoot\@oddfoot}
\makeatother

\begin{document}
\begin{frontmatter}

\title{A variable-offset joint formulation for beams with arbitrary cross-sections using a null space method}

\author[inst1]{Myung-Jin Choi\corref{cor1}}\ead{choi@lbb.rwth-aachen.de}


\cortext[cor1]{Corresponding author}

\affiliation[inst1]{organization={Chair of Structural Analysis and Dynamics, RWTH Aachen University},
            addressline={Mies-van-der-Rohe Str.\,1}, 
            city={Aachen},
            postcode={52074}, 
            country={Germany}}

\author[inst2,inst3,inst4]{Roger A. Sauer}
\author[inst1]{Simon Klarmann}
\author[inst1]{Sven Klinkel}


\affiliation[inst2]{organization={Institute for Structural Mechanics, Ruhr University Bochum},
            addressline={Universit{\"a}tsstra{\ss}e 150}, 
            city={Bochum},
            postcode={44801}, 
            country={Germany}}

\affiliation[inst3]{organization={Faculty of Civil and Environmental Engineering, Gda\'{n}sk University of Technology},
            addressline={ul. Narutowicza 11/12}, 
            city={Gda\'{n}sk},
            postcode={80-233}, 
            country={Poland}}

\affiliation[inst4]{organization={Department of Mechanical Engineering, Indian Institute of Technology Guwahati},
            addressline={Assam 781039}, 
            postcode={22222}, 
            country={India}}            

\begin{highlights}
    \item A local constraint formulation for an implicit interface with rotational degrees-of-freedom, and variable offset coordinates, considering arbitrarily shaped cross-sections
    \item A discrete null space method for size reduction and improved conditioning of the system matrix
    \item Application to an extensible director-based beam formulation, as well as a brick formulation
\end{highlights}
%
\begin{abstract}
In this paper, we present a variational formulation of local configurational constraints that couple multiple beams with arbitrarily shaped cross-sections. Since this formulation requires no explicit interface to rotational degrees-of-freedom, it applies to any beam kinematics and finite element discretization. Here, we define the offset coordinates in a moving frame to constrain or release the relative position between connected beams. The present method is based on a first-order approximation of the Lagrange multiplier field in the cross-section, which limits the transferability of the joint to the resultant force and moment only. The multipliers are eliminated using a discrete null space method, which provides size reduction and improved conditioning of the system matrix. Further, we apply the developed formulation to a beam element based on extensible directors and to a brick element in nonlinear elastostatics. Several numerical examples are presented.
\end{abstract}

\end{frontmatter}

\setcounter{remark}{0}
\section{Introduction}
In this work, we develop an efficient and stable approach to implementing various types of beam joints with the following three objectives: First, it should apply to beam kinematics, without an explicit interface to rotational degrees-of-freedom, as for example in so-called \textit{solid beam} formulations that use brick elements with translational degrees-of-freedom only, see, e.g., \citet{frischkorn2013solid}. \textcolor{blue}{Another type of such a formulation is a beam kinematics based on unconstrained directors}, see, e.g., \citet{antman1966dynamical}, \citet{rhim1998vectorial}, 
\citet{coda2009solid}, \citet{carrera2010refined}, \citet{durville2012contact}, \citet{moustacas2019enrichissement},
\citet{choi2021isogeometric}, \citet{choi2022isogeometric}). Second, the constraint is formulated in continuous form so that it can be applied at any point along the beam's axis, regardless of the finite element discretization. This can be especially useful for spatial discretization methods using non-interpolatory basis functions such as spline functions (e.g., see \citet{shafqat2024robust}, \citet{leonetti2025locking}) in the framework of what is known as \textit{isogeometric analysis} (IGA) \citep{hughes2005isogeometric}. Third, an additional treatment is necessary for the saddle-point system resulting from using the Lagrange multiplier method, since it suffers from (i) size-increase due to multipliers, and (ii) ill-conditioning. To address these issues, one can employ penalty-based formulations, such as augmented- and perturbed-Lagrangian methods, to eliminate the multipliers; see \citet{bertsekas2014constrained}, \citet{simo1985perturbed}, and references therein. In particular, for the static condensation of multipliers in the perturbed Lagrangian method, we refer to the recent discussion in \citet{duong2023variationally} and \citet[Section 3.5]{boungard2024master}. Here, our goal is to enforce the constraints exactly, without using an additional penalty parameter or iteration loop. For each objective, our current approach is in three steps. 

First, we introduce a constraint formulation using a Lagrange multiplier method, for an \textit{implicit interface} to the purely rotational and translational degrees-of-freedom of the cross-section at an arbitrary point along the axis. This method was presented in \citet{markovic2004micro} and later employed in \citet{klarmann2020homogenization}. Here, we refer to these additional translational and rotational degrees-of-freedom as \textit{interface variables}. \textcolor{blue}{Second, we parameterize the interface variables in terms of the joint's kinematic variables; translational, rotational, and offset degrees-of-freedom.} By fixing and releasing the offset variables, one can model a \textit{rigid} and \textit{sliding} offset joints, respectively. Further, having an offset can be useful for avoiding unphysical overlap between connected beams and for correctly representing the joint's stiffness \citep{di2024corotational}. In the present work, we further show that we may also release those offset variables, which can be useful in engineering applications, e.g., modeling of a prismatic joint. Third, for size-reduction and improved conditioning of the system matrix, we employ a \textit{null space} method, which has been well addressed in the literature, e.g., see \citet[Section 6]{benzi2005numerical}, and the relevant works include \citet{betsch2005discrete, betsch2006discrete} for rigid multibody dynamics, \citet{leyendecker2008discrete} for flexible multibody dynamics, and \citet{munoz2008modelling} and \citet{hesch2009mortar} for contact problems. After a finite element discretization, we construct an unconstrained (reduced) solution space for the system of linear algebraic equations. This so-called \textit{discrete null space} method requires finding a null space matrix for a discrete constraint Jacobian, which, in general, requires (i) an identification of independent rows (constraints), as well as (ii) independent columns (coordinates), from using a matrix factorization \citep{wehage1982generalized}. 

An accurate and robust method for enforcing constraints in mechanical systems has been the subject of extensive research, especially for multibody systems, see the literature review by \citet{bauchau2008review}. Given nodal degrees-of-freedom of translation and rotation, \citet{jelenic1996non} expressed the slave degrees-of-freedom in terms of the master and released ones, which is then utilized for reducing the system of linear equations. This method was later extended to \textcolor{blue}{dynamic} problems in \citet{jelenic2001dynamic}, and further applied to more complex joints in \citet{munoz2003master}. This so-called master-slave elimination approach is advantageous, since it provides (i) exact constraint enforcement, and (ii) a reduced number of unknown variables, see also the relevant work by \cite{boungard2024master} for general nonlinear multi-point constraints. However, this method requires an explicit interface to the nodal translational and rotational degrees-of-freedom as well as the determination of the master-slave pairs. This is not the case in the present formulation, which otherwise has the same advantages. It is therefore beneficial in the following perspectives: 
\begin{itemize}
    \item \textbf{General applicability}: No explicit interface to rotational degrees-of-freedom at nodes is required,
    \item \textbf{Unbiased formulation}: No master-slave relation between beams is required.
    \item \textbf{Straightforward calculation of constraint force/moment}: The constraint force and moment can be directly obtained from the Lagrange multipliers.
\end{itemize}
A constraint formulation for various joint types has also been extensively developed within the absolute nodal coordinate (ANC) framework, which employs nodal position and slope vectors as kinematic unknowns. This approach enables inter-element continuity of the displacement gradients even for standard $C^0$-finite element basis functions \citep{sugiyama2003formulation}, see also \citet{gerstmayr2013review} for an overview. The \textit{unbiasedness} of our present method can also be interpreted as the \textit{virtual body} concept for joints in \citet{bae2000implementation}, which provides a general framework for implementing various types of joints. This concept formulates the constraints not directly between the connected beams, but between each beam and a virtual body. It has been further extended to constrain the in-plane cross-sectional strains in \citet{sugiyama2011spatial}. Compared with those previous works, our contribution has the following novelties:
\begin{itemize}
    \item Constraints are formulated in an \textbf{averaged (integral) form} for the cross-section, which avoids artificial stress concentration, see \citet{gerstmayr20063d} for a relevant discussion for three-dimensional finite elements with linear strains in a co-rotated frame.
    \item For \textbf{continuous kinematical strains} along the length in the present beam formulation, we employ higher-order continuous spline basis functions (IGA), instead of having additional nodal slope variables.
    \item \textbf{Arbitrarily shaped cross-section} with corrected bending, torsional, and shear stiffness, without having additional correction factors. Here, we employ the \textit{local concept} of enriching warping strains from \citet{wackerfuss2011nonlinear}, which we adapt to account for existing constant in-plane cross-sectional strains from the extensible (unconstrained) directors. \textcolor{blue}{This will be explained in Appendix\,\ref{app_construct_warp_basis}.}
\end{itemize}        
In the present work, we have also imposed some simplifying restrictions; their generalization remains future work. These include: (i) no sliding along the beam's axis is assumed. That is, the joint constraints are applied at a given (fixed) material point on the beam's axis during the deformation. For a relevant further extension, one may refer to an incorporation of variable arc-length coordinate in \citet[Section 7.1]{sugiyama2003formulation}, and a variable length element in \citet{hong2011modeling}. Further, (ii) no relative rotation between beam and joint is assumed, i.e., revolute joints or hinges are not considered, which would require a further relaxation of the present constraints by an additional rotational variable, see, e.g., \citet{jelenic1996non}. (iii) We have also assumed no frictional resistance for the released offset coordinates. This can also be extended further by incorporating additional physics, see, e.g., the modeling of a spherical joint with clearance and lubrication in \citet{tian2009dynamics}. 

A constraint formulation for interface coupling conditions has also been extensively investigated, especially in \textit{rotation-free} formulations of structures such as thin beams and shells. For example, \citet{greco2021non} extracted rotational degrees-of-freedom from two end control points of a B{\'e}zier curve, from using the fact that the curve is tangent to the control net as well as interpolatory at the end point. This method was further extended in \citet{greco2024objective} to implement a cylindrical joint by introducing a pivot vector. \textcolor{blue}{A multiplicative decomposition was also applied in \citet{choi2023selectively} to extract rotational degrees-of-freedom from the extensible directors at the end control points, where it was also assumed that the spline curve is interpolatory at the end points.} \textcolor{blue}{Further, one can find an additional bending-resisting coupling element in \citet{kiendl2010bending}, enforcement of various interface conditions between shells using penalty and Lagrange multiplier methods in \citet{duong2017new}, a least-square approach and null space-based static condensation in \citet{schuss2019multi}, and a master-slave relation between local tangent frames in
\citet{bauer2020weak}.} In contrast to those previous formulations which limit the deformability of cross-sections, our formulation aims to deal with a \textit{partially clamped}\footnote{We adopt the terminology of partially and fully clamped conditions from \citet{hussein2009clamped}} joint, allowing deformation in the cross-sections. \textcolor{blue}{Furthermore, in the present work, to investigate the influence of clamped boundary and joint conditions, we have also implemented reference brick solutions with fully clamped conditions, which will be discussed in the numerical examples.}

The remainder of this paper is organized as follows: In Section\,\ref{jct_cfg_constraints}, we present a variational formulation for the interface constraints and a parameterization of the interface variables in terms of the joint's kinematic variables. \textcolor{blue}{In Section\,\ref{null_space_sec}, we present a discrete null space method for reducing the resulting system of algebraic equations.} In Section\,\ref{num_ex}, several numerical examples are presented. Section\,\ref{conclusions_jct} concludes the paper.  
\setcounter{remark}{0}
\section{Variable-offset joint formulation}
\label{jct_cfg_constraints}
In this section, we present a variational formulation for variable-offset joints. We first present the formulation in a general form, applicable to any beam kinematics. An application to an extensible director-based beam formulation can be considered as a special case. We begin with a general kinematic description of a beam. 
\subsection{Preliminary: a general kinematic description}
\begin{figure}[h]
	\centering
	\begin{subfigure}[b]{0.4\textwidth}\centering
		\includegraphics[width=\linewidth]{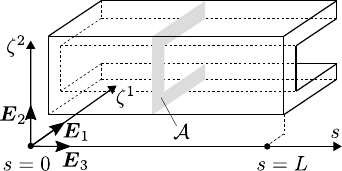}
		\caption{Reference configuration}
		\label{redraw_beam_kin_ref_3d}			
	\end{subfigure}
    \quad
    \begin{subfigure}[b]{0.55\textwidth}\centering
		\includegraphics[width=\linewidth]{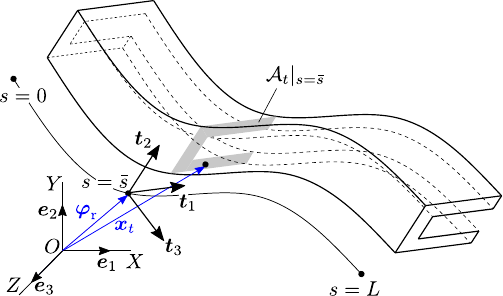}
		\caption{Current configuration}
		\label{redraw_beam_kin_cur_3d}			
	\end{subfigure}    
	\caption{A schematic illustration of the coordinate system, and the domain parameterization, in the (a) reference, and (b) physical domains.}
	\label{redraw_beam_kin_ref}	
\end{figure}

A beam is a three-dimensional slender body composed of a family of \textit{cross-sections} connected by an \textit{axis} curve. Therefore, the position of a material point in the beam's current configuration at time $t$ can be expressed by the mapping $\boldsymbol{x}_t\!\left(\zeta^1,\zeta^2,\zeta^3\right):\mathcal{A}\times\left[0,L\right]\rightarrow\Bbb{R}^3$, whose domain is called \textit{reference configuration}, where $\mathcal{A}\subset\Bbb{R}^2$ denotes the initially planar cross-sectional domain, and $L$ denotes the axis' initial length. The position of a material point in the reference configuration is given by
\begin{equation}
	{{\boldsymbol{X}}}=\zeta^i{\boldsymbol{E}}_i,
\end{equation}
with the standard Cartesian basis in $\mathbb{R}^3$, $\left\{\boldsymbol{E}_1,\boldsymbol{E}_2,\boldsymbol{E}_3\right\}$, where the cross-sectional domain $\mathcal{A}\ni{\left(\zeta^1,\zeta^2\right)}$ is spanned by two base vectors $\boldsymbol{E}_1$ and $\boldsymbol{E}_2$, and $\zeta^3\equiv s$ defines an arc-length coordinate, see Fig.\,\ref{redraw_beam_kin_ref_3d} for an illustration. Note that those coordinates $\zeta^1,\zeta^2$, and $s$ are time-independent. Here and hereafter, unless stated otherwise, repeated Latin indices like $i$ and $j$ imply summation over $1$ to $3$, and repeated Greek indices like $\alpha$ and $\beta$ imply summation over $1$ to $2$.
\begin{remark} \label{rem_area_mnt} \small For a given initial cross-section, we define the zeroth-, first-, and second-order moments of the cross-sectional area, as 
\begin{align*}
    A\coloneq \int_\mathcal{A} {{\rm{d}}{A}},\,\,I^\alpha\coloneq \int_\mathcal{A} {{\zeta ^\alpha }\,{\rm{d}}{A}},\,\,\mathrm{and}\,\,I^{\alpha\beta}\coloneq\int_\mathcal{A} {{\zeta ^\alpha }{\zeta ^\beta }\,{\rm{d}}{A}},
\end{align*}
respectively, where the indices $\alpha$ and $\beta$ range from $1$ to $2$, and $\mathrm{d}{A}\coloneq\mathrm{d}\zeta^1\mathrm{d}\zeta^2$ denotes the infinitesimal area. Here, the zeroth-order moment, $A$, represents the initial cross-sectional area.
\end{remark}
\noindent In the three-dimensional (physical) space, we have the standard Cartesian basis $\left\{\boldsymbol{e}_1,\boldsymbol{e}_2,\boldsymbol{e}_3\right\}$, with the corresponding coordinates $\left(X,Y,Z\right)\in\Bbb{R}^3$, whose origin is denoted by $O$. The position of a material point in the current cross-sectional domain $\mathcal{A}_t\subset\Bbb{R}^3$ is denoted by $\boldsymbol{x}_t\equiv\boldsymbol{x}_t\left(\zeta^1,\zeta^2,s\right)$. It should be noted that the cross-section may not remain plane, but be deformable. In the present joint formulation, our objective is to constrain only the \textit{purely rigid component} of a given cross-section at $s=\bar s$, a so-called \textit{partially clamped} joint. Therefore, we first need to introduce kinematical \textit{interface variables}, $\left(\boldsymbol{\varphi}_\mathrm{r},\boldsymbol{R}_\mathrm{r}\right)\in\Bbb{R}^3\times\mathrm{SO(3)}$ to describe the rigid motion. Here, SO(3) stands for the three-dimensional special orthogonal group, and $\boldsymbol{R}_\mathrm{r}\coloneqq\left[\boldsymbol{t}_1,\boldsymbol{t}_2,\boldsymbol{t}_3\right]$ represents the orientation of an attached orthonormal (rigid) frame, $\left\{\boldsymbol{t}_1,\boldsymbol{t}_2,\boldsymbol{t}_3\right\}$, such that $\boldsymbol{t}_i=\boldsymbol{R}_\mathrm{r}\boldsymbol{E}_i$ $(i=1,2,3)$, whose translation is given by the position vector $\boldsymbol{\varphi}_\mathrm{r}$, see Fig.\,\ref{redraw_beam_kin_cur_3d} for an illustration. In the following section, we explain how those interface variables can be implicitly connected to a general kinematic description of the cross-section.
\subsection{Implicit interface to the rigid cross-sectional motion}
Here, we introduce a constraint formulation for an interface to the rotational and translational degrees-of-freedom of the cross-section at an arbitrary point along the axis. This idea has been presented in \citet{markovic2004micro}, and also employed in \citet{klarmann2020homogenization}. For the given cross-section at $s=\bar s$, we consider a purely rigid component of the cross-section's motion, whose material point position is represented by
\begin{align}
	\label{mat_pt_pos_constrained_cs}
{\boldsymbol{x}}_\mathrm{r} &= \boldsymbol{\varphi}_\mathrm{r} + \zeta^\alpha\boldsymbol{t}_\alpha,\,\,(\zeta^1,\zeta^2)\in\mathcal{A}. 
\end{align}
where $\boldsymbol{\varphi}_\mathrm{r}\in\Bbb{R}^3$ and the orthonormal base vectors $\boldsymbol{t}_\alpha\in\Bbb{R}^3$ $(\alpha=1,2)$ represent the translational and rotational part of the cross-section's motion, respectively. Then, we introduce a local constraint,
\begin{align}
	\label{cfg_loc_cnst_general}
	{\boldsymbol{\Phi }} &\coloneqq \boldsymbol{x}_{\bar s} - \boldsymbol{x}_\mathrm{r} =\boldsymbol{0},\,\,(\zeta^1,\zeta^2)\in\mathcal{A},
\end{align}
where we have defined $\boldsymbol{x}_{\bar s}\coloneqq\boldsymbol{x}_t(\zeta^1,\zeta^2,\bar s)$ for brevity. Here, a point-wise (strong) enforcement of the constraint in Eq.\,(\ref{cfg_loc_cnst_general}) leads to a \textit{fully rigid} cross-section. One can relax this constraint such that only the purely rigid components of the cross-sectional motion are constrained, where the released higher-order components are associated with cross-sectional deformations. This relaxation will be further explained in the following section. 
\subsubsection{Lagrange multiplier method}
The constraint in Eq.\,(\ref{cfg_loc_cnst_general}) can be defined by the stationarity condition of the functional 
\begin{align}
	\label{lm_funct_pt_wise}
	\mathcal{J} &\coloneqq \int_\mathcal{A} {{\boldsymbol{\Lambda }} \cdot {\boldsymbol{\Phi }}\,{\rm{d}}A},
\end{align}
with respect to the Lagrange multiplier $\boldsymbol{\Lambda}\equiv\boldsymbol{\Lambda}(\zeta^1,\zeta^2)\in\Bbb{R}^3$. By assuming the following decomposition \citep{markovic2004micro,klarmann2018geometrisch},
\begin{align}
	\label{first_order_assump_lag_mult}
	{\boldsymbol{\Lambda }} = {{\boldsymbol{\lambda }}} + {{\boldsymbol{\mu}}}\times\boldsymbol{r},
\end{align}
where the Lagrange multipliers $\boldsymbol{\lambda}$ and $\boldsymbol{\mu}\in\Bbb{R}^3$ represent the constant stress, and a moment-like quantity which leads to a force acting perpendicular to the position vector $\boldsymbol{r}\coloneqq{{\zeta ^\alpha}{{\boldsymbol{t}}_\alpha }}$ in the cross-section, respectively. This, in other words, means that the transferability of the subsequent joint formulation is limited only to the resultant force
\begin{subequations}
\begin{align}
    {\Bbb{n}} \coloneqq \int_\mathcal{A} {{\boldsymbol{\Lambda }}\,{\rm{d}}A}  = A{\boldsymbol{\lambda }} + {\boldsymbol{\mu }} \times {I^\alpha }{{\boldsymbol{t}}_\alpha },
\end{align}
and the resultant moment
\begin{align}
    \label{resultant_mnt_cnst}
    {{\Bbb{m}}} \coloneqq \boldsymbol{t}_\alpha\times{{\Bbb{m}}^\alpha },
\end{align}
from the resultant couple
\begin{align}
    {{\Bbb{m}}^\alpha } \coloneqq \int_\mathcal{A} {{\zeta ^\alpha }{\boldsymbol{\Lambda }}\,{\rm{d}}A}  = {I^\alpha }{\boldsymbol{\lambda }} + {\boldsymbol{\mu }} \times {I^{\alpha \beta }}{{\boldsymbol{t}}_\beta },\,\,\mathrm{for}\,\,\alpha=1,2.
\end{align}
\end{subequations}
\noindent Substituting Eq.\,(\ref{first_order_assump_lag_mult}) into Eq.\,(\ref{lm_funct_pt_wise}) leads to 
\begin{align}
	\label{lm_funct_beam_kin}
	\mathcal{J} & ={\boldsymbol{\lambda }} \cdot {\boldsymbol{\phi}} + {\boldsymbol{\mu }} \cdot {\boldsymbol{\psi}},
\end{align}
where we have defined the following relaxed (resultant) constraint functions in the cross-section,
\begin{subequations}
\label{cnst_func_general}
\begin{align}
	\label{cnst_func_p_phi}
	{\boldsymbol{\phi}} &\coloneqq \int_\mathcal{A} {{\boldsymbol{\Phi}}\,{\rm{d}}{A}} = \int_\mathcal{A} {\boldsymbol{x}}_{\bar s} \,{\rm{d}}A - {I^\alpha }{{\boldsymbol{t}}_\alpha } - A{{\boldsymbol{\varphi }}_{\rm{r}}} =\boldsymbol{0},
\end{align}
and
\begin{align}
	{\boldsymbol{\psi}} &\coloneqq {{{\boldsymbol{t}}}_\alpha } \times \boldsymbol{\psi}^\alpha =\boldsymbol{0},\label{cnst_func_p_tht}
\end{align}
\end{subequations}
with 
\begin{align}
    \label{cnst_func_psi_general}
	{\boldsymbol{\psi}^\alpha } \coloneqq \int_\mathcal{A} \zeta^\alpha{{\boldsymbol{\Phi}}\,{\rm{d}}{A}} = \int_\mathcal{A} {{\zeta ^\alpha }{\boldsymbol{x}}}_{\bar s} \,{\rm{d}}A - {I^{\alpha \beta }}{{\boldsymbol{t}}_\beta } - {I^\alpha }{{\boldsymbol{\varphi }}_{\rm{r}}}.	
\end{align}
It is seen that, from Eq.\,(\ref{cfg_loc_cnst_general}), we finally have three constraints in Eq.\,(\ref{cnst_func_p_phi}) for the translations, and another three in Eq.\,(\ref{cnst_func_p_tht}) for the rotations, which provides an implicit interface to the variables $\left(\boldsymbol{\varphi}_\mathrm{r},\boldsymbol{R}_\mathrm{r}\right)\in\Bbb{R}^3\times\mathrm{SO(3)}$. In the following, we parameterize those six constrained cross-sectional degrees-of-freedom by the newly introduced joint kinematic variables and offsets.
\textcolor{blue}{
\begin{remark} \label{rem_non_zero_Phi} \small The constraints $\boldsymbol{\phi}=\boldsymbol{\psi}=\boldsymbol{0}$ in Eqs.\,(\ref{cnst_func_p_phi}) and (\ref{cnst_func_p_tht}) do not necessarily enforce the point-wise satisfaction of $\boldsymbol{\Phi}(\zeta^1,\zeta^2)=\boldsymbol{0},\,\,(\zeta^1,\zeta^2)\in\mathcal{A}$. That is, they allow all the functions $\boldsymbol{x}_{\bar s}$ whose integral satisfy Eqs.\,(\ref{cnst_func_p_phi}) and (\ref{cnst_func_p_tht}). Here, the value $\boldsymbol{\Phi}\ne\boldsymbol{0}$ measures the deviation from the purely rigid motion due to $\left(\boldsymbol{\varphi}_\mathrm{r},\boldsymbol{R}_\mathrm{r}\right)$. 
\end{remark}
}
\begin{figure}[h]
	\centering
	\begin{subfigure}[b]{0.475\textwidth}\centering
		\includegraphics[width=\linewidth]{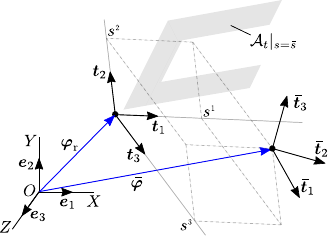}
	\end{subfigure}
	\caption{A schematic illustration of the interface variables ($\boldsymbol{\varphi}_\mathrm{r},\boldsymbol{R}_\mathrm{r}$), the joint's kinematic variables ($\boldsymbol{\bar\varphi},\boldsymbol{\bar R}$), and offset coordinates $s^i$.}
	\label{redraw_beam_kin_ref_current_shift}
\end{figure}
\subsection{Joint kinematics}
We consider a joint with attached set of orthonormal base vectors $\boldsymbol{\bar t}_i\in\Bbb{R}^3$ ($i=1,2,3$), a so-called \textit{joint frame}, whose translation is described by the position vector $\boldsymbol{\bar \varphi}\in\Bbb{R}^3$, which can be assumed to describe the center of mass. First, for the translational motion, we introduce three variable offset coordinates $\left\{s^1,s^2,s^3\right\}$,
\begin{align}
	\label{jct_connect_translate_gen}
    {\boldsymbol{\varphi}}_\mathrm{r} \equiv \boldsymbol{\bar \varphi} - s^i\boldsymbol{t}_i,
\end{align}
as illustrated in Fig.\,\ref{redraw_beam_kin_ref_current_shift}. It is noted that the offset coordinates $s^i$ are defined in the local frame $\left\{{\boldsymbol{t}}_1,{\boldsymbol{t}}_2,{\boldsymbol{t}}_3\right\}$, in order to constrain or release the \textit{relative} position of connected beams with respect to the joint. By having the offsets, we have the following advantages: 
\begin{itemize}
    \item Correct rigidity of the joint: Unphysical initial overlap between connected beams can be avoided,
    \item Releasing a variable offset coordinate enables modeling a prismatic\footnote{This term refers to \citet{angeles1982spatial} and \citet{jelenic1996non}. Here, the `sliding' means a relative translation between a \textit{fixed} pair of cross-sections.} (sliding) joint.
\end{itemize}
Second, we also relate the rigid rotation of a given cross-section, described by the orthonormal frame $\left\{\boldsymbol{t}_1,\boldsymbol{t}_2,\boldsymbol{t}_3\right\}$, to the \textit{joint frame}, as
\begin{align}
	\label{rot_param_tvec_bar}
	\boldsymbol{t}_i = \boldsymbol{R}_0\,\boldsymbol{\bar t}_i,\,\,\boldsymbol{R}_0\in\text{SO(3)},\,\,\mathrm{for}\,\,i=1,2,3.
\end{align}
In this work, we consider that the rotation tensor $\boldsymbol{R}_0$ is constant, without allowing any relative rotation between the joint and the connected cross-section, which means a rigid joint. Here, the subscript $0$ means that the rotation tensor can be determined from the initial configuration. It is believed that one can straightforwardly extend the present formulation by introducing an additional degree-of-freedom for releasing the relative rotation, e.g., a hinge, revolute, or cylindrical joint. This is left for future work. Hereafter, for a compact notation, we introduce the following notation:
\begin{itemize}
    \item $\Bbb{c}\coloneqq\left[\boldsymbol{\phi}^\mathrm{T},\boldsymbol{\psi}^\mathrm{T}\right]^\mathrm{T}\in\Bbb{R}^m$ for the constraint functions, where $m=6$ denotes the number of constraints,
    \item $\Bbb{f}\coloneqq\left[\boldsymbol{\lambda}^\mathrm{T},\boldsymbol{\mu}^\mathrm{T}\right]^\mathrm{T}\in\Bbb{R}^m$ for the Lagrange multipliers,    
    \item $\boldsymbol{\bar q}\coloneqq\left[\boldsymbol{\bar\varphi}^\mathrm{T},\boldsymbol{\bar t}_1^{\,\mathrm{T}},\boldsymbol{\bar t}_2^{\,\mathrm{T}},\boldsymbol{\bar t}_3^\mathrm{T}\right]^{\,\mathrm{T}}\in\Bbb{R}^{m+6}$ , for the kinematic variables of a joint, where those additional six dimensions are constrained due to $\boldsymbol{\bar R}\coloneqq\left[ {{{{\boldsymbol{\bar t}}}_1},{{{\boldsymbol{\bar t}}}_2},{{{\boldsymbol{\bar t}}}_3}} \right]\in\text{SO(3)}$,
    \item $\boldsymbol{s}\coloneqq\left[s^1,s^2,s^3\right]^\mathrm{T}\in\Bbb{R}^3$ for the variable offset coordinates, defined in the local frame $\left\{{{{{\boldsymbol{t}}}_1},{{{\boldsymbol{t}}}_2},{{{\boldsymbol{t}}}_3}}\right\}$.
\end{itemize}
Using the relations in Eqs.\,(\ref{jct_connect_translate_gen}) and (\ref{rot_param_tvec_bar}), the interface variables in Eq.\,(\ref{cnst_func_general}) can be replaced by $\boldsymbol{\bar q}$, so that the constraint equations can be rewritten as
%
\begin{align}
    \label{cnst_cc_xs}
    \Bbb{c}\equiv\Bbb{c}\!\left(\boldsymbol{x}_{\bar s},\boldsymbol{\bar q},\boldsymbol{s}\right)=\boldsymbol{0}.
\end{align}
%
From the following variational formulation, the constraint force and moment at the joint, along with their equilibrium equations, can be identified. \textcolor{blue}{
\begin{remark} 
\label{rem_coupling_number_constraints} \small
The present joint formulation is developed in an unbiased manner, in which the coupling (interface) constraints are defined between the joint and a single connected beam, rather than directly describing the relation between the connected beams. This allows us, in what follows, to describe the formulation for a single connected beam, without loss of generality. A finite element assembly will then result in coupling between the connected beams at the joint.
\end{remark}}
\subsection{Variational formulation}
Here, for simplicity, we assume that $\bar s$ is given (fixed), i.e., $\delta \bar s=0$. This implies that the joint constraint is always applied to the same material point during the deformation. That is, no relative sliding of the joint along the beam's axis is considered. Taking the first variation of $\mathcal{J}$ in Eq.\,(\ref{lm_funct_beam_kin}), after substituting Eqs.\,(\ref{jct_connect_translate_gen}) and (\ref{rot_param_tvec_bar}), gives
\begin{align}
	\label{fvar_functional_J_cnst}
	\delta \mathcal{J} &=\delta {{\boldsymbol{x}}_{\bar s}} \cdot {{\boldsymbol{Q}}_{\rm{x}}}  + \delta {\boldsymbol{\bar q}} \cdot {{\boldsymbol{Q}}_{{\rm{\bar q}}}} + \delta {\boldsymbol{s}} \cdot {{\boldsymbol{Q}}_{\rm{s}}} + \delta {\Bbb{f}} \cdot {\Bbb{c}},
\end{align}
where we have defined the \textit{generalized constraint forces}
\begin{subequations}
    \label{general_cnst_f_Q}
    \begin{alignat}{5}
    {{\boldsymbol{Q}}_{\rm{x}}} &\coloneqq\dfrac{\partial\mathcal{J}}{\partial\boldsymbol{x}_\mathrm{\bar s}}&&= {\boldsymbol{G}}_{\rm{x}}^{\rm{T}}{\Bbb{f}} = \left({\int_\mathcal{A} \left(\bullet\right)\,{{\rm{d}}A} }\right) {\boldsymbol{\lambda }} + \left(\int_\mathcal{A} {{\zeta ^\alpha}\left(\bullet\right)\,{\rm{d}}A}\right)\left( {{\boldsymbol{\mu }}\times {{\boldsymbol{t}}_\alpha } } \right),\\
    {{\boldsymbol{Q}}_{{\rm{\bar q}}}}&\coloneqq\dfrac{\partial\mathcal{J}}{\partial\boldsymbol{\bar q}} &&= {\boldsymbol{G}}_{{\rm{\bar q}}}^{\rm{T}}{\Bbb{f}}={\left[ {{{{\Bbb{\bar n}}}^{\rm{T}}},{\Bbb{\bar m}}_{\rm{t}}^{1\,{\rm{T}}},{\Bbb{\bar m}}_{\rm{t}}^{{2\,\rm{T}}},{\Bbb{\bar m}}_{\rm{t}}^{{3\,\rm{T}}}} \right]^{\rm{T}}},\label{cnst_f_qbar} \\ 
    {{\boldsymbol{Q}}_{\rm{s}}} &\coloneqq\dfrac{\partial \mathcal{J}}{\partial\boldsymbol{s}}&&= {\boldsymbol{G}}_{\rm{s}}^{\rm{T}}{\Bbb{f}}=\left[\boldsymbol{t}_1,\boldsymbol{t}_2,\boldsymbol{t}_3\right]^\mathrm{T}\Bbb{n},   
    \end{alignat}
\end{subequations}
with the constraint Jacobians,
\begin{subequations}
    \label{cnst_jcb_x_q_s}
    \begin{align}
    {{\boldsymbol{G}}_{\rm{x}}} &\coloneqq \dfrac{{\partial {\Bbb{c}}}}{{\partial {{\boldsymbol{x}}_{\bar s}}}},\\
    {{\boldsymbol{G}}_{{\rm{\bar q}}}} &\coloneqq  \dfrac{{\partial {\Bbb{c}}}}{{\partial {{\boldsymbol{\bar q}}}}},\\
    {{\boldsymbol{G}}_{{\rm{s}}}} &\coloneqq  \dfrac{{\partial {\Bbb{c}}}}{{\partial {{\boldsymbol{s}}}}},    
        \end{align}
\end{subequations}
whose detailed expressions can be found in Appendix \ref{app_cnst_jcb_general}. Further, in Eq.\,(\ref{cnst_f_qbar}), we have defined the constraint force
\begin{align}
    \label{const_force_n_bar}
    {\Bbb{\bar n}} \coloneqq \dfrac{{\partial \mathcal{J}}}{{\partial {\boldsymbol{\bar \varphi }}}} =  - {\Bbb{n}},
\end{align}
and the constraint couple
\begin{align}
    {\Bbb{\bar m}}_{\rm{t}}^i \coloneqq 
    \dfrac{\partial\mathcal{J}}{\partial \boldsymbol{\bar t}_i}={\boldsymbol{R}}^\mathrm{T}_0\,{\Bbb{m}}_{\rm{t}}^i,
\end{align}
for $i=1,2,3$, with
\begin{align}
    {\Bbb{m}}_{\rm{t}}^i \coloneqq \dfrac{{\partial \mathcal{J}}}{{\partial {{\boldsymbol{t}}_i}}} = {s^i}{\Bbb{n}} - \delta _\alpha ^i\left( {{{\Bbb{m}}^\alpha } + {\boldsymbol{\mu }} \times {{\boldsymbol{\psi }}^\alpha }} \right).
\end{align}
\subsubsection{Rotational parameterization of the joint frame}
We further introduce three rotational parameters through a rotational pseudo vector $\boldsymbol{\bar \theta}\in\Bbb{R}^3$ \citep{argyris1982excursion}, due to the constraint $\left[\boldsymbol{\bar t}_1,\boldsymbol{\bar t}_2,\boldsymbol{\bar t}_3\right]\in\text{SO(3)}$, such that
\begin{align}
	\label{fvar_unit_t_i}
	\delta {{\boldsymbol{\bar t}}_i} = \delta {\boldsymbol{\bar \theta }} \times {{\boldsymbol{\bar t}}_i},\,\,\mathrm{for}\,\,i=1,2,3.
\end{align}
Then, Eq.\,(\ref{fvar_functional_J_cnst}) can be rewritten, as
\begin{align}
    \label{fvar_J_cnst_jcb_rewrite_fc}
	\delta \mathcal{J} &={\delta {{\boldsymbol{x}}_{\bar s}} \cdot {{\boldsymbol{Q}}_{\rm{x}}}} + {{\delta {\Bbb{\bar q}} \cdot \boldsymbol{Q}_{\Bbb{\bar q}}}} + {\delta {\boldsymbol{s}} \cdot {{\boldsymbol{Q}}_{\rm{s}}}} + {\delta {\Bbb{f}} \cdot {\Bbb{c}}},
\end{align}
where we have defined $\Bbb{\bar q}\coloneqq\left[\boldsymbol{\bar \varphi}^\mathrm{T},\boldsymbol{\bar \theta}^\mathrm{T}\right]^\mathrm{T}\in\Bbb{R}^m$, and \textcolor{blue}{the operator
\begin{align}
    {\boldsymbol{\bar\varXi}} \coloneqq \left[ {\renewcommand{\arraystretch}{1.5}\begin{array}{*{20}{c}}
{\bf{1}}_3&{\bf{0}}_{3\times3}&{\bf{0}}_{3\times3}&{\bf{0}}_{3\times3}\\
{{\bf{0}}_{3\times3}}&{{{\left[ {{{{\boldsymbol{\bar t}}}_1}} \right]}_ \times }}&{{{\left[ {{{{\boldsymbol{\bar t}}}_2}} \right]}_ \times }}&{{{\left[ {{{{\boldsymbol{\bar t}}}_3}} \right]}_ \times }}
\end{array}} \right]^\mathrm{T},
\end{align}
such that $\delta\boldsymbol{\bar q}=\boldsymbol{\bar\varXi}\delta\Bbb{\bar q}$.} Here and hereafter, ${\bf{1}}_n$ denotes the $n$-dimensional identity matrix, and ${\bf{0}}_{m\times{n}}$ denotes the $m$ by $n$ zero matrix. Then, we have the \textit{generalized constraint force}
\begin{align}
\boldsymbol{Q}_{\Bbb{\bar q}}\coloneqq\dfrac{\partial\mathcal{J}}{\partial\Bbb{\bar q}}=\boldsymbol{G}_\Bbb{\bar q}^\mathrm{T}\Bbb{f}=\left[\Bbb{\bar n}^\mathrm{T},\Bbb{\bar m}^\mathrm{T}\right]^\mathrm{T},
\end{align}
with the constraint Jacobian
\begin{align}
\boldsymbol{G}_\Bbb{\bar q}\coloneqq\dfrac{\partial\Bbb{c}}{\partial\Bbb{\bar q}}={{\boldsymbol{G}}_{{\rm{\bar q}}}}\,\boldsymbol{\bar \varXi}.
\end{align}
\textcolor{blue}{
It is noted that the generalized force $\boldsymbol{Q}_{\Bbb{\bar q}}$ consists of the constraint force $\Bbb{\bar n}$ in Eq.\,(\ref{const_force_n_bar}), and the \textit{constraint moment}
\begin{align}
    \label{def_constraint_mnt}
	\Bbb{\bar m} &\coloneqq \dfrac{\partial\mathcal{J}}{\partial{\boldsymbol{\bar\theta}}}={{\boldsymbol{\bar t}}_i}\times\Bbb{\bar m}^i_\mathrm{t}=\boldsymbol{R}^\mathrm{T}_0\left(\boldsymbol{t}_i\times\Bbb{m}^i_\mathrm{t}\right),
\end{align}
where the multiplication of $\boldsymbol{R}_0$ accounts for the initial relative orientation between the beam and joint.}
\begin{remark} \small The constraint moment in Eq.\,(\ref{def_constraint_mnt}) can be rewritten, as
\begin{align}
    \label{res_mnt_mb_def_0}
    \Bbb{\bar m}  = {{\boldsymbol{R}}^\mathrm{T}_0}\left( \underbrace{{\boldsymbol{\phi} \times {\boldsymbol{\lambda }} + {\boldsymbol{\psi }} \times {\boldsymbol{\mu }}}}_{\text{non-equilibrium}}
     - {\Bbb{m}} + \underbrace{{s^i}{{\boldsymbol{t}}_i} \times {\Bbb{n}}}_{\text{offset}} - \underbrace{\int_\mathcal{A} {{\boldsymbol{\Phi }} \times {\boldsymbol{\Lambda }}\,{\rm{d}}{A}}}_{{\text{higher-order modes}}}\right).
\end{align}
Here, in the parentheses, the first two terms vanish at equilibrium, due to ${\boldsymbol{\phi}}=\boldsymbol{\psi}=\boldsymbol{0}$, and the fourth and fifth terms represent the additional moment due to the offset, and the released higher-order cross-sectional deformation modes (i.e., $\boldsymbol{\Phi}\ne\boldsymbol{0}$, see Remark \ref{rem_non_zero_Phi}), respectively. Further, the third term implies that the constraint moment is a reaction to the resultant moment $\Bbb{m}$ in Eq.\,(\ref{resultant_mnt_cnst}).
\end{remark}
\subsubsection{External virtual work}
We define a \textit{work function} for a time-dependent external load applied to the joint,
\begin{align}
	\label{ext_work_func_jct}
	{\overline W}_\mathrm{ext} \coloneqq{\boldsymbol{\bar f}}_\mathrm{ext}(t) \cdot {\boldsymbol{\bar \varphi }} + {{\boldsymbol{\bar m}}_\mathrm{ext}^i}(t) \cdot {{\boldsymbol{\bar t}}_i},
\end{align}
where $t$ denotes time, and $\boldsymbol{\bar f}_\mathrm{ext}$ and $\boldsymbol{\bar m}_\mathrm{ext}^i$ $(i=1,2,3)$ denote the external force and director couple, respectively, which are, here, assumed deformation-independent. Then, taking the first variation of Eq.\,(\ref{ext_work_func_jct}), and using Eq.\,(\ref{fvar_unit_t_i}), we obtain the external virtual work, 
\begin{align}
	\delta{{\overline W}_\mathrm{ext}} = \delta\Bbb{\bar q}\cdot{\boldsymbol{Q}_\mathrm{ext}}(t),
\end{align}
where we have defined the prescribed external load ${\boldsymbol{Q}_\mathrm{ext}}\coloneqq\left[\boldsymbol{\bar f}_\mathrm{ext}^{\,\mathrm{T}},{\boldsymbol{\bar m}}_\mathrm{ext}^{\,\mathrm{T}}\right]^\mathrm{T}\in\Bbb{R}^m$, with the \textit{external moment},
\begin{align}
	\boldsymbol{\bar m}_\mathrm{ext}\coloneqq\boldsymbol{\bar t}_i\times\boldsymbol{\bar m}_\mathrm{ext}^i.
\end{align} 

\subsubsection{\textcolor{blue}{Equilibrium equations}}
\textcolor{blue}{
For a given load at time $t$, we have the following balance of momentum:
\begin{subequations}
    \label{var_eq_general}
	\begin{align}
		G_\mathrm{int}\!\left(\boldsymbol{x},\delta\boldsymbol{x}\right) + \delta\boldsymbol{x}_{\bar s}\cdot\boldsymbol{Q}_\mathrm{x}= G_\mathrm{ext}\!\left(\delta\boldsymbol{x}\right),\,\,\forall\delta\boldsymbol{x}\in\mathcal{V},\label{lin_mnt_balance}
	\end{align}
    subject to the joint equilibrium conditions
	\begin{alignat}{5}
		\boldsymbol{Q}_\mathrm{s} &=\boldsymbol{0},\quad&&\text{(free sliding)}\label{equilibrium_release_offset} \\
        \boldsymbol{Q}_{\Bbb{\bar q}} &= {\boldsymbol{Q}_\mathrm{ext}}(t),\quad&&(\text{joint equilibrium})
	\end{alignat}		
	and the interface constraints,
	\begin{alignat}{5}
		\Bbb{c} =\boldsymbol{0}.
	\end{alignat}	
\end{subequations}
\noindent In the present work, we consider a problem in nonlinear elastostatics, where ${G}_\mathrm{int}$ and ${G}_\mathrm{ext}$ represent the internal and external virtual works, respectively, whose expressions can be found in literature on finite element method for solids, e.g., \citet{bonet2010nonlinear}. Further, $\mathcal{V}$ denotes the typical variational space for the displacement field \citep{hughes2003finite}. In Eq.\,(\ref{equilibrium_release_offset}), we have assumed free sliding due to the released offset coordinates, i.e., no constraint force due to friction is considered. 
}
\subsection{Linearization}
For an iterative solution process, we need to linearize the governing variational equations in Eq.\,(\ref{var_eq_general}). The increment of the first variation in Eq.\,(\ref{fvar_functional_J_cnst}) can be expressed as
\begin{align}
	\label{del_J_tangent_k}
		\Delta \delta \mathcal{J} &= \underbrace{\delta {{\boldsymbol{x}}_{\bar s}} \cdot \Delta {{\boldsymbol{Q}}_{\rm{x}}} + \delta {\boldsymbol{\bar q}} \cdot \Delta {{\boldsymbol{Q}}_{{\boldsymbol{\bar q}}}} + \delta {\boldsymbol{s}} \cdot \Delta {{\boldsymbol{Q}}_{\rm{s}}}}_\text{material part} + \underbrace{\Delta\!\left(\delta {{{\boldsymbol{\bar t}}}_i}\right) \cdot{\Bbb{\bar m}}_\mathrm{t}^i}_\text{geometric part} + \delta {\Bbb{f}} \cdot \Delta {\Bbb{c}}\nonumber\\
		&=\left\{ {\renewcommand{\arraystretch}{1.5}\begin{array}{*{20}{c}}
				{\delta {\boldsymbol{x}}_\mathrm{\bar s}}\\
                {\delta {\Bbb{\bar q}}}\\
                {\delta {\boldsymbol{s}}}\\
                {\delta {\Bbb{f}}}
		\end{array}} \right\}.\underbrace{\left[ {\renewcommand{\arraystretch}{1.5}\begin{array}{*{20}{c}}
{{\bf{0}}_{3\times3}}&{{{\boldsymbol{k}}_{{\rm{x}\Bbb{\bar q}}}}}&{{\bf{0}}_{3\times3}}&{{{\boldsymbol{G}}^\mathrm{T}_{{\rm{x}}}}}\\
{{\boldsymbol{k}}_{{\rm{x}\Bbb{\bar q}}}^{\rm{T}}}&{{{\boldsymbol{k}}_{{\Bbb{\bar q}\Bbb{\bar q}}}}}&{{{\boldsymbol{k}}_{{\Bbb{\bar q}}\rm{s}}}}&{\boldsymbol{G}^\mathrm{T}_{\Bbb{\bar q}}}\\
{{\bf{0}}_{3\times3}}&{{\boldsymbol{k}}_{{\Bbb{\bar q}}\rm{s}}^{\rm{T}}}&{{{\bf{0}}_{3\times3}}}&{\boldsymbol{G}^\mathrm{T}_{\mathrm{s}}}\\
{{\boldsymbol{G}}_{\rm{x}}}&{{\boldsymbol{G}}_{{\Bbb{\bar q}}}}&{{\boldsymbol{G}}_{{\rm{s}}}}&{{\bf{0}}_{m\times{m}}}
\end{array}} \right]}_{\eqqcolon\boldsymbol{k}}\left\{ {\renewcommand{\arraystretch}{1.5}\begin{array}{*{20}{c}}
				{\Delta {\boldsymbol{x}}_\mathrm{\bar s}}\\
                {\Delta {\Bbb{\bar q}}}\\
                {\Delta {\boldsymbol{s}}}\\   
				{\Delta {\Bbb{f}}}
		\end{array}} \right\},
\end{align}
where we have defined
\begin{subequations}
\begin{alignat}{5}
    \boldsymbol{k}_{\mathrm{x}\Bbb{\bar q}} &\coloneqq \dfrac{{\partial {\boldsymbol{Q}_{\rm{x}}}}}{{\partial {\Bbb{\bar q}}}} &&= {\left( {\dfrac{{{\partial ^2}{\Bbb{c}}}}{{\partial {{\boldsymbol{x}}_{\bar s}}\partial {\Bbb{\bar q}}}}} \right)^{\!\rm{T}}}{\Bbb{f}},\\
    \boldsymbol{k}_{\mathrm{s}\Bbb{\bar q}} &\coloneqq \dfrac{{\partial {\boldsymbol{Q}_{\rm{s}}}}}{{\partial {\Bbb{\bar q}}}} &&= {\left( {\dfrac{{{\partial ^2}{\Bbb{c}}}}{{\partial {{\boldsymbol{s}}}\partial {\Bbb{\bar q}}}}} \right)^{\!\rm{T}}}{\Bbb{f}},\\
    \boldsymbol{k}_{\Bbb{\bar q}\Bbb{\bar q}} &\coloneqq \dfrac{{\partial {\boldsymbol{Q}_{\Bbb{\bar q}}}}}{{\partial {\Bbb{\bar q}}}} &&= {\left( {\dfrac{{{\partial ^2}{\Bbb{c}}}}{{\partial {\Bbb{\bar q}}\partial {\Bbb{\bar q}}}}} \right)^{\!\rm{T}}}{\Bbb{f}},
\end{alignat}
\end{subequations}
whose detailed expressions can be found in Appendix \ref{app_tangent_oper}. \textcolor{blue}{Here and hereafter, $\Delta(\bullet)$ denotes an increment.}
\begin{remark} 
\label{rem_sym_tstiff_geom} \small
Here, the term \textit{material} part implies that the tangent operators come from the derivative of the generalized constraint forces in Eq.\,(\ref{general_cnst_f_Q}). The \textit{geometric} part arises due to the nonlinearity of SO(3), and appears in $\boldsymbol{k}_{\Bbb{\bar q}\Bbb{\bar q}}$ only, see Eq.\,(\ref{k_tht_tht_decompose}), whose unsymmetric part vanishes at the (moment) equilibrium \citep{simo1986three}, see Appendix\,\ref{app_tangent_oper} for the proof.
\end{remark}
\noindent The tangent operator $\boldsymbol{k}$ in Eq.\,(\ref{del_J_tangent_k}) can be rewritten, as
\begin{align}
    \boldsymbol{k}\coloneqq\left[ {\begin{array}{*{20}{c}}
{\Bbb{k}}&{{{\boldsymbol{G}}^{\rm{T}}}}\\
{\boldsymbol{G}}&{{\bf{0}}_{m\times{m}}}
\end{array}} \right],
\end{align}
where
\begin{align}
    {\Bbb{k}} \coloneqq \left[ {\renewcommand{\arraystretch}{1.5}\begin{array}{*{20}{c}}
{\bf{0}}_{3\times{3}}&{{{\boldsymbol{k}}_{{\rm{x}\Bbb{\bar q}}}}}&{\bf{0}}_{3\times{3}}\\
{{\boldsymbol{k}}_{{\rm{x}\Bbb{\bar q}}}^{\rm{T}}}&{{{\boldsymbol{k}}_{{\Bbb{\bar q\bar q}}}}}&{{{\boldsymbol{k}}_{{\Bbb{\bar q}\rm{s}}}}}\\
{\bf{0}}_{3\times{3}}&{{\boldsymbol{k}}_{\Bbb{\bar q}{\rm{s}}}^{\rm{T}}}&{\bf{0}}_{3\times{3}}
\end{array}} \right],
\end{align}
and the constraint Jacobian,
\begin{align}
    {\boldsymbol{G}} \coloneqq \left[ {\begin{array}{*{20}{c}}
{{{\boldsymbol{G}}_{\rm{x}}}}&{{{\boldsymbol{G}}_{{\Bbb{\bar q}}}}}&{{{\boldsymbol{G}}_{\rm{s}}}}
\end{array}} \right].
\end{align}
\begin{remark} \small \textit{Application of first-order beam kinematics with extensible-directors}. \textcolor{blue}{In Appendix \ref{app_beam_form}, we explain the application of the first-order beam kinematics to the present general form of constraint formulation.} Note that the six constraint equations in Eq.\,(\ref{cnst_cc_xs}) constrain purely rigid motion only, which is the so-called partially clamped condition. That is, they allow in-plane cross-sectional deformations; two transverse normal strains, and one pure shear strain; see Fig.\,\ref{beam_cs_strn_mode_allow} for an illustration. Further, the enhanced warping strains ($\boldsymbol{\widetilde E}$), presented in Appendix \ref{app_construct_warp_basis}, are also freely allowed. 
\begin{figure}[h]
	\centering
	\begin{subfigure}[b]{0.625\textwidth}\centering
		\includegraphics[width=\linewidth]{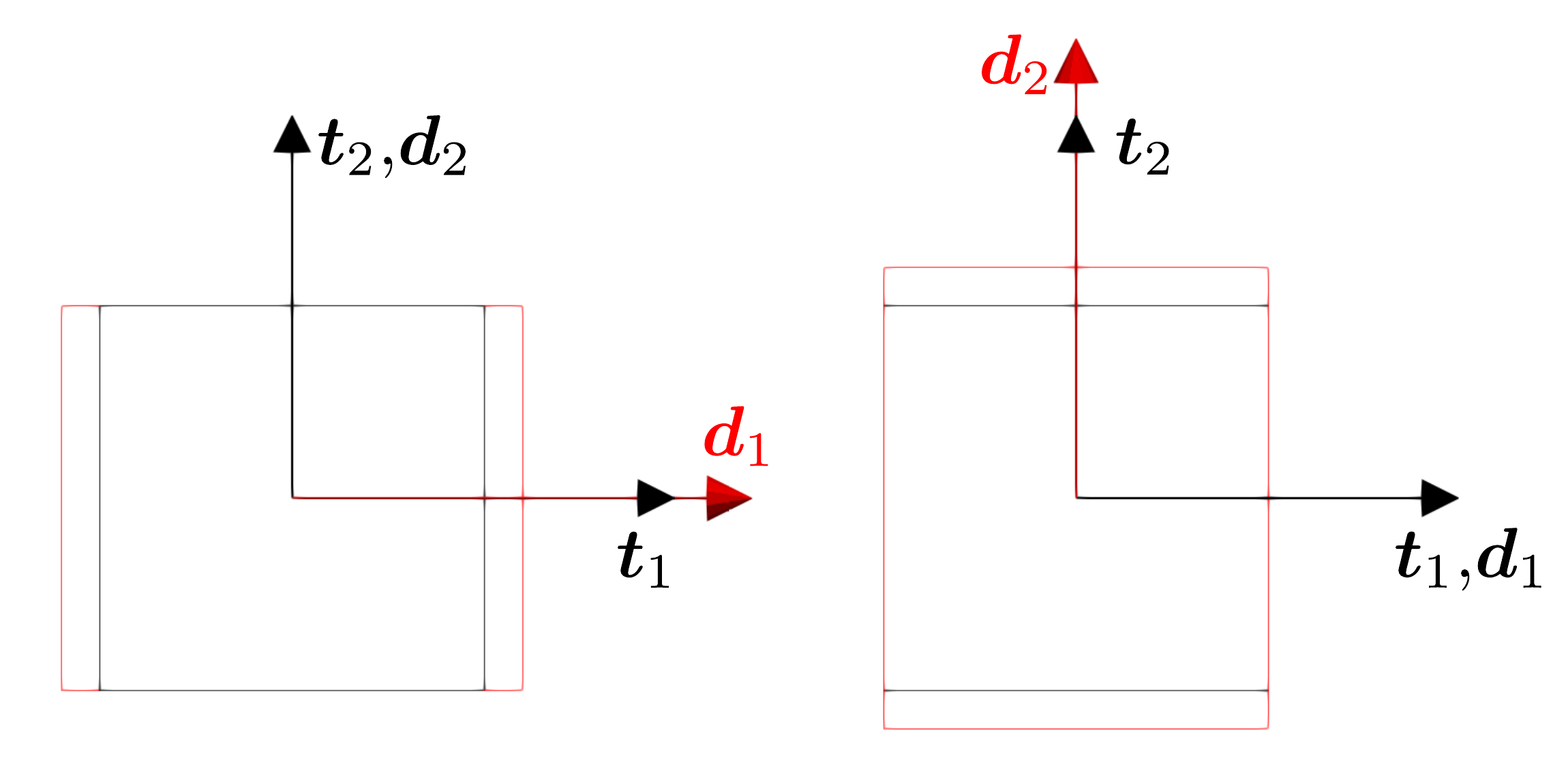}
		\caption{Two stretching modes}
		\label{cs_deform_mode_1}			
	\end{subfigure}        
	\quad
	\begin{subfigure}[b]{0.3125\textwidth}\centering
		\includegraphics[width=\linewidth]{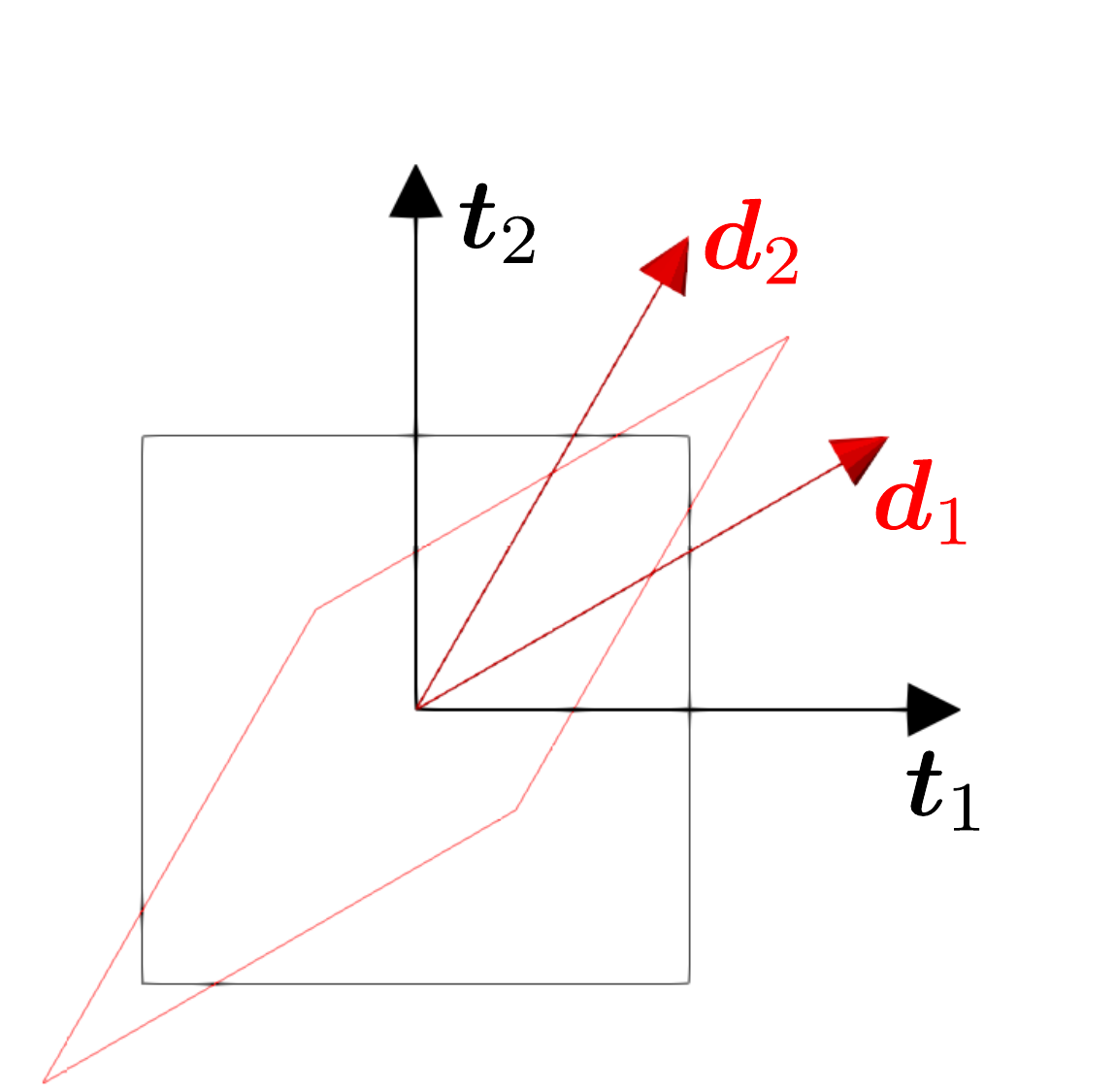}
		\caption{Pure shear mode}
		\label{cs_deform_mode_3}			
	\end{subfigure}	
	\caption{A schematic illustration of three kinematical (deformation) modes of the cross-section, not constrained by $\boldsymbol{\phi}=\boldsymbol{\psi}=\boldsymbol{0}$.}
	\label{beam_cs_strn_mode_allow}	
\end{figure}
\end{remark}
\begin{remark} 
\label{fe_disc} \small \textcolor{blue}{\textit{Objectivity and path-independence of the beam formulation.}
In the finite element discretization, the beam's configuration variable $\boldsymbol{y}\coloneqq\left[\boldsymbol{\varphi}^\mathrm{T},\boldsymbol{d}_1^{\,\mathrm{T}},\boldsymbol{d}_2^{\,\mathrm{T}}\right]^\mathrm{T}\in\Bbb{R}^d$ ($d=9$; see Appendix \ref{app_beam_form}) is approximated by \citep{choi2024objective}
\begin{align}
    \label{config_var_fe_disc}
    \boldsymbol{y}^h= \sum\limits_{I = 1}^{{n_e}} {{N^p_I}(\xi)\,{{\bf{y}}_I}},
\end{align}
where we have used NURBS basis functions $N^p_I$ of degree $p$, and $n_e=p+1$ denotes the number of local support basis functions in each element, and ${\bf{y}}_I\in\Bbb{R}^d$ denotes the control coefficients. It should be noted that the joint's kinematic variables are inherently discrete, where no further finite element approximation is applied, since the constraints are defined point-wise along the beam's axis. This implies that the objectivity and path-independence of the present finite element formulation, verified in \citet{choi2024objective}, are still valid, even after applying the joint constraints.}
\end{remark}
\noindent \textcolor{blue}{Applying the finite element discretization, we obtain the following system of linear equations:
\begin{subequations}
\label{reduced_eq_algeb_K_G_set}
    \begin{align}
        \label{reduced_eq_algeb_K_G}
    {{\bf{K}}_{\rm{r}}}\,\Delta {{{\bf{X}}_{\rm{f}}}}  + {{\bf{G}}_{\rm{r}} ^{\,\rm{T}}}\Delta {\bf{f}} = {{\bf{R}}_{\rm{r}}},
    \end{align}
    \begin{align}
        \label{constraint_linearized_after_row_red}
        {{{\bf{G}}_{\rm{r}}}}\,\Delta {{{\bf{X}}_{\rm{f}}}}  =  - {\bf{c}}.
    \end{align}
\end{subequations}
Here, $\bf{K}$ and $\bf{R}$ denote the global tangent stiffness matrix and residual vector, respectively, where $(\bullet)_\mathrm{r}$ denotes the reduced matrices and vectors due to the application of kinematic boundary conditions. Further, $\Delta{\bf X}_\mathrm{f}\in\Bbb{R}^{n_\mathrm{f}}$ denotes the increment of the free (reduced) kinematical unknown variables, and $n_\mathrm{f}$ denotes its number. Further, ${\bf f},{\bf{c}}\in\Bbb{R}^M$ denote the global arrays of the Lagrange multipliers and the constraints (${\Bbb f}$, ${\Bbb c}\in\Bbb{R}^m$) due to the finite element assembly, respectively, where $M$ denotes the total number of constraints; see also Remark \ref{rem_coupling_number_constraints}. 
\begin{remark}
\label{rem_rot_update}\small
Using the rotational increment vector $\Delta\boldsymbol{\bar\theta}$ in $\Delta{\bf X}_\mathrm{f}$, we update the joint's orientation $\boldsymbol{\bar R}\in\mathrm{SO(3)}$. For a singularity-free rotational update procedure with a minimal secondary storage, we utilize four quaternion parameters per joint, which refers to \citet{simo1986three}. 
\end{remark}
\noindent During the solution process of the system of algebraic equations (\ref{reduced_eq_algeb_K_G_set}), numerical instability may arise due to
\begin{itemize}
    \item Redundant constraints, i.e., deficiency of row rank in the constraint Jacobian ${\bf{G}}_\mathrm{r}$,
    \item Ill-conditioned system of algebraic equations due to scale difference between $\bf{K}_\mathrm{r}$ and ${\bf{G}}_\mathrm{r}$.
\end{itemize} 
In the following section, we present a discrete null space method to resolve these issues and further reduce the system size.}
\setcounter{remark}{0}
\section{Discrete null space method for stability and efficiency}
\label{null_space_sec}
In this section, we present a solution procedure for the system of algebraic equations in Eq.\,(\ref{reduced_eq_algeb_K_G_set}), with two additional treatments for improving numerical stability. First, we identify the redundant constraints and partition the unknown kinematic variables into independent/dependent ones. This process relies on a factorization of the (column-reduced) constraint Jacobian (${\bf G}_\mathrm{r}$), see \citet{wehage1982generalized}. Second, we further apply a null space method for improved conditioning as well as reduced size of the resulting system, which is based on the \textit{fundamental basis} of null space \citep{wolfe1962reduced}, see also \citet[Section 6]{benzi2005numerical} and \citet{rees2018comparative} for literature reviews. These two procedures are explained in detail below. For computational efficiency, we exploit the assumed locality of the constraints, that is, the number of constraints ($M$) is much smaller than that of free kinematic degrees-of-freedom ($n_\mathrm{f}$).
\subsection{Partitioning of the constraint Jacobian matrix}
By extracting the independent $r$ rows and columns, and rearranging, we have a partitioned constraint Jacobian
\begin{align}
    \label{row_partition_Gr}
    {{\overline{{\bf{G}}_{\rm{r}}}} } \coloneqq {\bf{P}}_\mathrm{row}\,{{\bf{G}}_{\rm{r}}}\,{\bf{P}}^\mathrm{T}_\mathrm{col} = \left[ {\renewcommand{\arraystretch}{1.5}\begin{array}{*{20}{l}}
\left[{{{\bf{G}}_{\rm{1}}}}\right]_{{r\times{n_\mathrm{f}}}}\\
\left[{{{\bf{G}}_{\rm{2}}}}\right]_{{(M-r)\times{n_\mathrm{f}}}}
\end{array}} \right],
\end{align}
where $r$ denotes the rank of the constraint Jacobian ${\bf{G}}_\mathrm{r}\in\Bbb{R}^{M\times{n_\mathrm{f}}}$. Here, ${\bf G}_1$ and ${\bf G}_2$ represent the \textit{independent} and \textit{dependent} rows, partitioned by the row-based permutation matrix, ${\bf P}_\mathrm{row}$. Further, those columns have also been rearranged by column-based permutation matrix ${\bf P}^\mathrm{T}_\mathrm{col}$, such that
\begin{align}
    {{\bf{G}}_{\rm{1}}} = \left[ {\begin{array}{*{20}{c}}
{{{\left[ {{{\bf{G}}_{\rm{u}}}} \right]}_{r \times r}}}&{{{\left[ {{{\bf{G}}_{\rm{v}}}} \right]}_{r \times \left( {{n_{\rm{f}}} - r} \right)}}}
\end{array}} \right],
\end{align}
where ${\bf{G}}_\mathrm{u}$ denotes the \textit{invertible} submatrix. This column-based permutation corresponds to the partitioning of the free kinematic variables $\Delta{{\bf X}_\mathrm{f}}$ into \textit{independent} and \textit{dependent} parts, denoted by ${\bf u}\in\Bbb{R}^{r}$ and ${\bf v}\in\Bbb{R}^{n_\mathrm{f}-r}$, respectively, as \citep{wehage1982generalized}
\begin{align}
    {\bf{P}}_\mathrm{col}\,\Delta{\bf X}_\mathrm{f} = \left\{ {\begin{array}{*{20}{c}}
{{\bf u}}\\
{{\bf v}}
\end{array}} \right\}.
\end{align} 
\begin{remark} \small For a given constraint Jacobian ${\bf G}_\mathrm{r}$, we find its rank $r$, and the permutation matrices ${\bf P}_\mathrm{row}$ and ${\bf P}_\mathrm{col}$ for the partitioning in Eq.\,(\ref{row_partition_Gr}) by a factorization, for which we have utilized the subroutine \textit{dgeqp3}, provided by the \citet{onemkl_ref_2025}. Here, we performed two column-based factorizations of ${\bf G}_\mathrm{r}^\mathrm{T}$ and ${\bf G}_\mathrm{r}$ for finding ${\bf P}_\mathrm{row}$ and ${\bf P}_\mathrm{col}$, respectively. Since $M\ll{n_\mathrm{f}}$, we may consider that the time and space complexities of this operation are both linear with respect to $n_\mathrm{f}$. Note that \textit{dgeqp3} requires a dense matrix as an input, which requires changing the format of ${\bf G}_\mathrm{r}$ into a dense matrix. The efficiency of this procedure can be further improved, since ${\bf G}_\mathrm{r}$ is sparse, due to the locality of constraints. This remains future work. 
\end{remark}
\noindent Since the matrices ${\bf P}_\mathrm{row}$ and ${\bf P}_\mathrm{col}$ are orthogonal, from Eq.\,(\ref{row_partition_Gr}) we have 
\begin{align}
    \label{Gr_bar_EGP}
    {\bf{G}}_\mathrm{r} = {{\bf P}}_\mathrm{row}^\mathrm{T}\,\overline {{{\bf{G}}_{\rm{r}}}}\, {{\bf{P}}}_\mathrm{col}.
\end{align}
Substituting Eq.\,(\ref{Gr_bar_EGP}) into Eqs.\,(\ref{reduced_eq_algeb_K_G}) and (\ref{constraint_linearized_after_row_red}), we have
\begin{subequations}
\label{reduced_eq_linearized_saddle}
\begin{align}
    \label{reduced_eq_algeb_K_G_subs_P}
{{\bf{K}}_{\rm{r}}}\Delta {{{\bf{X}}_{\rm{f}}}} + \left(\overline {{{\bf{G}}_{\rm{r}}}}\,{\bf{P}}_\mathrm{col}\right)^{\mathrm{T}} {{\bf P}}_\mathrm{row}\,\Delta {\bf{f}} = {{\bf{R}}_{\rm{r}}},
\end{align}
\begin{align}
    \label{constraint_linearized_after_row_red_subs_P}
    \left(\overline {{{\bf{G}}_{\rm{r}}}}\,{\bf{P}}_\mathrm{col}\right)\Delta {{{\bf{X}}_{\rm{f}}}}  =  - {\bf P}_\mathrm{row}\,{\bf{c}}.
\end{align}
\end{subequations}
In the following, we construct a solution space for $\Delta {\bf X}_\mathrm{f}$ by the null space $\mathcal{N}\!\left(\overline {{{\bf{G}}_{\rm{r}}}}\,{\bf{P}}_\mathrm{col}\right)$, so that the number of unknown variables is reduced by the rank $r$, and the constraints in Eq.\,(\ref{constraint_linearized_after_row_red_subs_P}) can be exactly satisfied. Here, $\mathcal{N}\!\left(\bullet\right)$ denotes the null space of a given matrix.
\subsection{Null space method}
We solve the system of $n_\mathrm{f}$ equations in Eq.\,(\ref{reduced_eq_algeb_K_G_subs_P}), subject to $r$ independent constraints in Eq.\,(\ref{constraint_linearized_after_row_red_subs_P}), using a null space method. Here, we parameterize the solution $\Delta{\bf{X}}_\mathrm{f}\in\Bbb{R}^{n_\mathrm{f}}$, as
\begin{align}
    \label{sol_h_nh}
    \Delta {{{\bf{X}}_{\rm{f}}}} = {\bf{Z}}{\bf{U}}_\mathrm{h}  + {\bf{U}}_\mathrm{p},
\end{align}
where ${\bf Z}\in\Bbb{R}^{{n_\mathrm{f}}\times(n_\mathrm{f}-r)}$ denotes the so-called \textit{null space matrix}, whose columns span $\mathcal{N}\!\left(\overline {{{\bf{G}}_{\rm{r}}}}\, {{\bf{P}}}_\mathrm{col}\right)$, i.e., 
\begin{align}
    \label{null_GPtZ_0}
    \left(\overline {{{\bf{G}}_{\rm{r}}}}\, {{\bf{P}}_\mathrm{col}}\right){\bf{Z}}=\boldsymbol{0}.
\end{align}
Here, the reduced unknown coefficient vector ${\bf{U}}_\mathrm{h}\in\Bbb{R}^{n_\mathrm{f}-r}$ represents the \textit{homogeneous} solution. Further, due to the residual of constraint functions at an inequilibrium, i.e., the right-hand side of Eq.\,(\ref{constraint_linearized_after_row_red_subs_P}), we also need 
a \textit{particular} solution, ${\bf{U}}_\mathrm{p}\in\Bbb{R}^{n_\mathrm{f}}$. Substituting Eq.\,(\ref{sol_h_nh}) into Eq.\,(\ref{constraint_linearized_after_row_red_subs_P}), and using Eq.\,(\ref{null_GPtZ_0}), we have
\begin{align}
    \label{part_sol_gov_eq}
    \left(\overline {{{\bf{G}}_{\rm{r}}}}\, {{\bf{P}}_\mathrm{col}}\right){\bf{U}}_\mathrm{p}  =  - {\bf{P}}_\mathrm{row}\,{\bf{c}}.    
\end{align}
In the following, we detail how to find ${\bf Z}$ and ${\bf{U}}_\mathrm{p}$, and the subsequent solution procedure.
\subsubsection*{Step 1: Find a null space matrix and a particular solution}
From using the fact that the matrix ${{\bf{G}}_\mathrm{u}}$ is invertible, and ${\bf P}_\mathrm{col}$ is orthogonal, we can simply find a null space matrix from 
\begin{align}
    \label{null_matrix_Z}
    {\bf{Z}} \coloneqq {\bf{P}}^\mathrm{T}_\mathrm{col}\underbrace{\left[ {\renewcommand{\arraystretch}{1.5}\begin{array}{*{20}{c}}
{ - {{{{{\bf{G}}^{\,-1}_{\rm{u}}}} }}\,{{{\bf{G}}_{\rm{v}}}} }\\
{{{\bf{1}}_{{n_{\rm{f}}} - r}}}
\end{array}} \right]}_{\eqqcolon{\bar{\bf Z}}},
\end{align}
due to ${\overline{{{\bf G}}_\mathrm{r}}}\,{\bar{\bf Z}}={\bf 0}$. This matrix $\bf Z$ is often called a \textit{fundamental basis} \citep{wolfe1962reduced}, see also \citet{wehage1982generalized}. In order to find a particular solution, we first rewrite Eq.\,(\ref{part_sol_gov_eq}), as
\begin{align}
    \label{part_sol_gov_eq_recall}
    \overline {{{\bf{G}}_{\rm{r}}}}\hspace{-3pt}\underbrace{{\bf P}_\mathrm{col}{\bf U}_\mathrm{p}}_{={\left[ {{\bf{u}}_{\rm{p}}^{\rm{T}},{\bf{v}}_{\rm{p}}^{\rm{T}}} \right]^{\rm{T}}}} = - \hspace{-5pt}\underbrace{{\bf{P}}_\mathrm{row}\,{\bf{c}}}_{=\left[{\bf c}^\mathrm{T}_1,{\bf c}^\mathrm{T}_2\right]^\mathrm{T}},    
\end{align}
where ${\bf u}_\mathrm{p}\in\Bbb{R}^r$ and ${\bf v}_\mathrm{p}\in\Bbb{R}^{n_\mathrm{f}-r}$ are the independent and dependent parts of the particular solution, respectively. Further, ${\bf c}_1\in\Bbb{R}^{r}$ and ${\bf c}_2\in\Bbb{R}^{M-r}$ denote the independent and dependent constraints, respectively. By simply choosing ${\bf v}_\mathrm{p}=\boldsymbol{0}$, we can obtain a particular solution, as
\begin{align}
    \label{particular_sol_up_null}
    {{\bf{U}}_{\rm{p}}} =  {\bf{P}}^\mathrm{T}_\mathrm{col}\left\{ {\renewcommand{\arraystretch}{1.5}\begin{array}{*{20}{c}}
-{{\bf{G}}_{\rm{u}}^{\,-1}}{{\bf c}}_1\\
{\bf{0}}
\end{array}} \right\}.
\end{align}
\subsubsection*{Step 2: Solve the reduced system of algebraic equations}
Further, substituting the null space matrix from Eq.\,(\ref{null_matrix_Z}) into Eq.\,(\ref{reduced_eq_algeb_K_G_subs_P}), we have the reduced problem: For a given ${\bf{U}}_\mathrm{p}\in\Bbb{R}^{n_\mathrm{f}}$, find ${\bf U}_\mathrm{h}\in\Bbb{R}^{n_\mathrm{f}-r}$ such that
\begin{align}
    \label{sys_alg_eq_reduced_null}
    \underbrace{\left( {{{\bf{Z}}^{\rm{T}}}{{\bf{K}}_{\rm{r}}}{\bf{Z}}} \right)}_{\eqqcolon{\overline{{\bf K}_\mathrm{r}}}}{{\bf{U}}_{\rm{h}}} = \underbrace{{{\bf{Z}}^{\rm{T}}}\!\left( {{{\bf{R}}_{\rm{r}}} - {{\bf{K}}_{\rm{r}}}{{\bf{U}}_{\rm{p}}}} \right)}_{\eqqcolon{\overline{{\bf R}_\mathrm{r}}}},
\end{align}
where we have the reduced system matrix ${\overline{{\bf K}_\mathrm{r}}}$ and residual vector $\overline{{\bf R}_\mathrm{r}}$. Note that the reduced matrix ${\overline{{\bf K}_\mathrm{r}}}$ is positive definite, but may be unsymmetric at inequilibrium, see Remark \ref{rem_sym_tstiff_geom}. 
\begin{remark} \small \textit{Sparsity due to the fundamental basis}. Even though the first $r$ rows of $\bar{\bf Z}$ may be dense, the remaining ${n_\mathrm{f}}-r$ rows are still sparse. As we assume $r\ll{n_\mathrm{f}}$, the matrix $\bf Z$ is sparse, which enables us to keep the final reduced system matrix $\overline{{\bf K}_\mathrm{r}}$ sparse.
\end{remark}
\noindent Then, from using Eq.\,(\ref{sol_h_nh}), we can obtain $\Delta {\bf X}_\mathrm{f}$, using ${\bf U}_\mathrm{h}$ from Eq.\,(\ref{sys_alg_eq_reduced_null}), and the particular solution from Eq.\,(\ref{particular_sol_up_null}).
\subsubsection*{Step 3: Calculate the Lagrange multipliers}
For the given solution $\Delta{\bf X}_\mathrm{f}$, we still need to calculate the increment of the Lagrange multipliers, $\Delta{\bf f}$. Therefore, we first rewrite Eq.\,(\ref{reduced_eq_algeb_K_G_subs_P}), using the orthogonality of ${\bf P}_\mathrm{col}$, as
\begin{align}
    \label{reduced_eq_algeb_K_G_recover_mult}
\overline {{{\bf{G}}_{\rm{r}}}}^{\,\mathrm{T}}\hspace{-9pt}\underbrace{{{\bf{P}}_\mathrm{row}}\,\Delta {{\bf f}}}_{={\left[ {\Delta {\bf{f}}_1^{\,\rm{T}},\Delta {\bf{f}}_2^{\,\rm{T}}} \right]^{\rm{T}}}} = \underbrace{{\bf P}_\mathrm{col}\left({{\bf{R}}_{\rm{r}}} - {{\bf{K}}_{\rm{r}}}\Delta {{{\bf{X}}_{\rm{f}}}}\right)}_{=\left[{\bf r}_\mathrm{u}^\mathrm{T},{\bf r}_\mathrm{v}^\mathrm{T}\right]^\mathrm{T}},
\end{align}
where $\Delta{\bf f}_1\in\Bbb{R}^{r}$ and $\Delta{\bf f}_2\in\Bbb{R}^{M-r}$ denote the independent and dependent multipliers, respectively, partitioned by the operator ${\bf P}_\mathrm{row}$. Further, the right-hand side can be also partitioned, due to the operation ${\bf P}_\mathrm{col}$, into the independent and dependent parts, denoted by ${\bf r}_\mathrm{u}\in\Bbb{R}^{r}$ and ${\bf r}_\mathrm{v}\in\Bbb{R}^{n_\mathrm{f}-r}$, respectively. For the redundant constraints, we may simply choose $\Delta{\bf f}_2=\boldsymbol{0}$, from which we have 
\begin{align}
    \label{null_update_mult_inc}
    \Delta {\bf{f}} = {{{\bf P}}_\mathrm{row}^{\rm{T}}}\left\{ {\renewcommand{\arraystretch}{1.5}\begin{array}{*{20}{c}}
{ {\bf{G}}_{\rm{u}}^{ - {\rm{T}}}{{\bf{r}}_{\rm{u}}}}\\
{\bf{0}}
\end{array}} \right\}.
\end{align}
The above solution procedure is summarized in Algorithm\,\ref{algorithm_null_space}.
\begin{algorithm}[H]
\caption{Solution process using the null space method}\label{alg:cap}
\label{algorithm_null_space}
\begin{algorithmic}[1]
\State Find the operators ${\bf P}_\mathrm{row}$ and ${\bf P}_\mathrm{col}$ by factorizing ${\bf{G}}_{\rm{r}}^{\,\mathrm{T}}$ and ${{{{\bf{G}}_{\rm{r}}}}}$, respectively, using \textit{dgeqp3}.
\State $\mathcal{F}\left({\bf G}_\mathrm{u}\right)$ defines the factorized invertible matrix ${\bf G}_\mathrm{u}$. 
\State Calculate the null space matrix ${\bf Z}$ using Eq.\,(\ref{null_matrix_Z}), with $\mathcal{F}\left({\bf G}_\mathrm{u}\right)$.
\State Calculate the particular solution ${\bf U}_\mathrm{p}$ using Eq.\,(\ref{particular_sol_up_null}), with $\mathcal{F}\left({\bf G}_\mathrm{u}\right)$.
\State Solve the reduced system of equations in Eq.\,(\ref{sys_alg_eq_reduced_null}) to obtain ${\bf U}_\mathrm{h}$, and calculate the solution increment $\Delta{\bf X}_\mathrm{f}$ using Eq.\,(\ref{sol_h_nh}),
\State Calculate $\Delta\bf{f}$, using Eq.\,(\ref{null_update_mult_inc}), with $\mathcal{F}\left({\bf G}_\mathrm{u}\right)$.
\State Update the kinematic variables ${\bf X}_\mathrm{f}$ and the Lagrange multipliers ${\bf f}$ by $\Delta{\bf X}_\mathrm{f}$ and $\Delta{\bf{f}}$, respectively. For the rotational update, see Remark \ref{rem_rot_update}. 
\end{algorithmic}
\end{algorithm}
\setcounter{remark}{0}
\section{Numerical examples}
\label{num_ex}
We present three numerical examples having the following objectives.
\begin{itemize}
    \item \textbf{Ex.\,1.\,Torsion of a straight beam with Z-section.} We verify the accuracy and efficiency of the present beam formulation for large torsional deformation of an open cross-section. It is noted that we utilize the present constraint formulation to enforce the prescribed rotation and the clamped boundary condition. 
    \item \textbf{Ex.\,2.\,Right-angled frame with a prismatic joint.} We show the modeling of a prismatic joint, using the present variable-offset formulation. This incorporates the rigid joint as a special case by eliminating the offset degrees-of-freedom. We also verify the accuracy and efficiency of the beam solutions by comparison with the brick solutions.
    \item \textbf{Ex.\,3.\,Framed shallow dome.} We verify the present beam formulation for modeling a built-up structure in an instability problem. \textcolor{blue}{Here, we also investigate the effect of a correct representation of joint stiffness by the present offset-joint formulation.}
\end{itemize}
\textcolor{blue}{
\begin{remark}
\label{rem_partial_clamped} \small
The present joint formulation constrains only purely rigid motion, while allowing cross-sectional warping. For beam solutions, our investigation is limited to this \textit{partially clamped} joint and boundary condition. For a fully clamped condition, additional interface constraints on the warping strains are necessary. This remains future work. On the other hand, for bricks, we can easily constrain the warping by adding volumetric patches having much higher stiffness. This \textit{fully clamped} condition for bricks has also been implemented to verify the beam solution. 
\end{remark}
\noindent In all examples, we have used uniform load increments, and $n_\mathrm{load}$ denotes the total number of load steps. In both beam and brick formulations, we use NURBS (Non-Uniform Rational B-Spline) basis functions in the framework of IGA, where we use the following notations:
\begin{itemize}
    \item For the beam formulation, $p$ and $n_\mathrm{el}$ denote the degree of basis functions and the number of elements along the axis, respectively. For the enhanced cross-sectional strain field, we use $q_\mathrm{a}$, $n^\mathrm{a}_\mathrm{el}={n^\mathrm{W}_\mathrm{el}}\times{n^\mathrm{H}_\mathrm{el}}\,(\times{N})$, and $m^\mathrm{a}_\mathrm{cp}$ for the degree of basis functions, the number of elements, and the total number of control points, respectively, where $N$ denotes the number of surface patches in the cross-section. Further, we use the degree of basis functions $p_\mathrm{a}=p-1$ for the approximation of enhanced strain parameters along the axis. The relevant description of the finite element discretization of the enhanced strain field can be found in Appendix \ref{app_warp_cor_local_shape}.
    \item For the brick formulation, we use $\mathrm{deg.}=(p_\mathrm{L},p_\mathrm{W},p_\mathrm{H})$ and $n_\mathrm{el}={n^\mathrm{L}_\mathrm{el}}\times{n^\mathrm{W}_\mathrm{el}}\times{n^\mathrm{H}_\mathrm{el}}\,(\times N)$ for denoting the degrees of basis functions and the number of elements, respectively, and $N$ denotes the number of patches in the cross-section.
\end{itemize}
Here, `L', `W', and `H' represent the directions along the beam's length, width, and height, respectively.
\begin{remark} \small In the present numerical examples, the constraint Jacobian matrices always have full row rank. That is, all the constraints are independent of each other. This is because the connected beams at a joint have independent constrained degrees-of-freedom.
\end{remark}
}
\subsection{Torsion of a straight beam with Z-section}
\begin{figure}[H]
	\centering
	\begin{subfigure}[b]{0.75\textwidth}\centering
		\includegraphics[width=\linewidth]{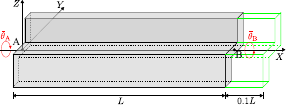}
	\end{subfigure}
	\caption{Torsion of a Z-section beam: Initial geometry and boundary conditions. \textcolor{blue}{The green-colored part represents the added volume for constraining warping in the fully clamped brick model.}}
	\label{z_shape_init_config_bdc}
\end{figure}

\begin{figure}[H]
	\centering
	\begin{subfigure}[b]{0.3\textwidth}\centering
	\includegraphics[width=\linewidth]{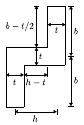}
		\caption{Case 1 (no fillet)}		
		\label{z_shape_without_fillet}			
	\end{subfigure}
    \quad\quad
	\begin{subfigure}[b]{0.3\textwidth}\centering
	\includegraphics[width=\linewidth]{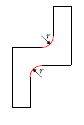}
		\caption{Case 2 (added fillet)}		
		\label{z_shape_with_fillet}			
	\end{subfigure}
	\caption{Torsion of a Z-section beam: Two different cases for the cross-section. In (b), we have added a circular fillet of radius $r$ at the two inner sharp corners.}
	\label{Z_shape_twist_init_config_bdc_cs_dim}
\end{figure}

\begin{figure}[H]
	\centering
	\begin{subfigure}[b]{0.3125\textwidth}\centering
		\includegraphics[width=\linewidth]{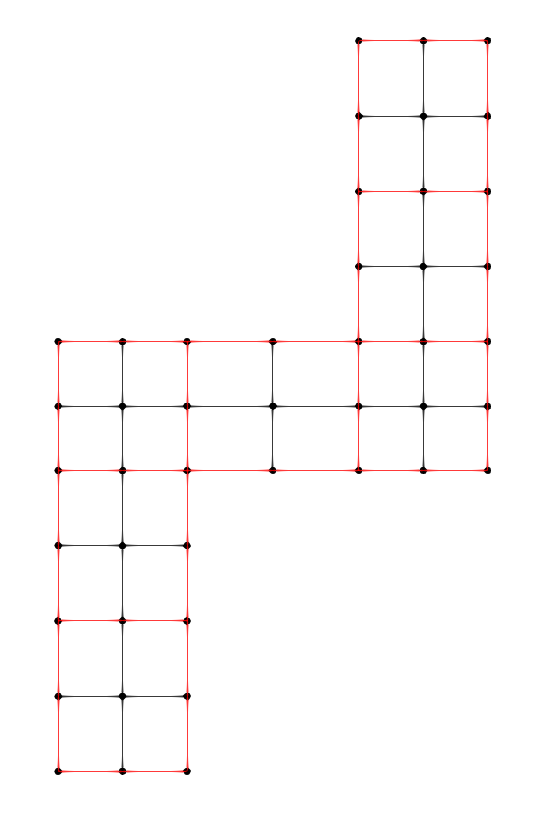}
		\caption{Case 1 (no fillet)}		
		\label{shape_twist_beam_cs_mesh_2x2}	
	\end{subfigure}
	\begin{subfigure}[b]{0.3125\textwidth}\centering
		\includegraphics[width=\linewidth]{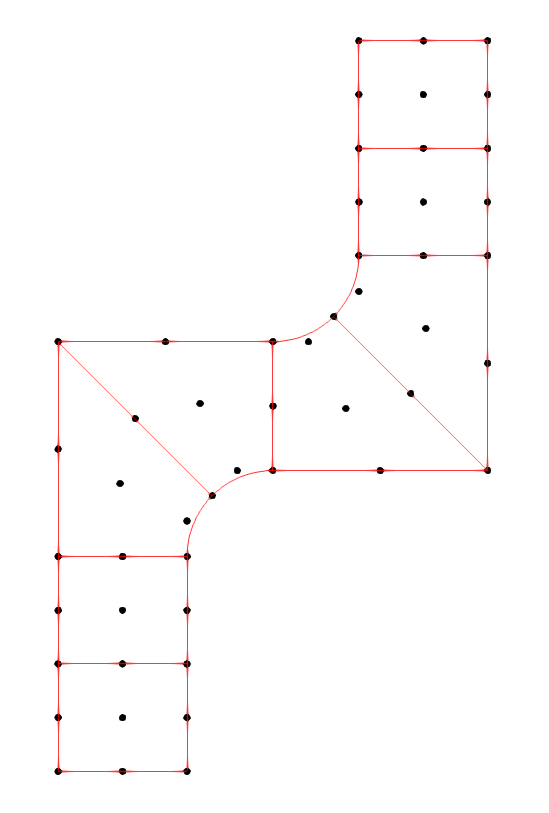}
		\caption{Case 2 (added fillet)}		
		\label{shape_twist_beam_cs_mesh_3x3}			
	\end{subfigure}
	\caption{Torsion of a Z-section beam: Cross-sectional discretization for the enhanced strain field. Here, the red-colored lines represent the boundaries of the surface patches. \textcolor{blue}{For each case, we have used (a) linear ($q_a=1$, $n^\mathrm{a}_\mathrm{el}=2\times2(\times7)$, $m^\mathrm{a}_\mathrm{cp}=45$), and (b) quadratic NURBS elements ($q_a=2$, $n^\mathrm{a}_\mathrm{el}=1\times1(\times8)$, $m^\mathrm{a}_\mathrm{cp}=51$). Black dots represent the control points. $q_\mathrm{a}$, $n^\mathrm{a}_\mathrm{el}$, and $m^\mathrm{a}_\mathrm{cp}$ denote the degree of basis functions, the number of elements, and the number of control points; see Remark \ref{app_remark_eas_notation}.}}
	\label{Z_shape_twist_beam_cs_mesh}
\end{figure}

\begin{figure}[H]
	\centering
	\begin{subfigure}[b]{0.4875\textwidth}\centering
			\includegraphics[width=\linewidth]{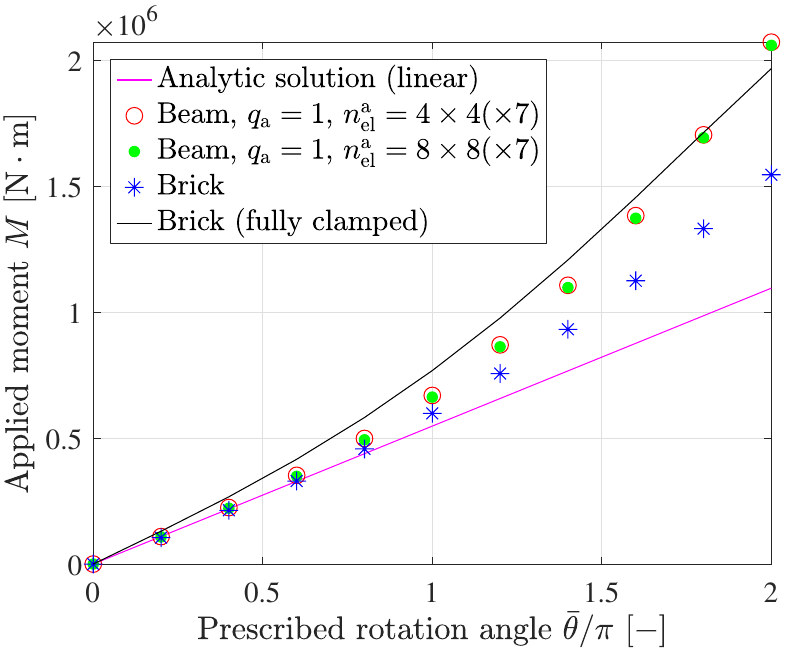}
			\caption{Case 1 (no fillet)}		
			\label{Z_torsion_mnt_case_1}			
	\end{subfigure}
	\begin{subfigure}[b]{0.4875\textwidth}\centering
			\includegraphics[width=\linewidth]{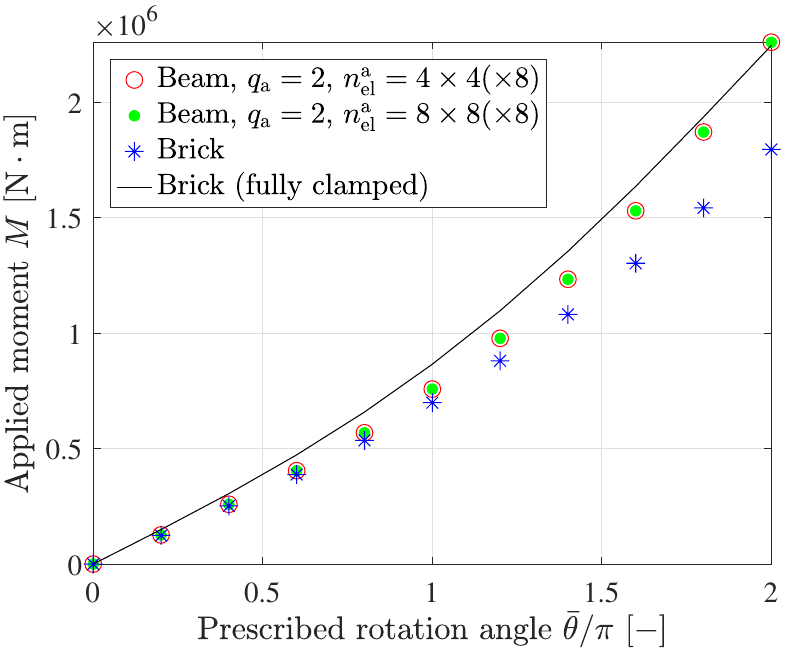}
			\caption{Case 2 (added fillet)}		
			\label{Z_torsion_mnt_case_2}			
	\end{subfigure}    
	\caption{Torsion of a Z-section beam: Change of the applied moment along the prescribed rotation. Here, the beam solutions are obtained from using $p=3$, $n_\mathrm{el}=10$.}
	\label{Z_shape_twist_comparison_equilibrium_path}
\end{figure}

\begin{figure}[H]
	\centering
	\begin{subfigure}[b]{0.4875\textwidth}\centering
		\includegraphics[width=\linewidth]{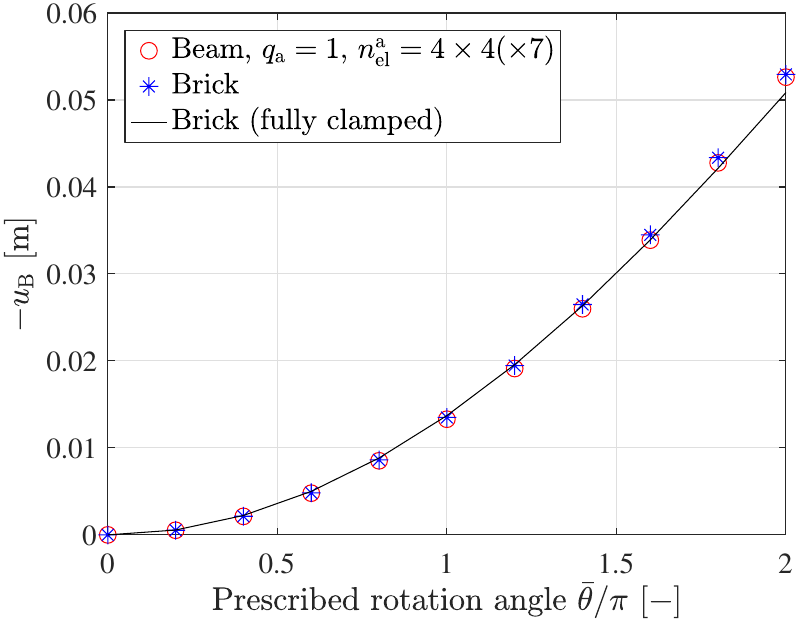}
		\caption{Case 1 (no fillet)}	
		\label{axial_disp}			
	\end{subfigure}	    
	\begin{subfigure}[b]{0.4875\textwidth}\centering
		\includegraphics[width=\linewidth]{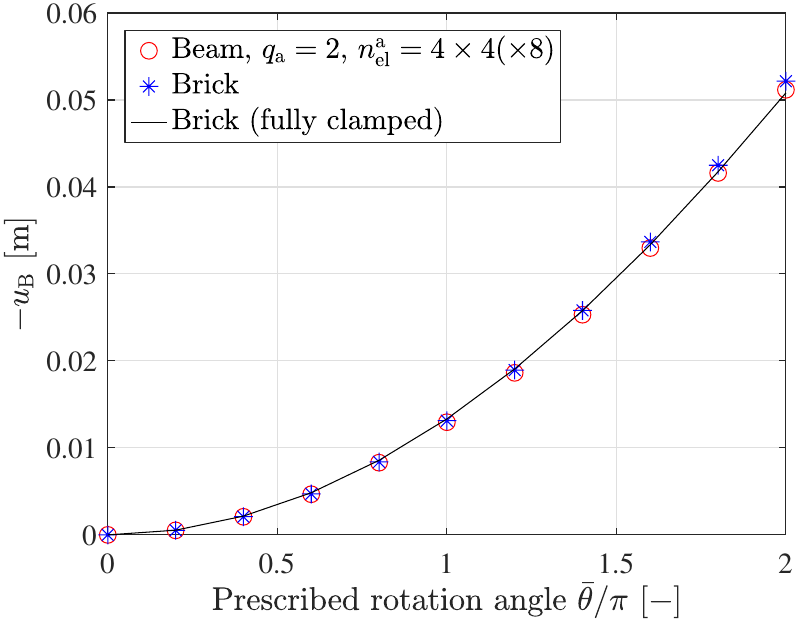}
		\caption{Case 2 (added fillet)}	
		\label{axial_disp_case2}			
	\end{subfigure}
	\caption{Torsion of a Z-section beam: $X$-directional tip displacement at B, denoted by $u_\mathrm{B}$. Here, the beam solutions are obtained from using $p=3$, $n_\mathrm{el}=10$. }
	\label{Z_shape_twist_comparison_ax_disp}
\end{figure}

\begin{figure}[H]
	\centering
	\begin{subfigure}[b]{0.95\textwidth}\centering
		\includegraphics[width=\linewidth]{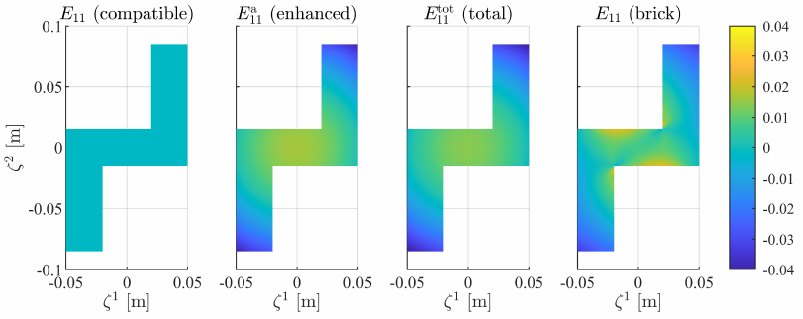}
		\label{Z_shape_twist_comparison_strain_field_E11}			
	\end{subfigure}
	\caption{Torsion of a Z-section beam: Comparison of the transverse normal strain component, $E_{11}$.}
	\label{Z_shape_twist_comparison_strain_field_diff_E11}
\end{figure}

\begin{figure}[H]
	\centering
	\begin{subfigure}[b]{0.95\textwidth}\centering
		\includegraphics[width=\linewidth]{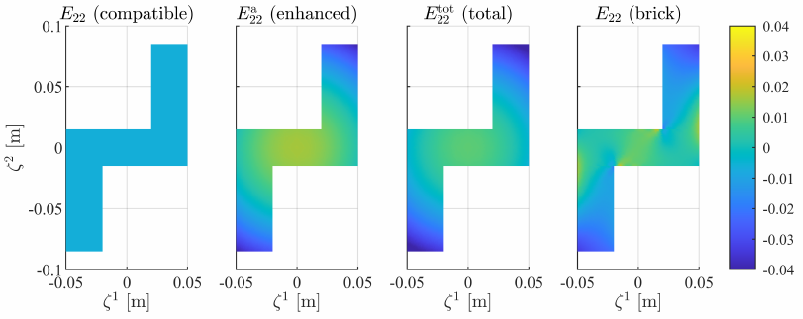}
		\label{Z_shape_twist_comparison_strain_field_E22}			
	\end{subfigure}
	\caption{Torsion of a Z-section beam: Comparison of the transverse normal strain component, $E_{22}$.}
	\label{Z_shape_twist_comparison_strain_field_diff_E22}
\end{figure}

\begin{figure}[H]
	\centering
	\begin{subfigure}[b]{0.95\textwidth}\centering
		\includegraphics[width=\linewidth]{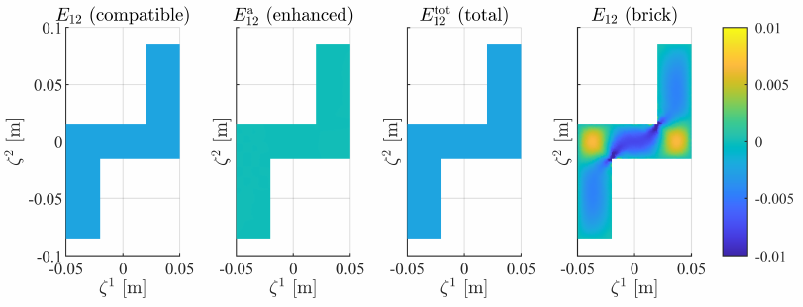}
		\label{Z_shape_twist_comparison_strain_field_E12}			
	\end{subfigure}
	\caption{Torsion of a Z-section beam: Comparison of the in-plane shear strain component, $E_{12}$.}
	\label{Z_shape_twist_comparison_strain_field_diff_E12}
\end{figure}

\begin{figure}[H]
	\centering
	\begin{subfigure}[b]{0.65\textwidth}\centering
		\includegraphics[width=\linewidth]{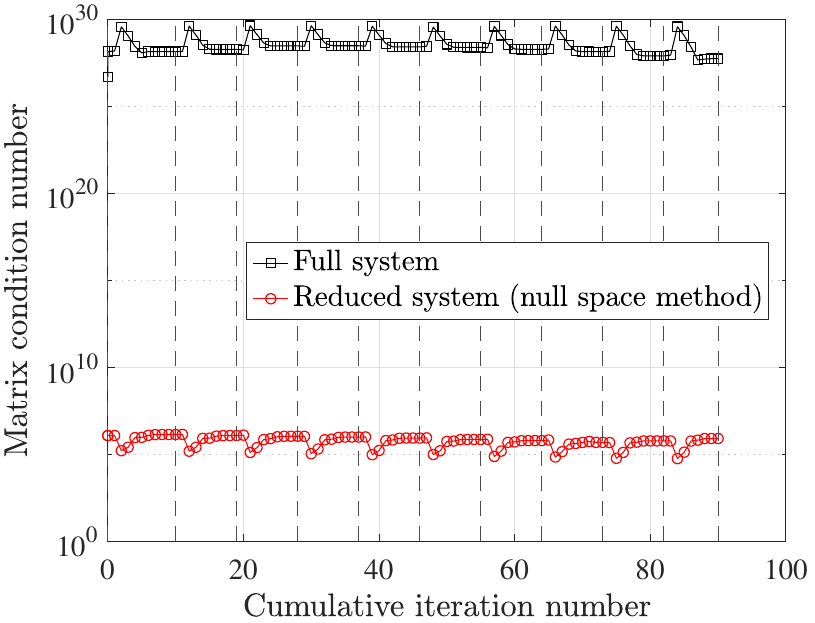}
		\label{Z_shape_twist_comparison_condition number}			
	\end{subfigure}
	\caption{Torsion of a Z-section beam: Comparison of the condition number of the global tangent stiffness matrix in Case 1 (no fillet). Here, we have used $p=3$, $p_\mathrm{a}=2$, and $n_\mathrm{el}=10$, $n^\mathrm{a}_\mathrm{el}=4\times4(\times7)$. The estimated condition number is calculated using the subroutine \textit{dgecon} from the \citet{onemkl_ref_2025}. \textcolor{blue}{Vertical dashed lines define load increments.}}
	\label{Z_shape_twist_comparison_condition_number}
\end{figure}

We consider a straight beam with Z-shaped cross-section under torsion, see Fig.\,\ref{z_shape_init_config_bdc} for an illustration. The twisting motion is prescribed by a rotation of angle $\bar \theta_\mathrm{A}=0$ and $\bar \theta_\mathrm{B}=2\pi$ about the $X$-axis for the cross-sections at the end points of the axis, denoted by A and B. The displacement at A is also constrained to avoid rigid-body translations. Here, $M$ denotes the applied moment corresponding to the prescribed rotation $\bar\theta_\mathrm{B}$. The beam has length $L=1\,\mathrm{m}$, and we consider the following two cases for the cross-section. 
\begin{itemize}
    \item Case 1: Dimensions are given by the three parameters $b=9\,\mathrm{cm}$, $t=3\,\mathrm{cm}$, and $h=7\,\mathrm{cm}$, see Fig.\,\ref{z_shape_without_fillet}. For this geometry, we refer to \citet[Section 5.1]{wackerfuss2011nonlinear}. 
    \item Case 2: Two fillets of radius $r$ are added to the inner corners of Case 1, see Fig.\,\ref{z_shape_with_fillet}.
\end{itemize}
In these two cases, we compare our beam solution to brick solutions. For the material model, we consider a St.\,Venant-Kirchhoff type hyperelasticity, with Young's modulus $E=210\,\mathrm{GPa}$ and Poisson's ratio $\nu=0.3$. \textcolor{blue}{As mentioned in Remark\,\ref{rem_partial_clamped}, for a reference brick solution using fully-clamped boundary conditions, we have introduced the following additional treatments:
\begin{itemize}
    \item We fix all degrees-of-freedom in the cross-section at the end point A.
    \item To constrain warping at the loaded cross-section at point B, we have extended the beam by $0.1L$ and increased Young's modulus by $100$ there; see Fig.\,\ref{z_shape_init_config_bdc} for an illustration.
\end{itemize}
}
\noindent In Fig.\,\ref{Z_shape_twist_beam_cs_mesh}, we show examples of the cross-sectional discretization for Cases 1 and 2, using linear and quadratic NURBS elements, respectively. 
\noindent \textcolor{blue}{In Tables \ref{tab_Z_shape_comp_ndof_case1} and \ref{tab_Z_shape_comp_ndof_case2}, we compare the degrees-of-freedom (DOFs) between beam and brick formulations. In the beam model, the number of global DOFs is much smaller than in the brick model, since the enhanced strain parameters are statically condensed and treated as internal variables. In Case 2, we have more internal and global DOFs for the beam and brick models than in Case 1 because we include an additional patch in the cross-section with the same mesh density. Further, the increased DOFs in the fully clamped brick models are due to the added volume patches. A convergence test for the brick solution is provided in Appendix \ref{app_z_shape_conv_test_brick}. In the following, we use brick solutions with $\mathrm{deg.}=(3,3,3)$ and $n_\mathrm{el}=50\times5\times5(\times{N})$ as reference solutions, where $N=7$ and 8 for Cases 1 and 2, respectively.}
\begin{table}[H]
  \centering
  \footnotesize
  \caption{Torsion of a Z-section beam (Case 1): Comparison of the degrees-of-freedom. Here, the element counts are defined per patch.}
    \begin{tabular}{lcccccccc}
    \toprule
    \multicolumn{1}{l}{} & \multicolumn{2}{c}{Degrees} &       & \multicolumn{2}{c}{Elements} &       & \multicolumn{2}{c}{DOFs} \\
\cmidrule{2-3}\cmidrule{5-6}\cmidrule{8-9}    \multicolumn{1}{c}{} & L     & \multicolumn{1}{c}{W,H} &       & L     & \multicolumn{1}{c}{W,H} &       & Global  & Internal \\
    \midrule
    \multicolumn{1}{l}{Beam} & 3     & 1     &       & 10    & 4     &       & 108   & 17220 \\
    \multicolumn{1}{l}{Brick} & 3     & 3     &       & 50    & 5     &       & 63591 & $-$ \\
    Brick (fully clamped) & 3     & 3     &       & 50    & 5     &       & 65997 & $-$ \\
    \bottomrule
    \end{tabular}%
  \label{tab_Z_shape_comp_ndof_case1}%
\end{table}%
\begin{table}[H]
  \centering
  \footnotesize
  \caption{Torsion of a Z-section beam (Case 2): Comparison of the degrees-of-freedom. Here, the element counts are defined per patch.}
    \begin{tabular}{lcccccccc}
    \toprule
    \multicolumn{1}{l}{} & \multicolumn{2}{c}{Degrees} &       & \multicolumn{2}{c}{Elements} &       & \multicolumn{2}{c}{DOFs} \\
\cmidrule{2-3}\cmidrule{5-6}\cmidrule{8-9}    \multicolumn{1}{c}{} & L     & \multicolumn{1}{c}{W,H} &       & L     & \multicolumn{1}{c}{W,H} &       & Global  & Internal \\
    \midrule
    \multicolumn{1}{l}{Beam} & 3     & 2     &       & 10    & 4     &       & 108   & 29340 \\
    \multicolumn{1}{l}{Brick} & 3     & 3     &       & 50    & 5     &       & 72495 & $-$ \\
    Brick (fully clamped) & 3     & 3     &       & 50    & 5     &       & 75237 & $-$ \\
    \bottomrule
    \end{tabular}%
  \label{tab_Z_shape_comp_ndof_case2}%
\end{table}%
\noindent \textcolor{blue}{In Figs.\,\ref{Z_torsion_mnt_case_1} and \ref{Z_torsion_mnt_case_2}, we compare the equilibrium paths from the beam and brick solutions for Cases 1 and 2, respectively. It is seen in Case 2 that, due to the fillet, the overall torsional stiffness is higher than in Case 1. The analytical solution in Fig.\,\ref{Z_torsion_mnt_case_1} assumes pure torsional motion and linear elasticity and is given by Eq.\,(\ref{app_z_shape_correct_K_asol}). This linear solution accounts for torsional stiffness very well in the small-strain range; however, it cannot capture the stiffening effect in the large-strain range. The fully clamped brick model exhibits higher torsional stiffness (black solid curves) than the others. By releasing the cross-sectional warping at the end cross-sections (i.e., under partially clamped conditions), the stiffness decreases significantly (blue stars). In the moderately large-strain range, the beam solution (red circles) agrees very well with the brick solution (blue stars). As the number of elements in the cross-section increases, the beam solution (green circles) approaches the brick solution under partially clamped conditions (blue stars), which is shown more clearly in Case 1 than in Case 2. However, in the very large-strain range, the beam solution exhibits more pronounced stiffening compared to the brick solutions.} 
\begin{observation} 
\label{observ_ex_z_torsion}\small
\textcolor{blue}{The discrepancy between the beam and brick solutions under the same partially clamped condition may arise due to the following reasons: First, for the application of the beam kinematics to the joint formulation in Appendix\,\ref{app_beam_form}, we incorporate up to first-order kinematic variables only, with the enhanced strains not considered. This is clearly different from the brick formulation, in which all degrees-of-freedom in the cross-section are involved. Second, the present EAS formulation for beams may not sufficiently account for the coupling between strain components. Further investigation on the discrepancy between the beam and brick solutions in the large-strain range remains future work.}
\end{observation}
\noindent \textcolor{blue}{In Fig.\,\ref{Z_shape_twist_comparison_ax_disp}, we compare the axial displacement at the tip B. The beam solution agrees very well with the brick solution. It is seen that, under the partially clamped condition, the beam solution slightly underestimates the axial displacement, which is consistent with the larger torsional stiffness of the beam solution shown in Fig.\,\ref{Z_shape_twist_comparison_equilibrium_path}.}
\noindent \textcolor{blue}{Further, at the final deformed configuration, we compare the cross-sectional strain fields (six components of the Green-Lagrange strain tensor) between beam and brick solutions at the middle of the axis, $s=L/2$, under partially clamped conditions in Case 1. Here, for the discretization of the beam's enhanced strain field, we use $q_a=1$ and $n^\mathrm{a}_\mathrm{el}=4\times4(\times7)$. First, the quadratic distribution of the longitudinal strain ($E_{33}$) can be captured by the compatible strains from the kinematics; see Appendix\,\ref{app_construct_warp_basis} for the relevant comments and references.}
\begin{figure}[H]
	\centering
	\begin{subfigure}[b]{0.95\textwidth}\centering
		\includegraphics[width=\linewidth]{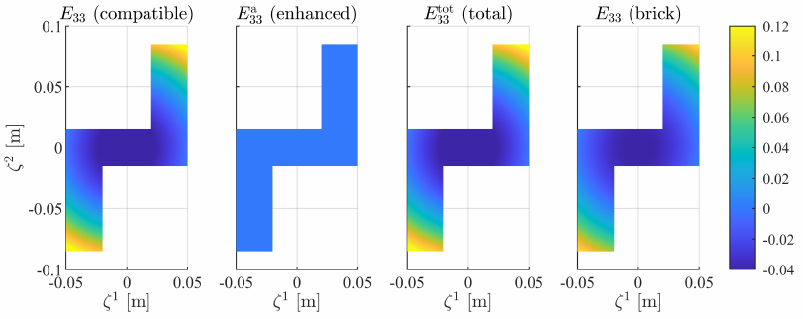}
		\label{Z_shape_twist_comparison_strain_field_E33}
	\end{subfigure}
	\caption{Torsion of a Z-section beam: Comparison of the longitudinal strain component, $E_{33}$.}
	\label{Z_shape_twist_comparison_strain_field_diff_E33}
\end{figure}
\begin{figure}[H]
	\centering
	\begin{subfigure}[b]{0.95\textwidth}\centering
		\includegraphics[width=\linewidth]{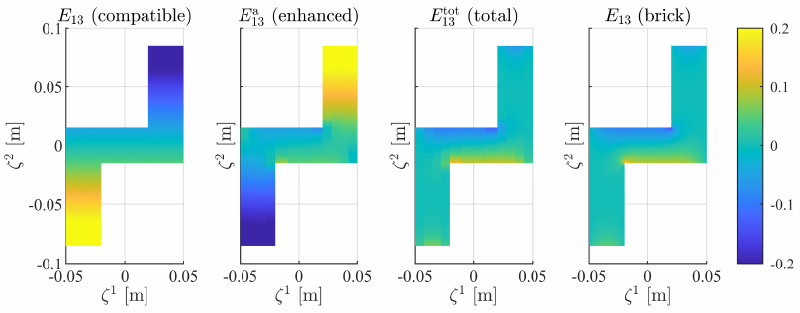}
		\label{Z_shape_twist_comparison_strain_field_E13}			
	\end{subfigure}
	\caption{Torsion of a Z-section beam: Comparison of the transverse shear strain component, $E_{13}$. The total strain in the beam formulation agrees very well with that of the brick formulation.}
	\label{Z_shape_twist_comparison_strain_E13}
\end{figure}
\noindent It is remarkable that the beam solution shows excellent agreement with the brick solutions for the transverse shear strain components $E_{13}$ and $E_{23}$. 
\begin{figure}[H]
	\centering
	\begin{subfigure}[b]{0.95\textwidth}\centering
		\includegraphics[width=\linewidth]{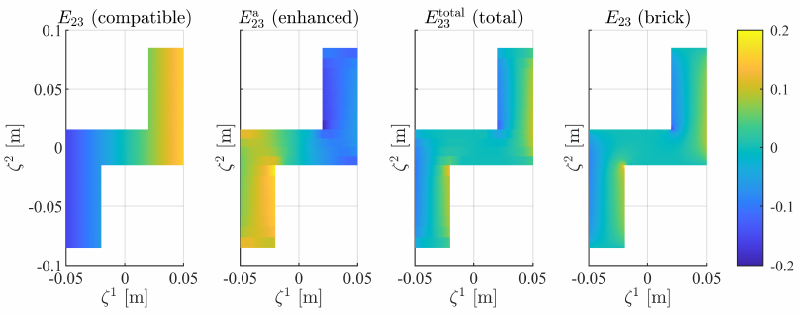}
		\label{Z_shape_twist_comparison_strain_field_E23}			
	\end{subfigure}
	\caption{Torsion of a Z-section beam: Comparison of the transverse shear strain component, $E_{23}$. The total strain in the beam formulation agrees very well with that of the brick formulation.}
	\label{Z-Z_shape_twist_comparison_strain_E23}
\end{figure}
\noindent\textcolor{blue}{In Figs.\,\ref{Z_shape_twist_comparison_strain_field_diff_E11} and \ref{Z_shape_twist_comparison_strain_field_diff_E22}, it is seen that the enhancement of the transverse normal strains ($E_{11}$, $E_{22}$) are limited to the quadratic field coupled with the longitudinal strain in Fig.\,\ref{Z_shape_twist_comparison_strain_field_diff_E33}, due to the nonzero Poisson's ratio. The enhancement cannot properly represent the localized strains. This is also observed in the result of the in-plane shear strain component $E_{12}$ (Fig.\,\ref{Z_shape_twist_comparison_strain_field_diff_E12}), where the enhancement is not activated. Since these in-plane components have smaller magnitudes than the others, they may not significantly affect the overall solution. Further investigation into the enhancement of higher-order in-plane strains remains future work.}
\noindent \textcolor{blue}{In Fig.\,\ref{Z_shape_twist_comparison_condition_number}, it is seen that the null space method significantly improves the condition of the tangent stiffness matrix throughout the whole solution process.}
\subsection{Right-angle frame with a prismatic joint}
\label{sec_right_frame}
\begin{figure}[H]
	\centering
	\begin{subfigure}[b]{0.6\textwidth}\centering
		\includegraphics[width=\linewidth]{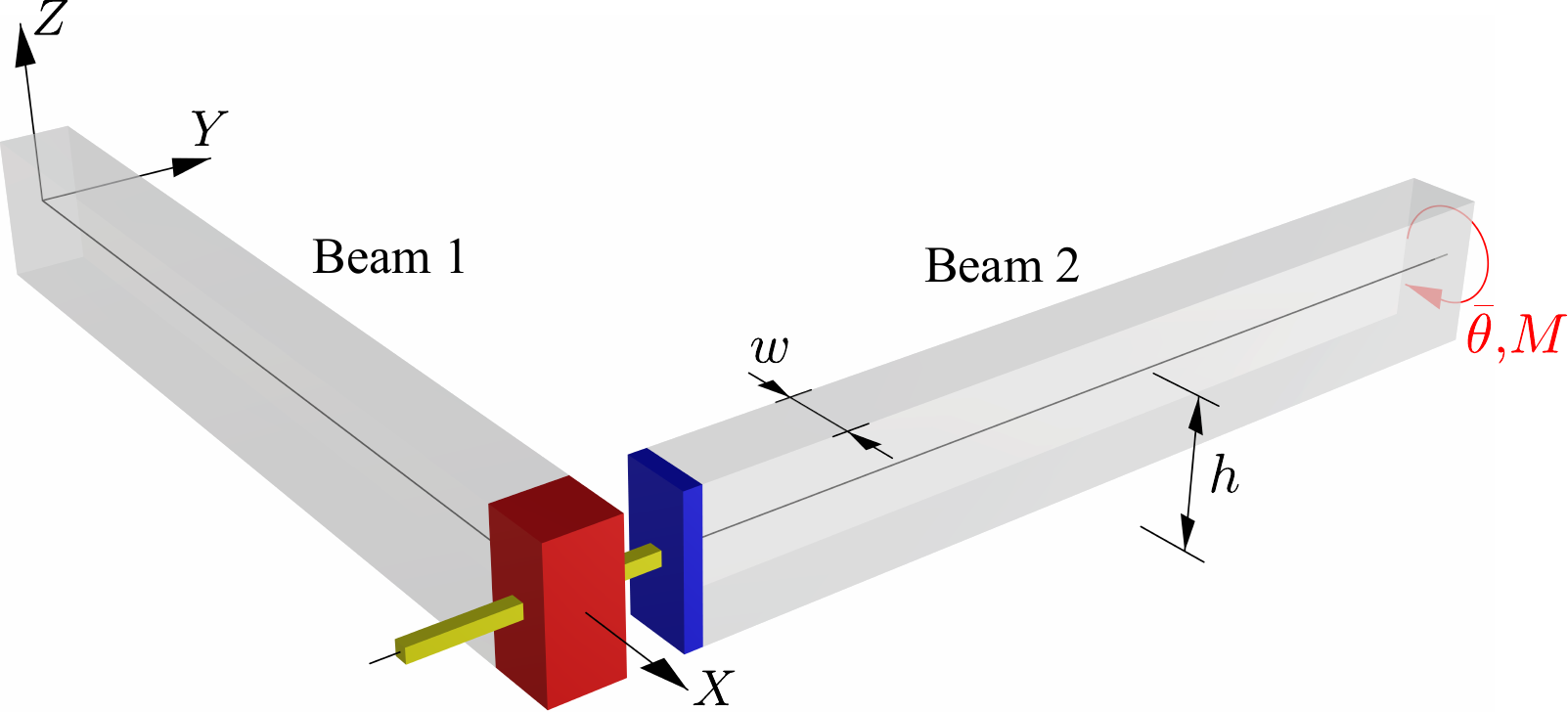}
		\caption{Prismatic joint}
		\label{prismatic_joint_perspective}			
	\end{subfigure}
	\begin{subfigure}[b]{0.35\textwidth}\centering
		\includegraphics[width=\linewidth]{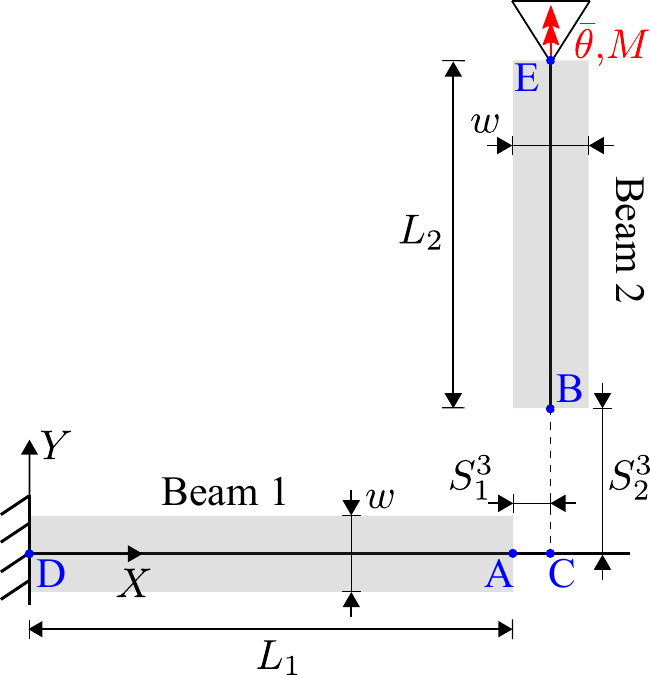}
		\caption{Top view}
		\label{prismatic_joint_planar}			
	\end{subfigure}    
	\caption{Right-angle frame. Initial configuration and boundary conditions. \textcolor{blue}{In (a), the red block represents the massless (virtual) body of a joint with its center denoted by C in (b), and the yellow block illustrates the prismatic part, which allows for relative sliding of Beam 2. The blue block illustrates a massless, rigid panel attached to Beam 2 to transfer the torsional moment to Beam 1 properly. In (b), $S^3_1$ and $S^3_2$ denote the initial offsets.}}
	\label{frame_init_geom_bdc}
\end{figure}

\begin{figure}[H]
	\centering
	\begin{subfigure}[b]{0.475\textwidth}\centering
		\includegraphics[width=\linewidth]{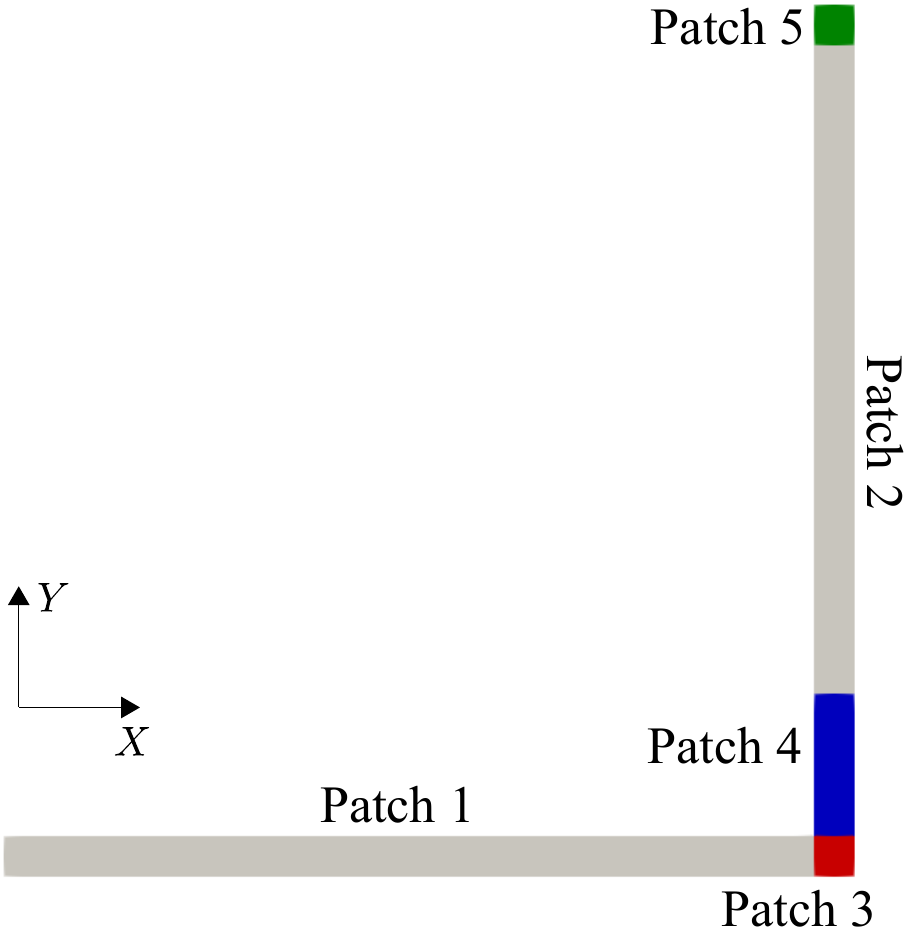}
		\caption{Rigid joint (fully clamped)}
		\label{right_frame_moment_patch_jct_vol_rigid}	
	\end{subfigure}
	\begin{subfigure}[b]{0.475\textwidth}\centering
		\includegraphics[width=\linewidth]{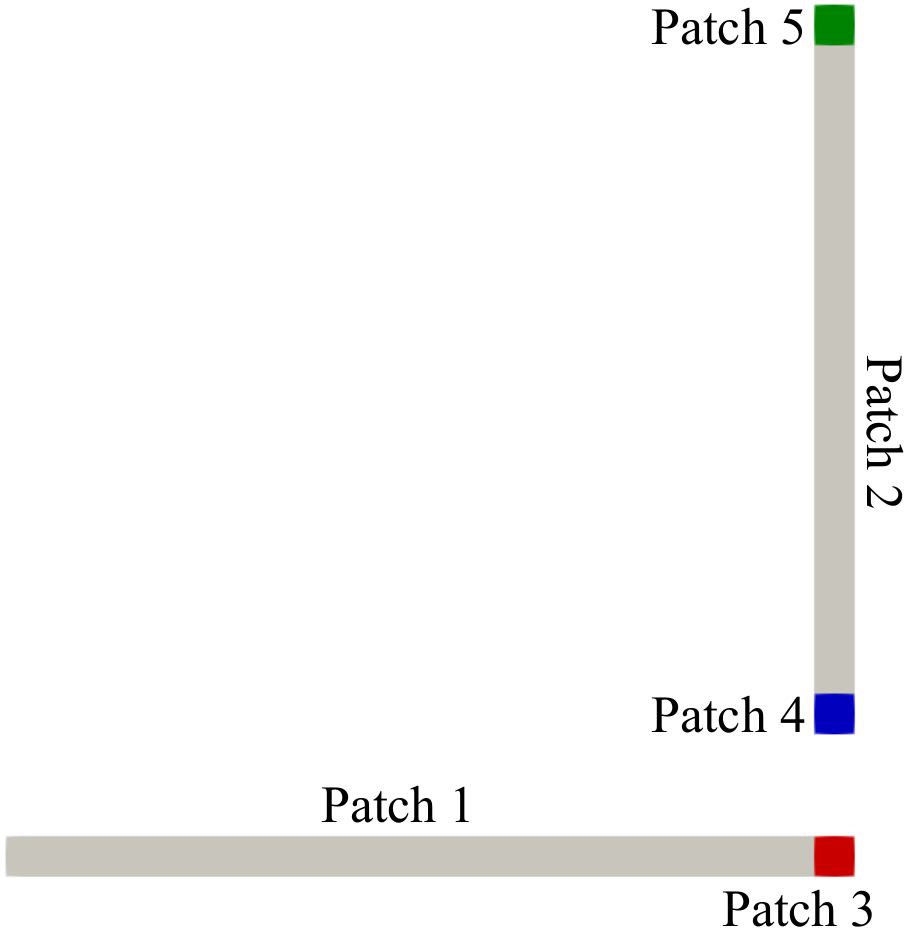}
		\caption{Prismatic joint (fully clamped)}
		\label{right_frame_moment_patch_jct_vol_prismatic}	
	\end{subfigure}    
	\caption{Right-angle frame. \textcolor{blue}{Modeling of fully clamped joints for reference brick solutions by adding Patches 3, 4, and 5 to constrain the cross-sectional warping.}}
	\label{right_frame_moment_patch_jct_vol}
\end{figure}

\begin{figure}[H]
	\centering
	\begin{subfigure}[b]{0.3\textwidth}\centering
		\includegraphics[width=\linewidth]{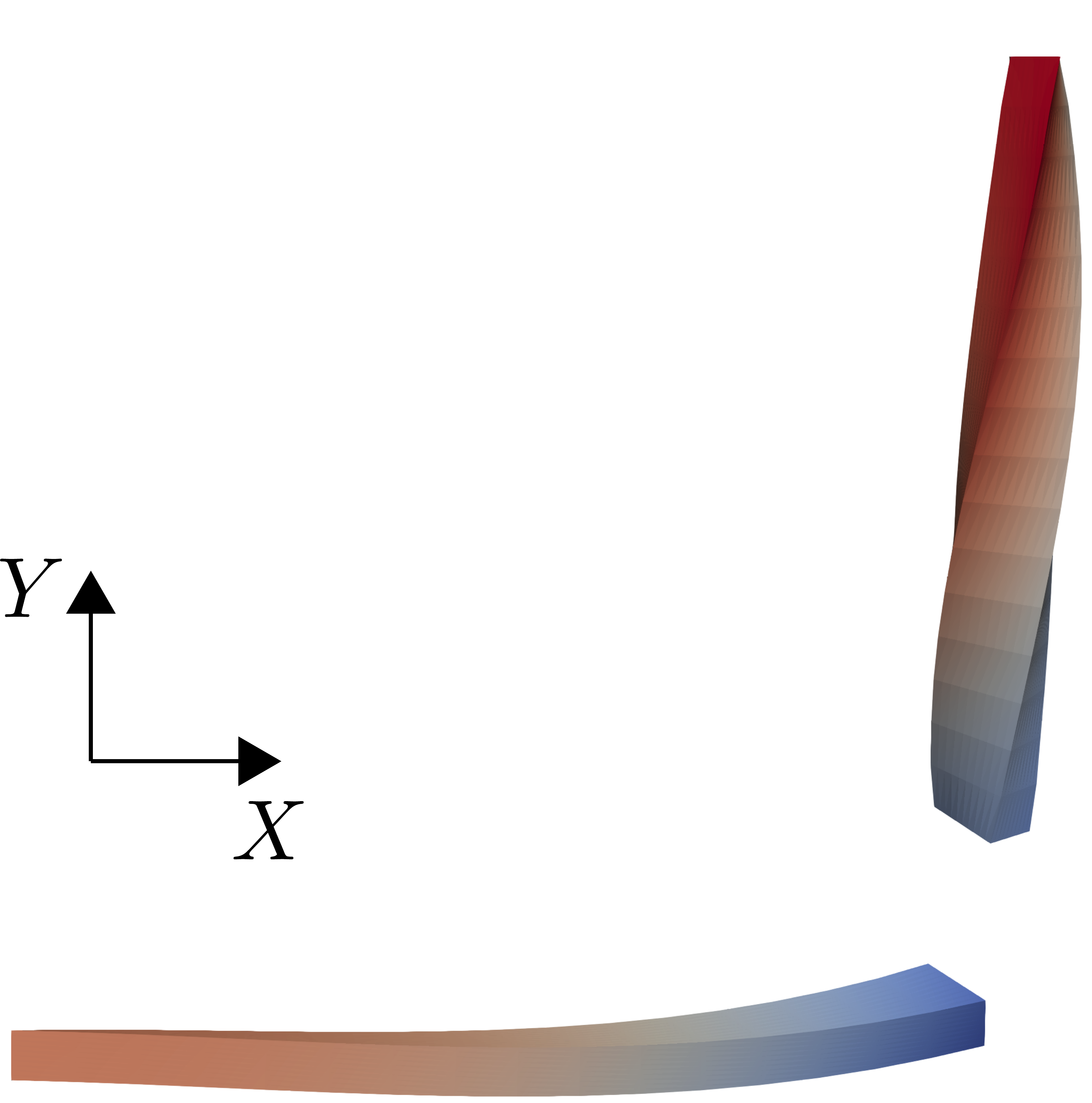}
		\caption{Beam}
		\label{right_angle_deformed_rigid_beam}	
	\end{subfigure}
	\begin{subfigure}[b]{0.3\textwidth}\centering
		\includegraphics[width=\linewidth]{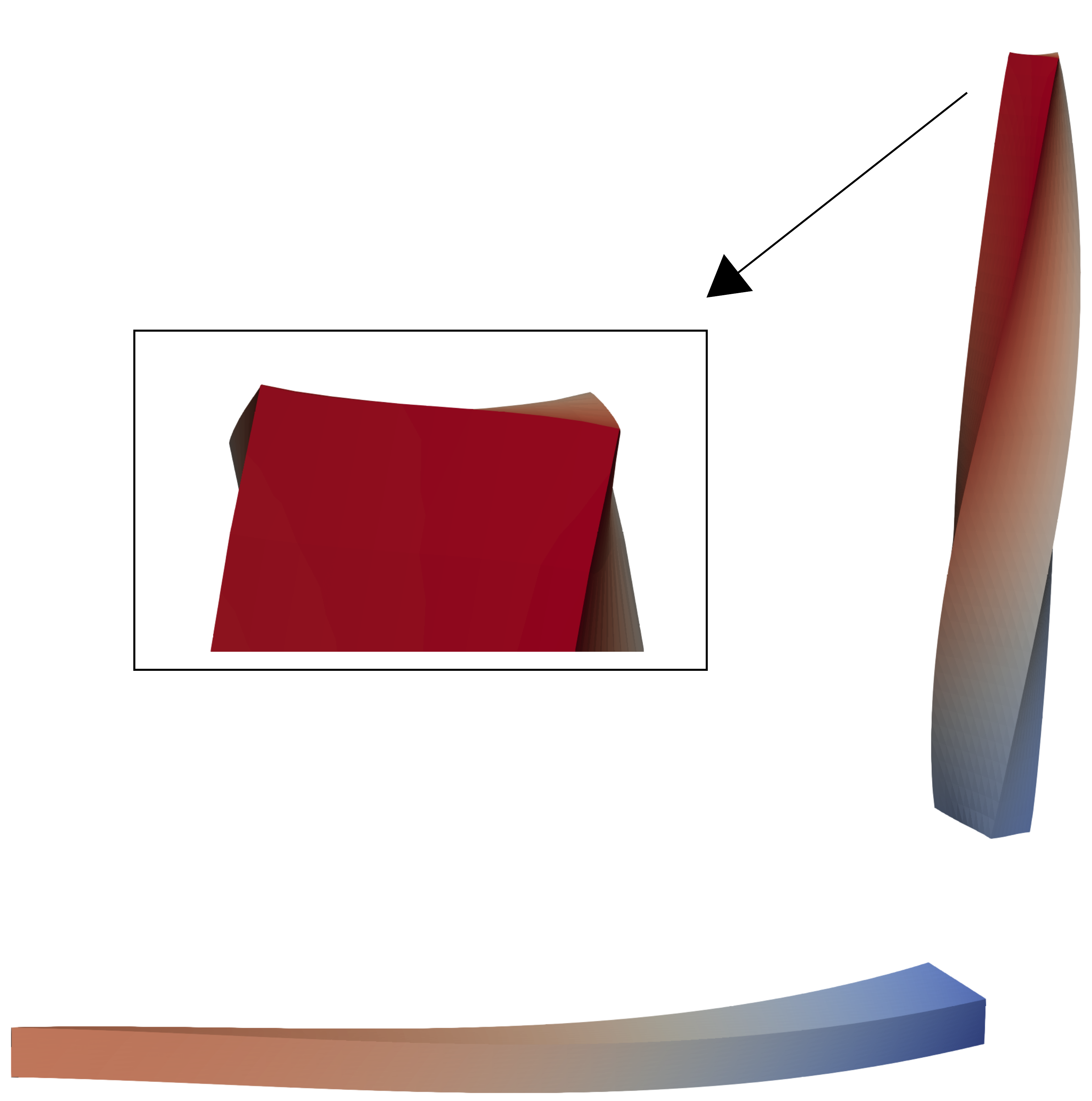}
		\caption{Brick}
		\label{right_angle_deformed_rigid_brick}	
	\end{subfigure}    
	\begin{subfigure}[b]{0.3\textwidth}\centering
		\includegraphics[width=\linewidth]{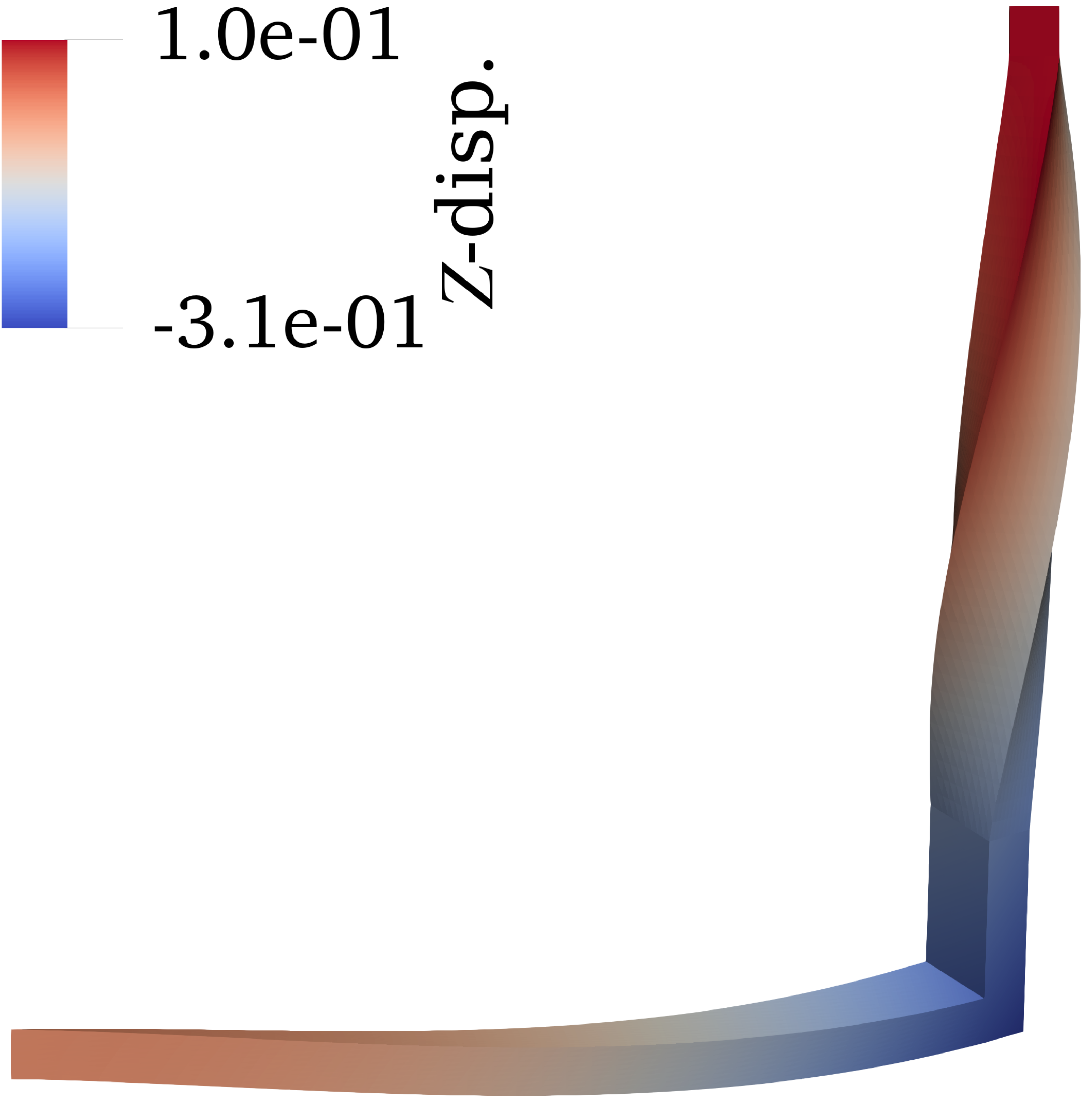}
		\caption{Brick (fully clamped)}
		\label{right_angle_deformed_rigid_brick_full}			
	\end{subfigure}                
	\caption{Right-angle frame (rigid joint). Comparison of the final deformed configurations at $\bar\theta=\pi$ in top view. The colors represent the $Z$-displacements. Table\,\ref{tab_frame_dof_info_rigid_prismatic} compares the degrees-of-freedom. In (c), we see a continuous displacement field in the added volume patches.}
	\label{right_angle_rigid_deformed}	
\end{figure}

\begin{figure}[H]
	\centering
	\begin{subfigure}[b]{0.3\textwidth}\centering
		\includegraphics[width=\linewidth]{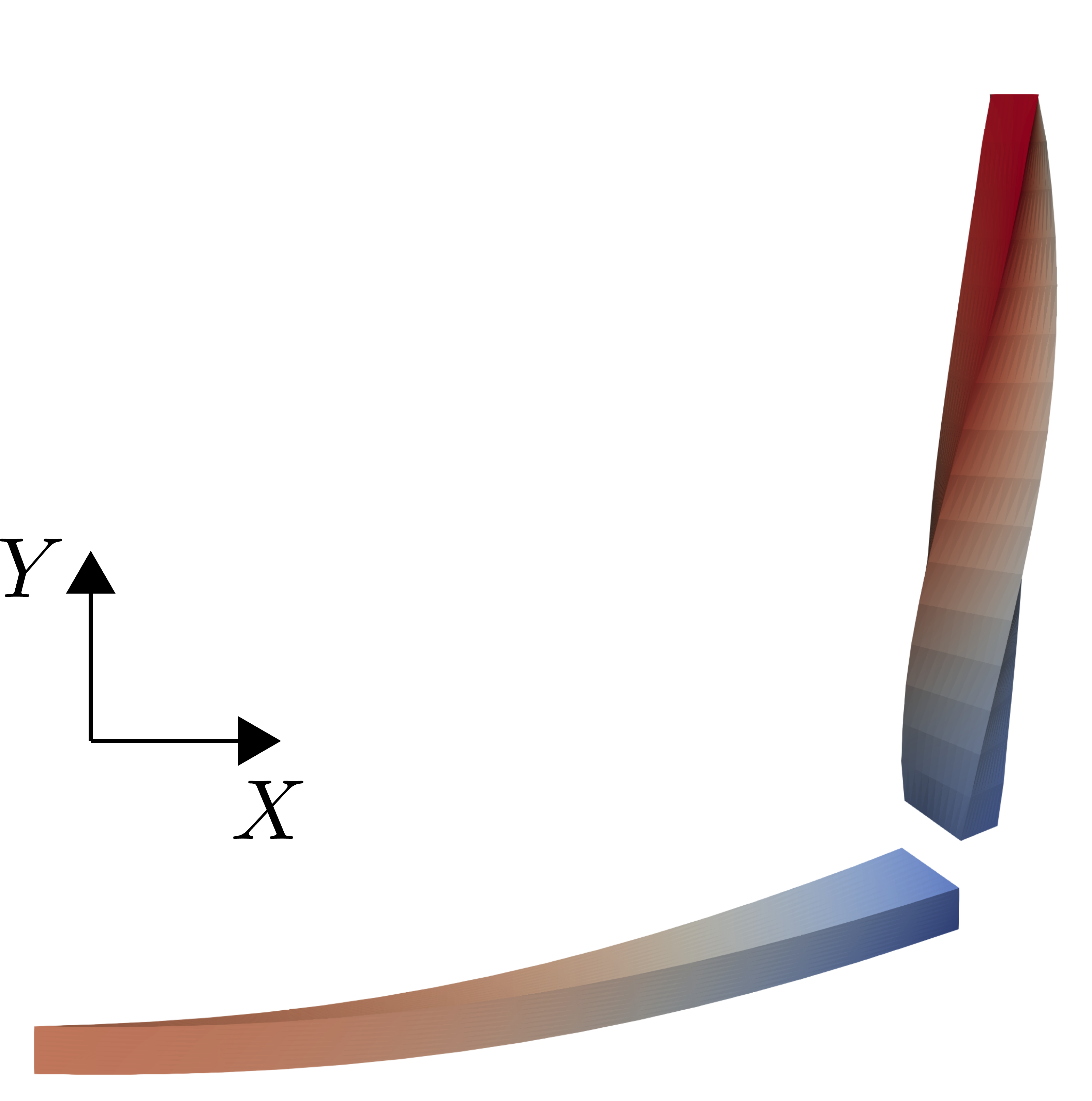}
		\caption{Beam}
		\label{right_angle_prismatic_deformed_beam}			
	\end{subfigure}           
	\begin{subfigure}[b]{0.3\textwidth}\centering
		\includegraphics[width=\linewidth]{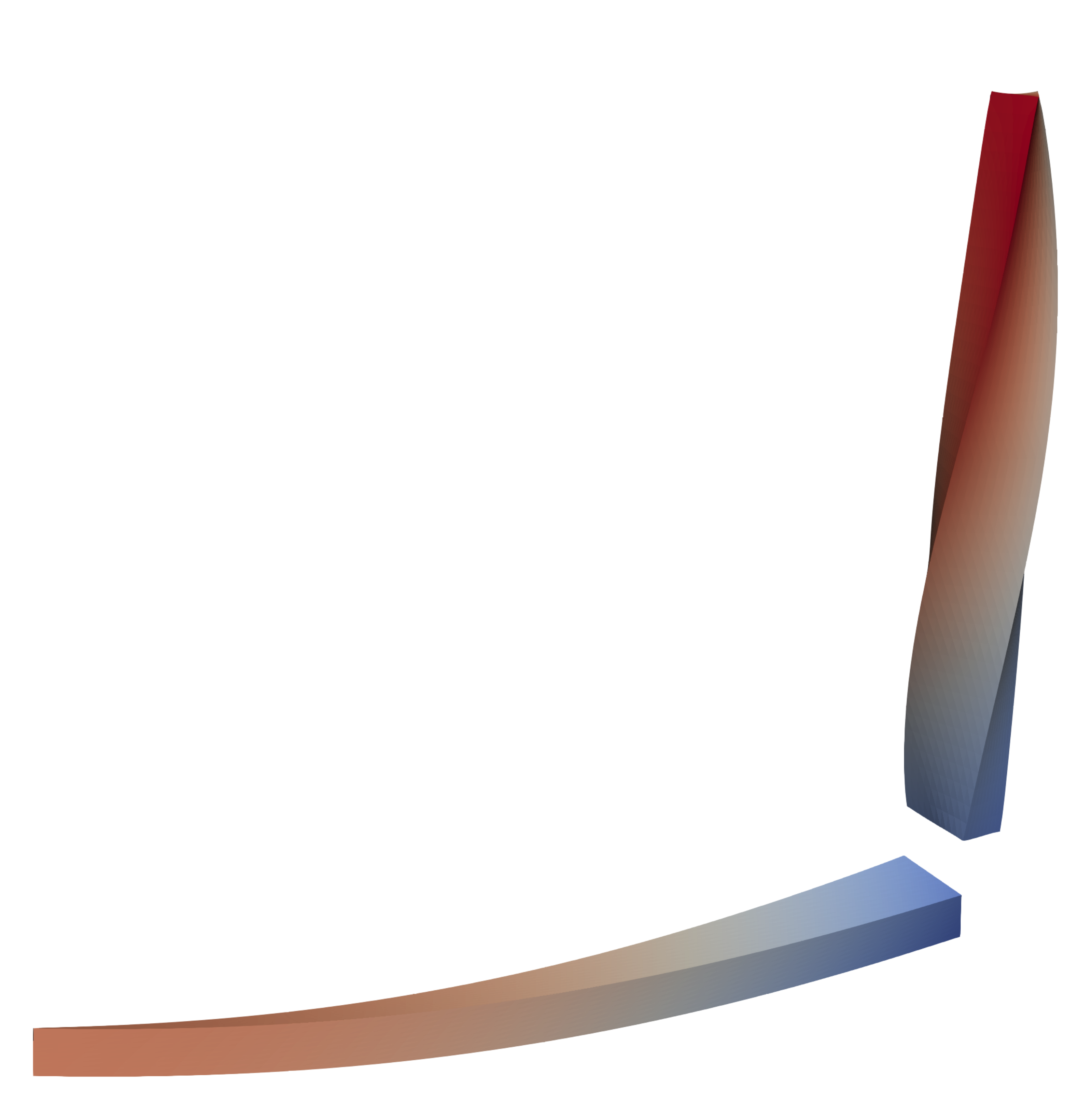}
		\caption{Brick}
		\label{right_angle_prismatic_deformed_brick}			
	\end{subfigure}        
	\begin{subfigure}[b]{0.3\textwidth}\centering
		\includegraphics[width=\linewidth]{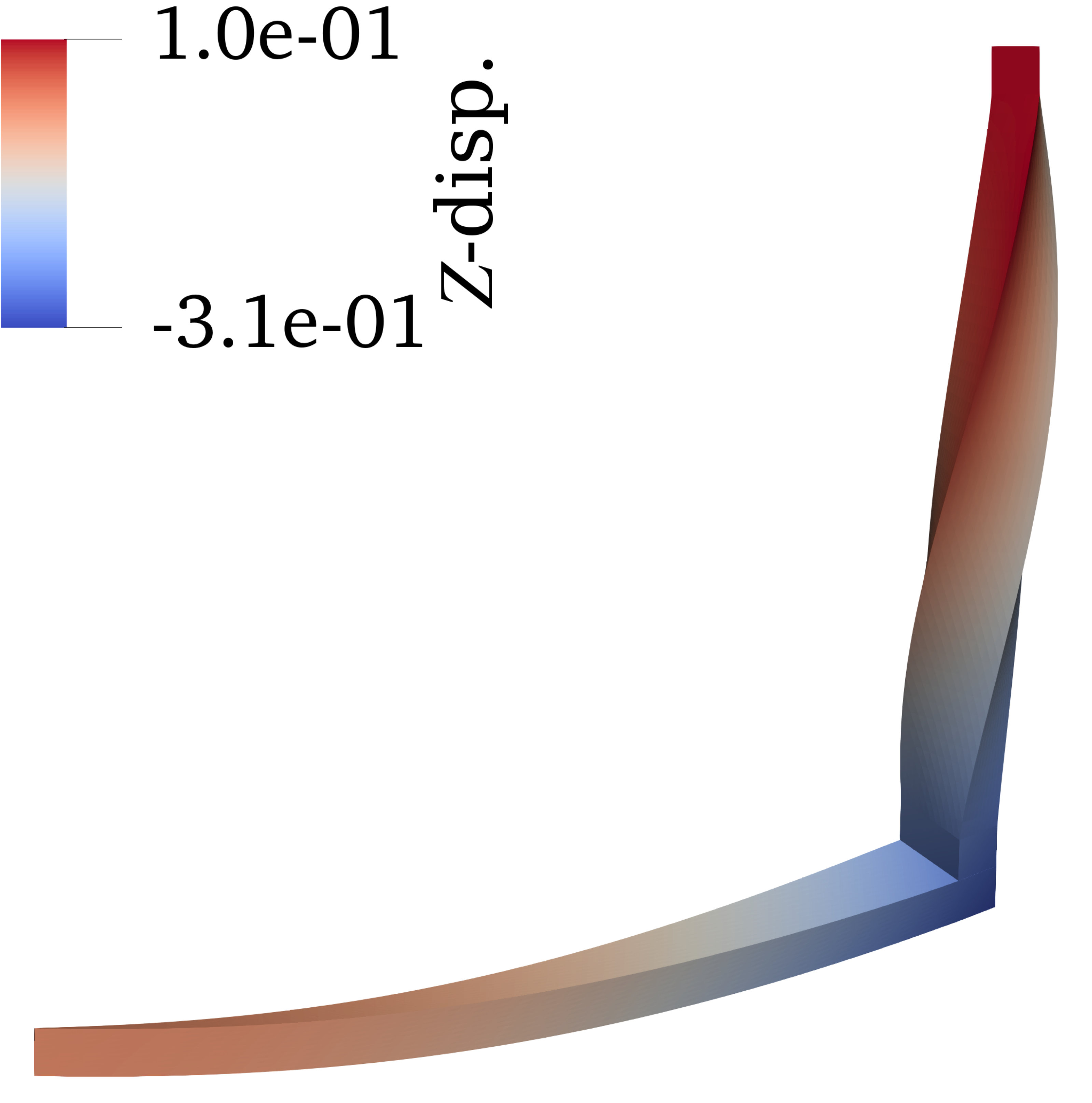}
		\caption{Brick (fully clamped)}
		\label{right_angle_prismatic_deformed_brick_fully_clamped}			
	\end{subfigure}               
	\caption{Right-angle frame (prismatic joint). Comparison of the final deformed configurations in top view. The contour plot represents the $Z$-displacements. Table\,\ref{tab_frame_dof_info_rigid_prismatic} compares the degrees-of-freedom.}
	\label{right_angle_prismatic_deformed}	
\end{figure}

\begin{figure}[H]
	\centering
	\begin{subfigure}[b]{0.475\textwidth}\centering
		\includegraphics[width=\linewidth]{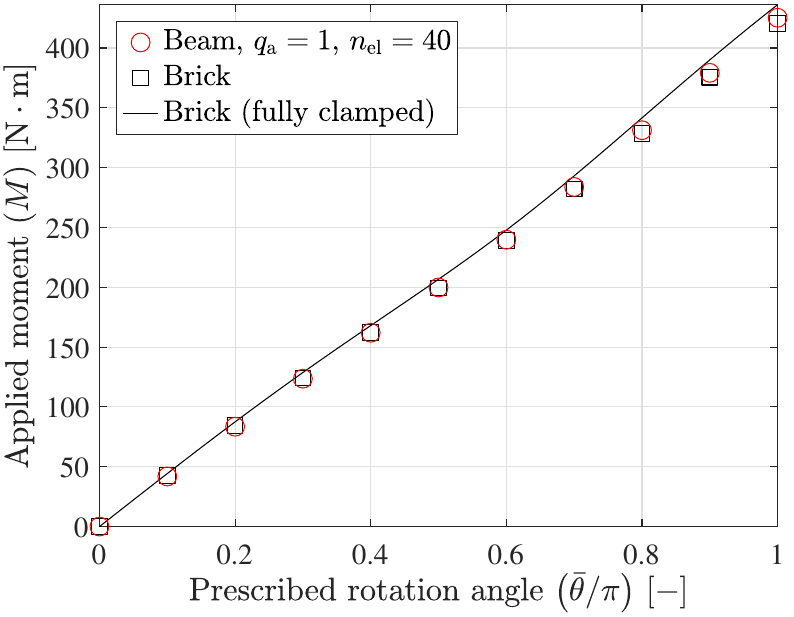}
		\caption{Rigid joint}
		\label{right_frame_moment_mnt_rigid}	
	\end{subfigure}
	\begin{subfigure}[b]{0.475\textwidth}\centering
		\includegraphics[width=\linewidth]{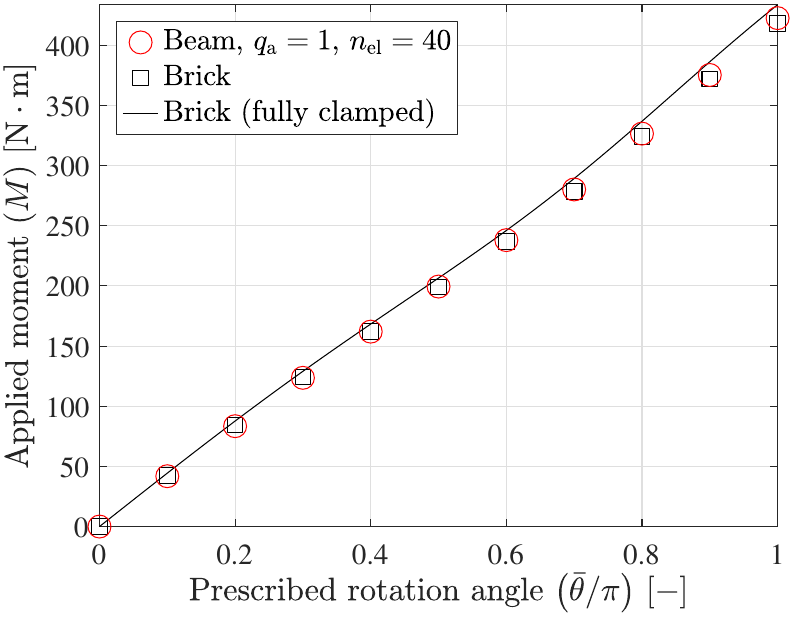}
		\caption{Prismatic joint}
		\label{right_frame_moment_mnt_prismatic}	
	\end{subfigure}    
	\caption{Right-angle frame. Comparison of the applied moment $M$.}
	\label{right_frame_moment_moment_curve}	
\end{figure}

\begin{figure}[H]
	\centering
	\begin{subfigure}[b]{0.475\textwidth}\centering
		\includegraphics[width=\linewidth]{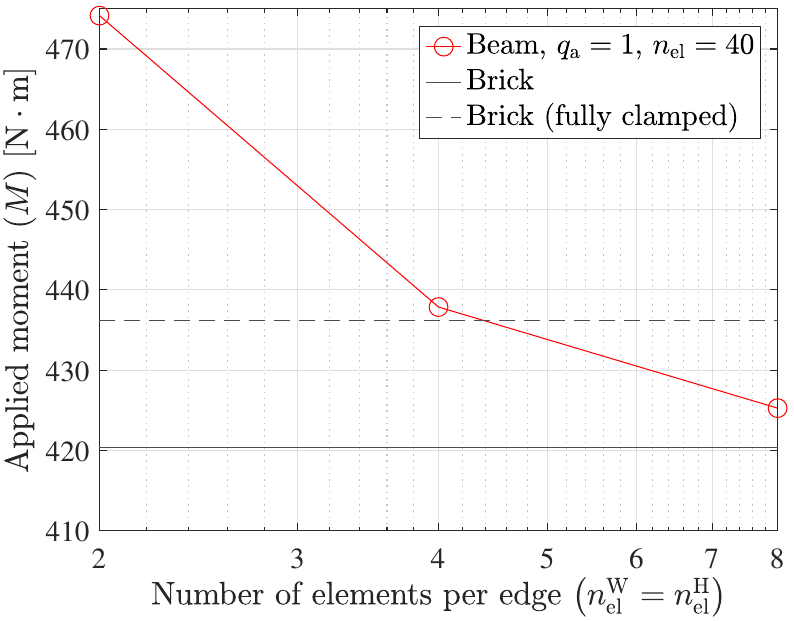}
		\caption{Refinement in the cross-section}
		\label{right_frame_rigid_convergence_moment_refine_cs}	
	\end{subfigure}
	\begin{subfigure}[b]{0.475\textwidth}\centering
		\includegraphics[width=\linewidth]{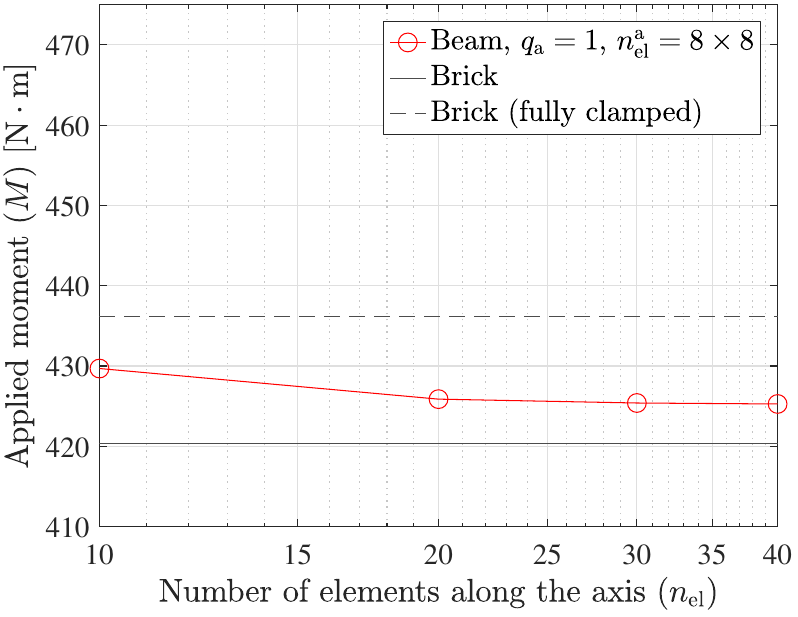}
		\caption{Refinement in the axis}
		\label{right_frame_rigid_convergence_moment_refine_axis}	
	\end{subfigure}
	\caption{Right-angle frame (Case 1: rigid joint). \textcolor{blue}{Convergence of the beam solution for the applied moment $M$ as the number of elements increases in the (a) cross-section $\left(n^\mathrm{a}_\mathrm{el}=n^\mathrm{W}_\mathrm{el}\times{n^\mathrm{H}_\mathrm{el}}\right)$ and (b) axis $\left(n_\mathrm{el}\right)$.}}
	\label{right_frame_rigid_conv}	
\end{figure}

\begin{figure}[H]
	\centering
	\begin{subfigure}[b]{0.475\textwidth}\centering
		\includegraphics[width=\linewidth]{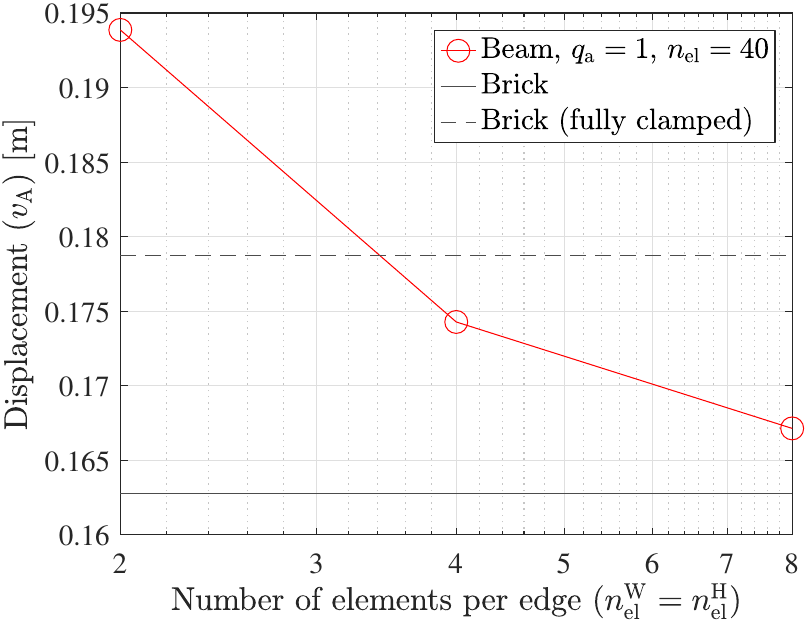}
		\caption{Refinement in the cross-section}
		\label{right_frame_rigid_conv_disp}	
	\end{subfigure}
	\begin{subfigure}[b]{0.475\textwidth}\centering
		\includegraphics[width=\linewidth]{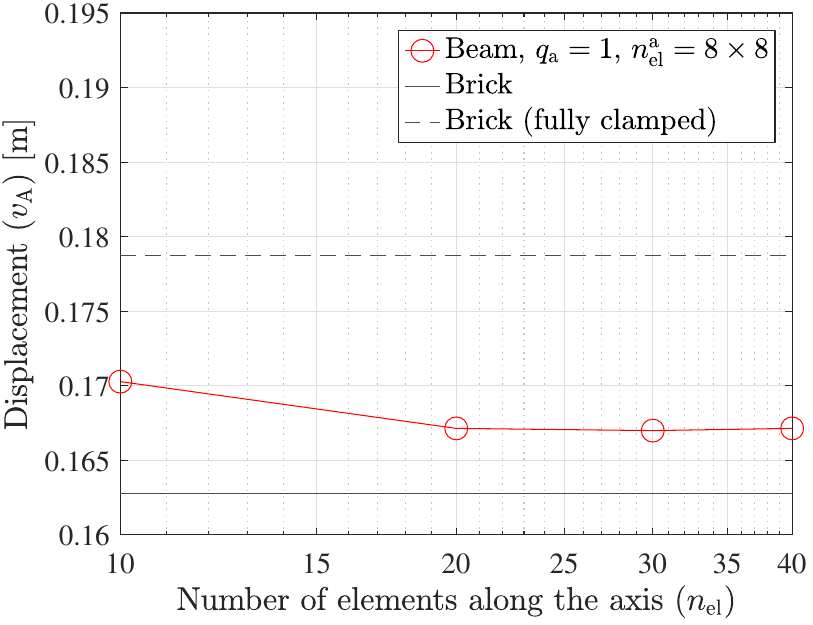}
		\caption{Refinement in the axis}
		\label{right_frame_rigid_conv_disp_cs}	        
	\end{subfigure}    
	\caption{Right-angle frame (Case 2: prismatic joint). \textcolor{blue}{Convergence of the beam solution using $p=3$ for the $Y$-displacement at the point A ($v_\mathrm{A}$), as the number of elements increases in the (a) cross-section and (b) axis.}}
	\label{right_frame_prismatic_conv}	
\end{figure}

\begin{figure}[H]
	\centering
	\begin{subfigure}[b]{0.475\textwidth}\centering
		\includegraphics[width=\linewidth]{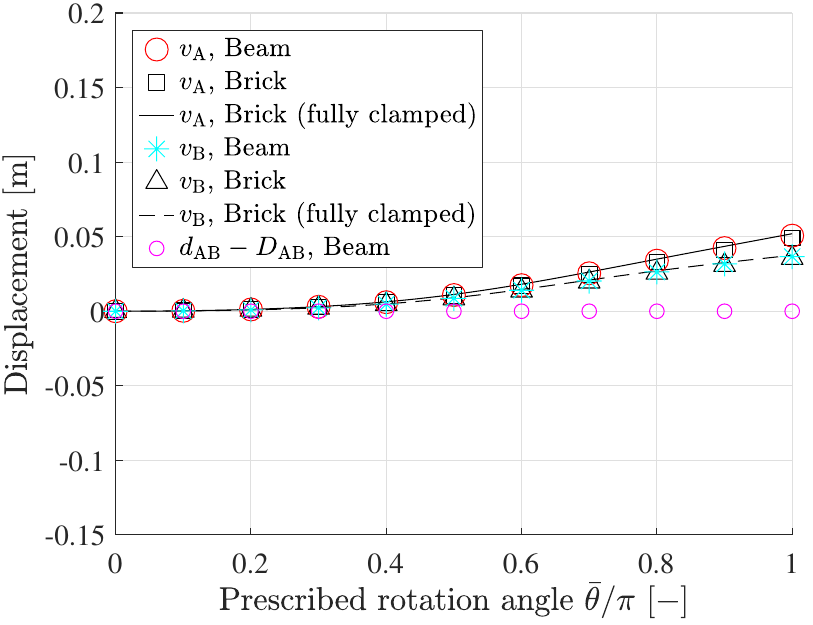}
		\caption{Rigid joint}
		\label{right_frame_moment_disp_rigid}	
	\end{subfigure}
	\begin{subfigure}[b]{0.475\textwidth}\centering
		\includegraphics[width=\linewidth]{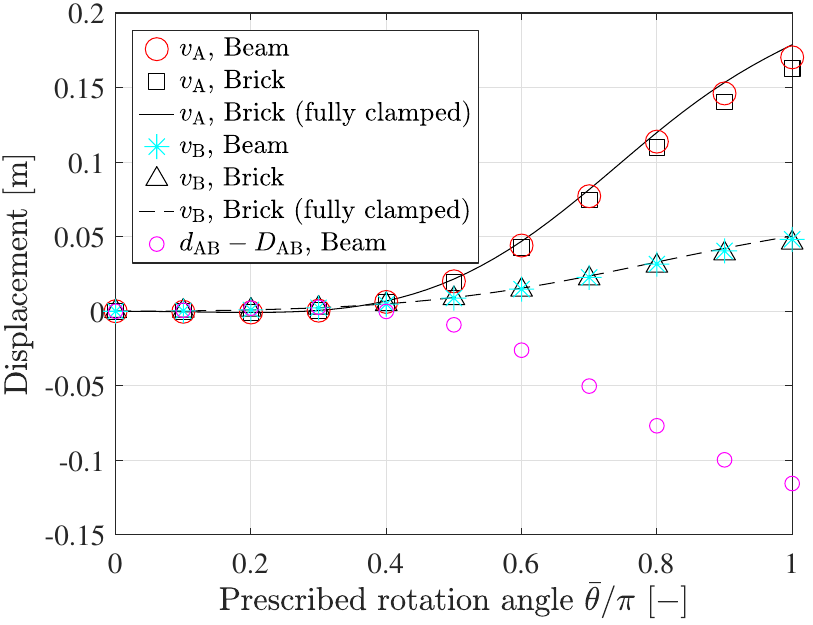}
		\caption{Prismatic joint}
		\label{right_frame_moment_disp_prismatic}	
	\end{subfigure}    
	\caption{Right-angle frame. Comparison of the displacements at the points A and B. $D_\mathrm{AB}$ and $d_\mathrm{AB}$ denote the initial and current distances between A and B, respectively. \textcolor{blue}{Table\,\ref{tab_frame_dof_info_rigid_prismatic} compares the degrees-of-freedom.}}
	\label{right_frame_pre_angle_disp}	
\end{figure}

\begin{figure}[H]
	\centering
	\begin{subfigure}[b]{0.625\textwidth}\centering
		\includegraphics[width=\linewidth]{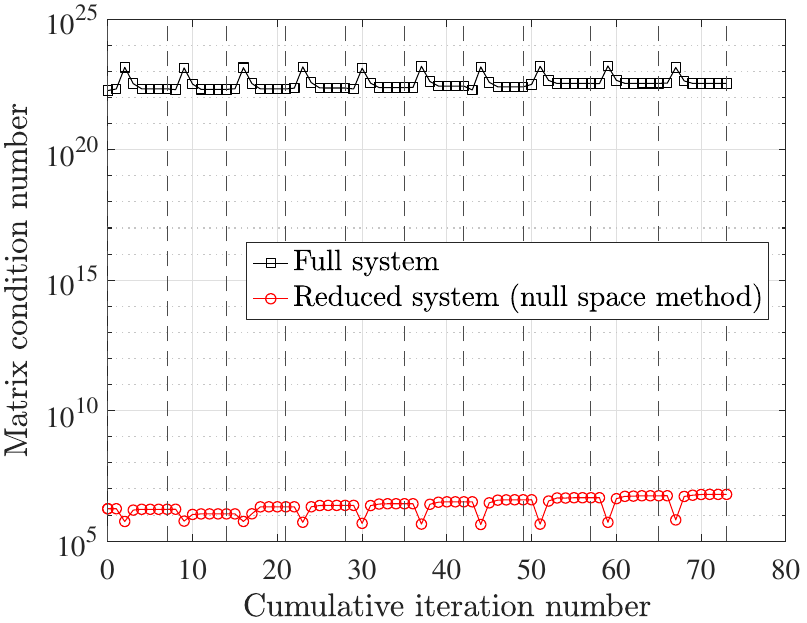}
		\label{right_frame_moment_cond_num}	
	\end{subfigure}
	\caption{Right-angle frame. Comparison of the condition number of the global tangent stiffness matrix between the full and reduced systems for Case 2 (prismatic joint). \textcolor{blue}{Vertical dashed lines define load increments. Here, we have used $p=3$ and $n_\mathrm{el}=10$ in the axis, and $q_\mathrm{a}=1$ and $n^\mathrm{a}_\mathrm{el}=10\times10$ in the cross-section.}}
	\label{right_frame_moment_cond_num_sparse}	
\end{figure}

We consider a frame structure consisting of two initially straight beams, referred to as Beams 1 and 2, connected by a prismatic joint, see Fig.\,\ref{prismatic_joint_perspective} for an illustration. The two beams have the same rectangular cross-section of height $h=0.1\,\mathrm{m}$ and width $w=0.05\,\mathrm{m}$, but different lengths $L_1=1\,\mathrm{m}$ and $L_2=0.8\,\mathrm{m}$, see Fig\,\ref{prismatic_joint_planar}. Beam 1 is aligned with the $X$-axis and is clamped at one end. On the other hand, Beam 2 is aligned with the $Y$-axis and a rotation of angle $\bar\theta=\pi$ about this axis is prescribed at one end, where $M$ denotes the corresponding applied moment. Here, we have introduced offsets $s_1^3$ and $s_2^3$ in the out-of-plane direction of the cross-sections at A and B, respectively, whose initial values are $S_1^3=w/2$ and $S_2^3=4w$. Further, $s_1^3$ is fixed (i.e., $s_1^3=S_1^3$), \textcolor{blue}{and we consider the following two cases for $s_2^3$.
\begin{itemize}
    \item Case 1 (rigid joint): Fixed offset $s_2^3=S_2^3$,
    \item Case 2 (prismatic joint): Released offset $s^3_2$ for a relative sliding between the two beams. 
\end{itemize}
}
For a material model, we consider St.\,Venant-Kirchhoff type hyperelasticity, with Young's modulus $E=1.2\times10^8\,\mathrm{Pa}$ and Poisson's ratio $\nu=0.3$. 
\noindent \textcolor{blue}{Further, as mentioned in Remark\,\ref{rem_partial_clamped}, for a reference brick solution using fully clamped boundary conditions, we have introduced the following additional treatments:
\begin{itemize}
    \item At the fixed end of Beam 1 (point D in Fig.\,\ref{prismatic_joint_perspective}), we fix all degrees-of-freedom.
    \item At the joint, we introduce additional volumes (Patches 3 and 4) to constrain cross-sectional warping, see Fig.\,\ref{right_frame_moment_patch_jct_vol}. For the rigid joint (Case 1), we further couple Patches 3 and 4 for rotational and translational continuity. In contrast, for the prismatic joint (Case 2), we decouple Patches 3 and 4 to allow for relative sliding. An alternative way to implement the prismatic joint is to use anisotropic stiffness in Patch 4, with axial stiffness much lower than the other stiffnesses; see \citet[Table 5]{jelenic1996non}.
    \item At the loaded end (point E), we also add Patch 5 to constrain the warping. 
    \item We increase the Young's modulus in the additional volumes (Patches 3, 4, and 5) by a factor of 100, for constraining the cross-sectional warping.
\end{itemize}
}
\noindent In Fig.\,\ref{right_angle_rigid_deformed}, we compare the deformed configurations between the beam and brick solutions. \textcolor{blue}{The torsional motion of deflected Beam 2 induces the $Y$-directional deflection in Beam 1. Here, it is also noted that the deformed configurations are plotted using only the global degrees-of-freedom. For the beam solution, the cross-sectional warping due to the internal degrees-of-freedom is not visualized. Despite this, the deformed configuration from the beam formulation agrees very well with those from the brick formulations. For the brick solution with a partially clamped condition, out-of-plane cross-sectional warping is noticeable at the loaded end, as shown in the magnified view (small box) in Fig.\,\ref{right_angle_deformed_rigid_brick}. This warping is constrained in the brick solution under fully clamped conditions (see Fig.\,\ref{right_angle_deformed_rigid_brick_full}), implemented by increasing the stiffness of the additional volume (Patch 5). Furthermore, the added volumes (Patches 3 and 4) with increased stiffness properly represent the fully rigid joints. In Fig.\,\ref{right_angle_prismatic_deformed}, it is noticeable that Beam 1 slides toward Beam 2, due to the released offset.}
\noindent \textcolor{blue}{In Fig.\,\ref{right_frame_moment_moment_curve}, we compare the applied moment ($M$), as the prescribed rotation ($\bar \theta$) increases. The fully clamped condition for the brick solution yields a larger moment (i.e., stiffer behavior) than the partially clamped condition does. It is also noticeable that the beam solutions lie between those of the partially clamped and fully clamped brick models. This is consistent with the fact that the beam solution approaches the brick solution under partially clamped conditions, as the number of elements in the cross-section increases (see Figs.\,\ref{right_frame_rigid_convergence_moment_refine_cs} and \ref{right_frame_rigid_conv_disp}). Note that mesh refinement in the cross-section does not increase the number of global degrees-of-freedom due to static condensation; see Table \ref{tab_frame_dof_info_rigid_prismatic}. In the following, we use the brick solutions with $\mathrm{deg.}=(3,3,3)$ and $n_\mathrm{el}=50\times5\times5$ as reference solutions.}
\noindent \textcolor{blue}{In Figs.\,\ref{right_frame_rigid_conv_disp} and \ref{right_frame_rigid_conv_disp_cs}, we show the convergence of the beam solution for the $Y$-displacement at point A, as we refine the mesh in the cross-section and axis, respectively. Since the applied moment decreases with the mesh refinement, the $Y$-displacement decreases as well. It is also seen that the converged beam solutions are much closer to the partially clamped brick solution than to the fully clamped one. In both Figs.\,\ref{right_frame_rigid_conv} and \ref{right_frame_prismatic_conv}, the converged beam solutions do not perfectly agree with the brick solutions under partially clamped conditions. Further investigation on this issue remains future work; see the relevant comments in Observation\,\ref{observ_ex_z_torsion}.}
\noindent Further, in Fig.\,\ref{right_frame_pre_angle_disp}, we compare the displacements at points A and B as the prescribed rotation increases. For the rigid joint, the distance between points A and B is maintained (pink circles), even though their $Y$-displacements ($v_\mathrm{A}$ and $v_\mathrm{B}$) are slightly different due to the joint's rotation. In contrast, the prismatic joint allows for relative sliding between A and B (pink curve in Fig.\,\ref{right_frame_moment_disp_prismatic}), leading to larger displacement at B. It is also seen that the slightly larger moments in the brick (fully clamped) and beam solutions, seen in Fig.\,\ref{right_frame_moment_moment_curve}, yield the slightly larger $Y$-displacements than that of the partially clamped brick solution.
\noindent In Fig.\,\ref{right_frame_moment_cond_num_sparse}, we investigate the condition number of the global tangent stiffness matrix for the beam solution. It is seen that the full (original) system is extremely ill-conditioned, which is resolved here by using the present null space method. 
\noindent \textcolor{blue}{In Table\,\ref{tab_frame_dof_info_rigid_prismatic}, we compare the degrees-of-freedom between beam and brick formulations. In the beam formulation, the degrees-of-freedom of the enhanced warping strain field are eliminated from the global system via static condensation and treated as internal variables, resulting in a much smaller system of equations compared to the brick formulation. The brick solution with fully clamped conditions has more degrees-of-freedom than the partially clamped one, due to the additional volumes (Patches 3, 4, and 5). In Appendix \ref{app_num_ex_right_frame}, Tables \ref{tab_frame_brick_conv_rigid} and \ref{tab_frame_brick_conv_prismatic} show the convergence behavior for the brick solutions.}
\begin{table}[H]
  \centering
  \footnotesize
  \caption{\textcolor{blue}{Right-angle frame. Comparison of degrees-of-freedom. Here, `r' and `p' denote the Cases 1 (rigid joint) and 2 (prismatic joint), respectively, and `f' denotes the fully clamped conditions. Element counts are defined per patch.}}
    \begin{tabular}{lccccccccc}
    \toprule
    \multicolumn{1}{r}{} & \multicolumn{2}{c}{Degrees} &       & \multicolumn{2}{c}{Elements} &       & \multicolumn{3}{c}{Degrees-of-freedom} \\
\cmidrule{2-3}\cmidrule{5-6}\cmidrule{8-10}    \multicolumn{1}{r}{} & L     & \multicolumn{1}{c}{W,H} &       & L     & \multicolumn{1}{c}{W,H} &       & Global\,(r) & Global\,(p) & Internal\,(r,p) \\
    \midrule
    \multicolumn{1}{l}{Beam} & 3     & 1     &       & 10    & 10    &       & 234   & 234   & 28680 \\
    \multicolumn{1}{l}{Brick} & 3     & 3     &       & 50    & 5     &       & 20352 & 20352 & $-$ \\
    Brick (f) & 3     & 3     &       & 50    & 5     &       & {22656} & {22848} & $-$ \\
    \bottomrule
    \end{tabular}%
  \label{tab_frame_dof_info_rigid_prismatic}%
\end{table}%
\subsection{Framed shallow dome}
\label{sec_num_ex_shallow}
\begin{figure}[H]
	\centering
	\begin{subfigure}[b]{0.65\textwidth}\centering
		\includegraphics[width=\linewidth]{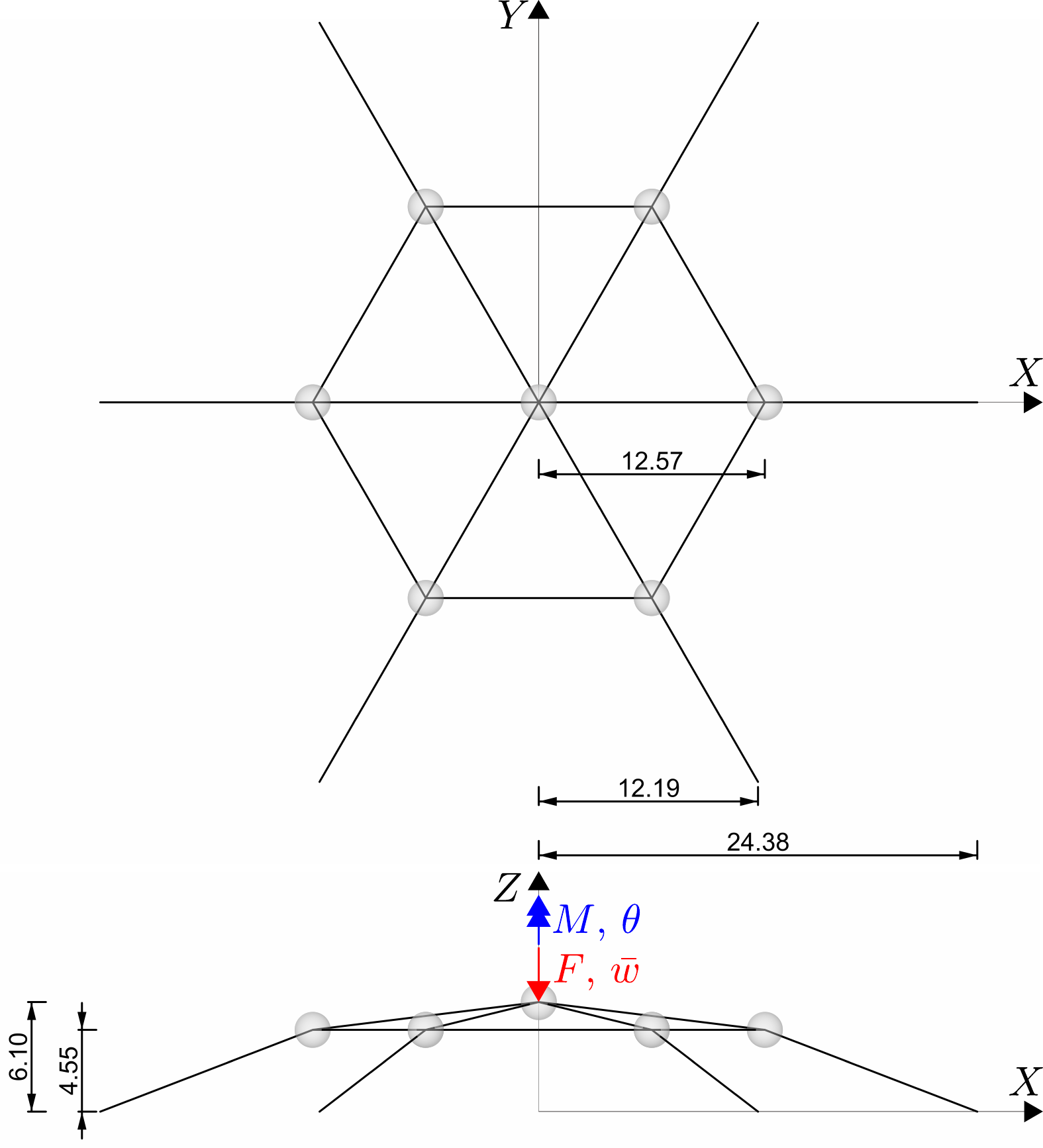}	
	\end{subfigure}
	\caption{Shallow dome. Initial configuration in planar views. Here, seven spheres of radius $1\,\mathrm{m}$ are introduced for the trimming at the junctions. All dimensions are in meters.}
	\label{shallow_dome_init_bdc_load_cond}	
\end{figure}

\begin{figure}[H]
	\centering
	\begin{subfigure}[b]{0.45\textwidth}\centering
		\includegraphics[width=\linewidth]{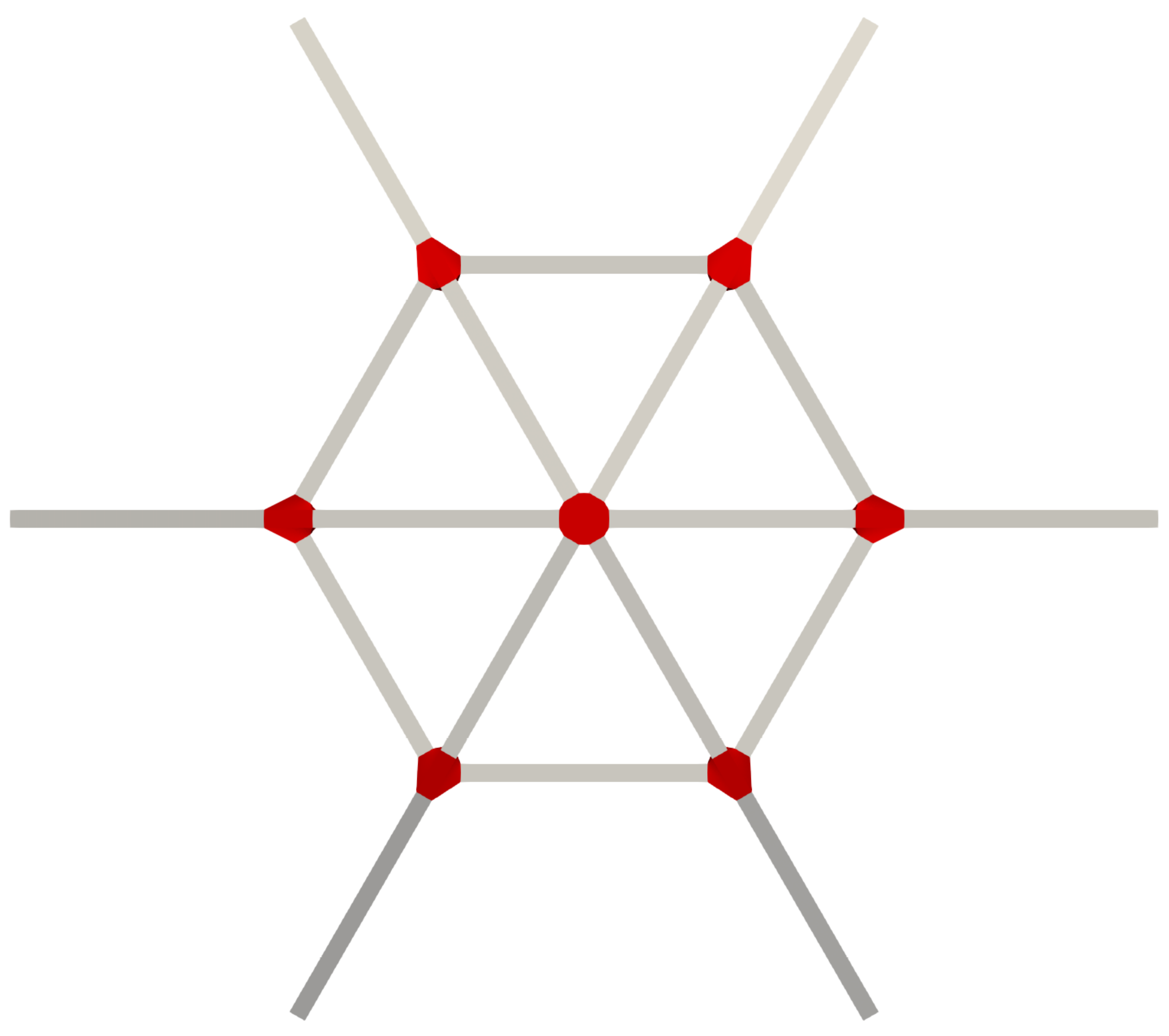}
		\caption{Top view}
	\end{subfigure}
	\begin{subfigure}[b]{0.45\textwidth}\centering
		\includegraphics[width=\linewidth]{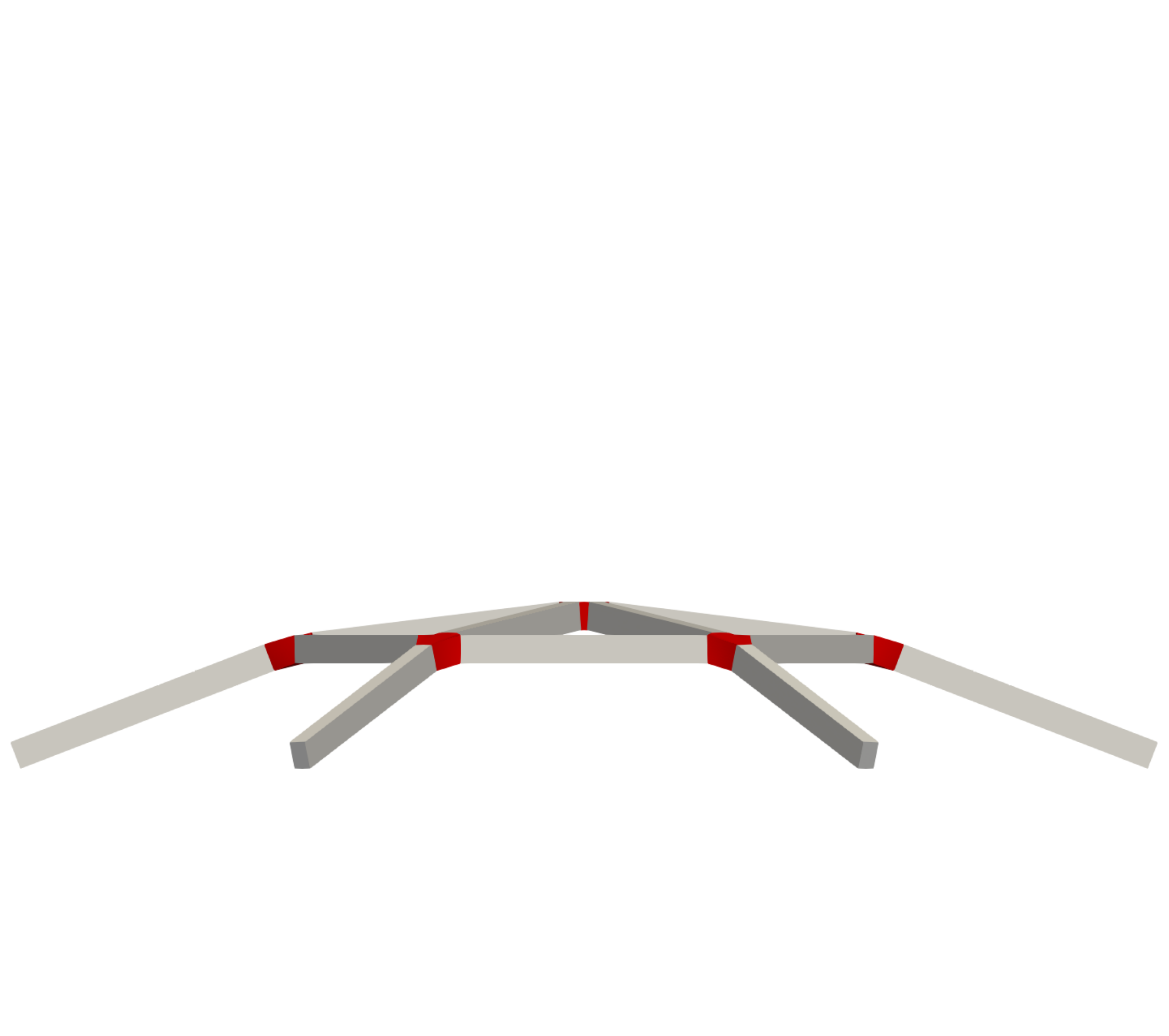}
		\caption{Side view}
	\end{subfigure}    
	\caption{Shallow dome. \textcolor{blue}{Initial configuration of the {fully clamped} brick model for a {reference solution}. The red parts are the added volume patches, whose geometry is defined by connecting the vertices of the beams (gray blocks).}}
	\label{shallow_dome_init_jct_vol}	
\end{figure}

\begin{figure}[H]
	\centering
	\begin{subfigure}[b]{0.325\textwidth}\centering
		\includegraphics[width=\linewidth]{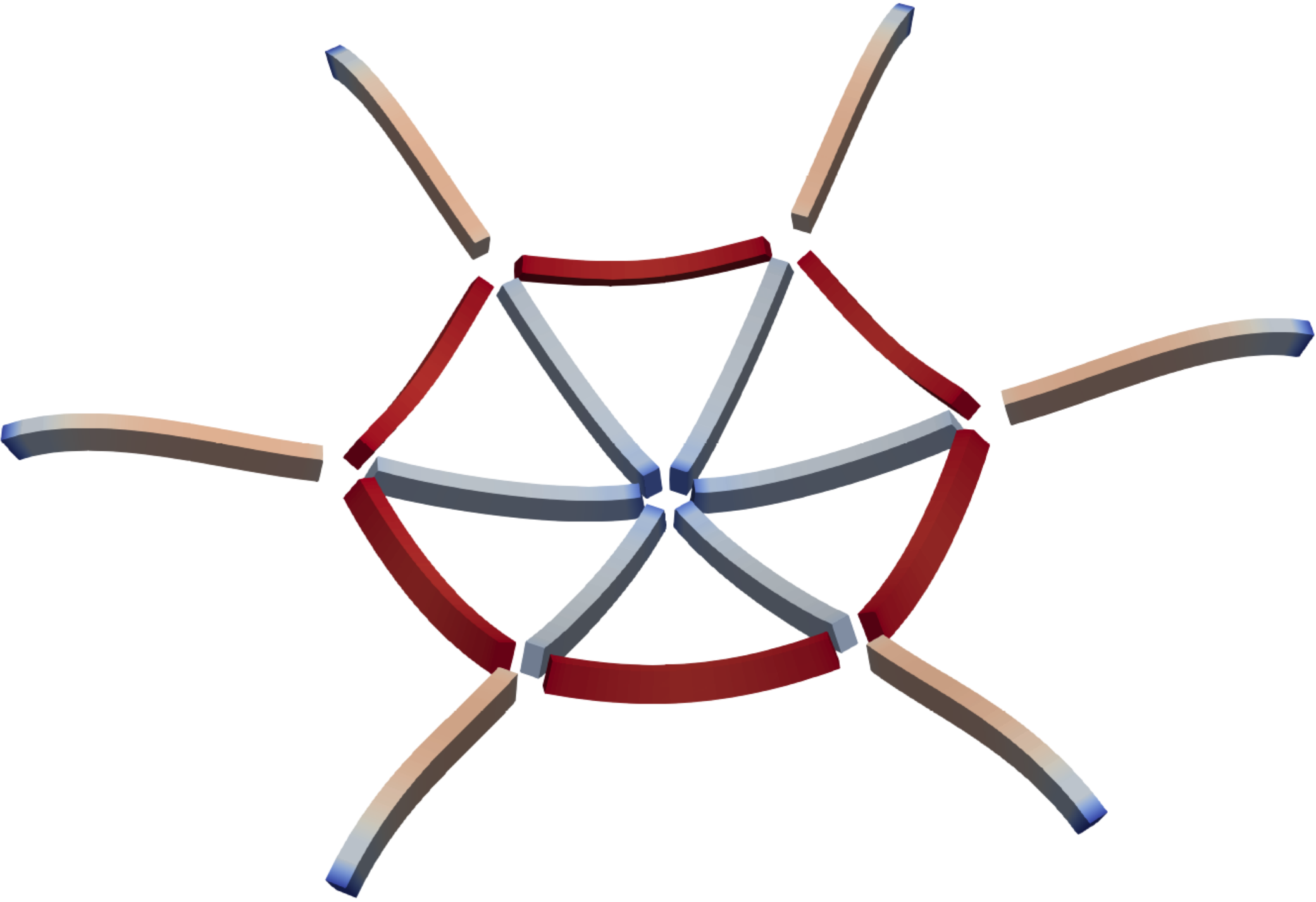}
		\caption{Beam}
		\label{shallow_case1_deformed_compare_beam}			
	\end{subfigure}        
	\begin{subfigure}[b]{0.325\textwidth}\centering
		\includegraphics[width=\linewidth]{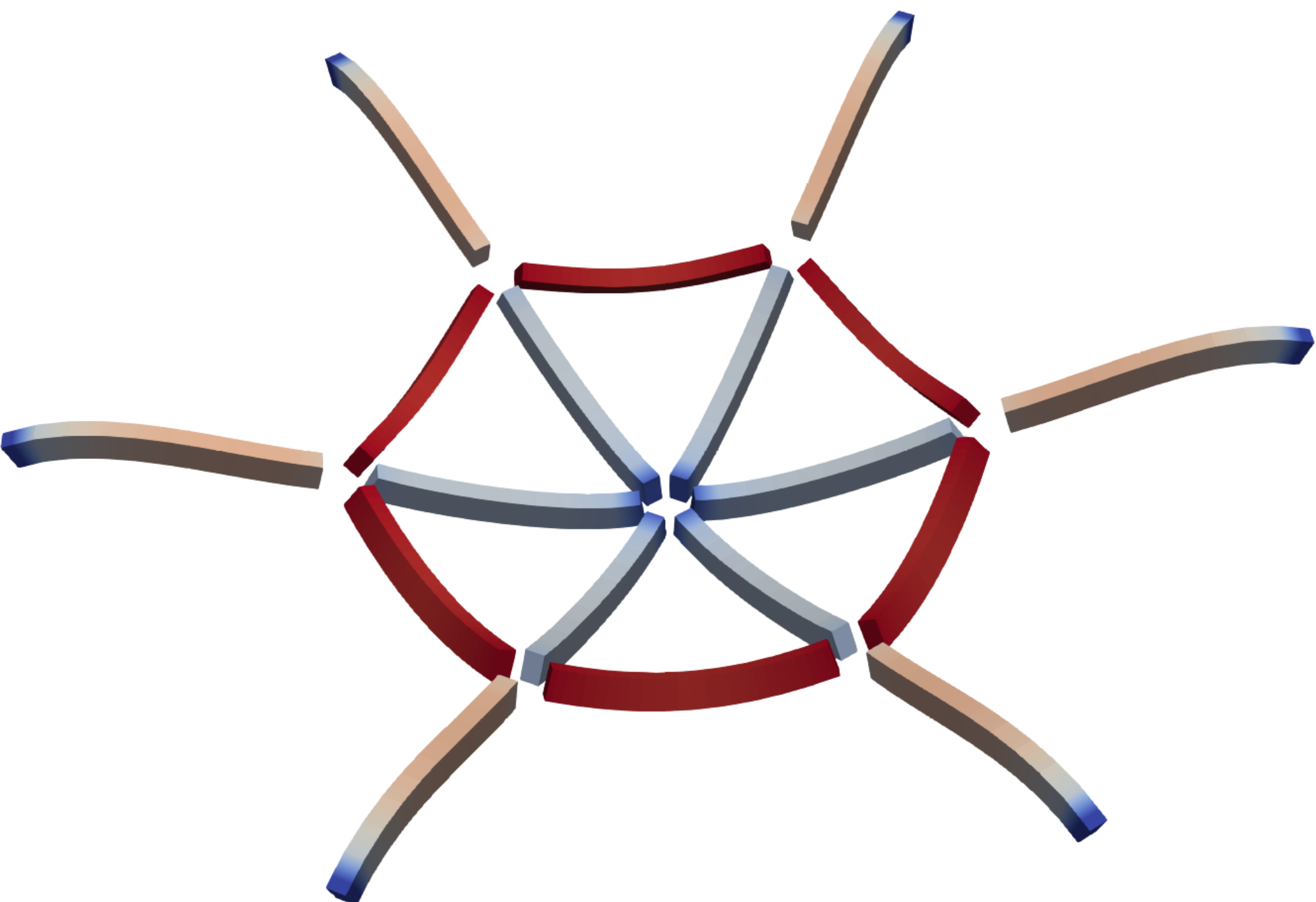}
		\caption{Brick}
		\label{shallow_case1_deformed_compare_brick}			
	\end{subfigure}
	\begin{subfigure}[b]{0.325\textwidth}\centering
		\includegraphics[width=\linewidth]{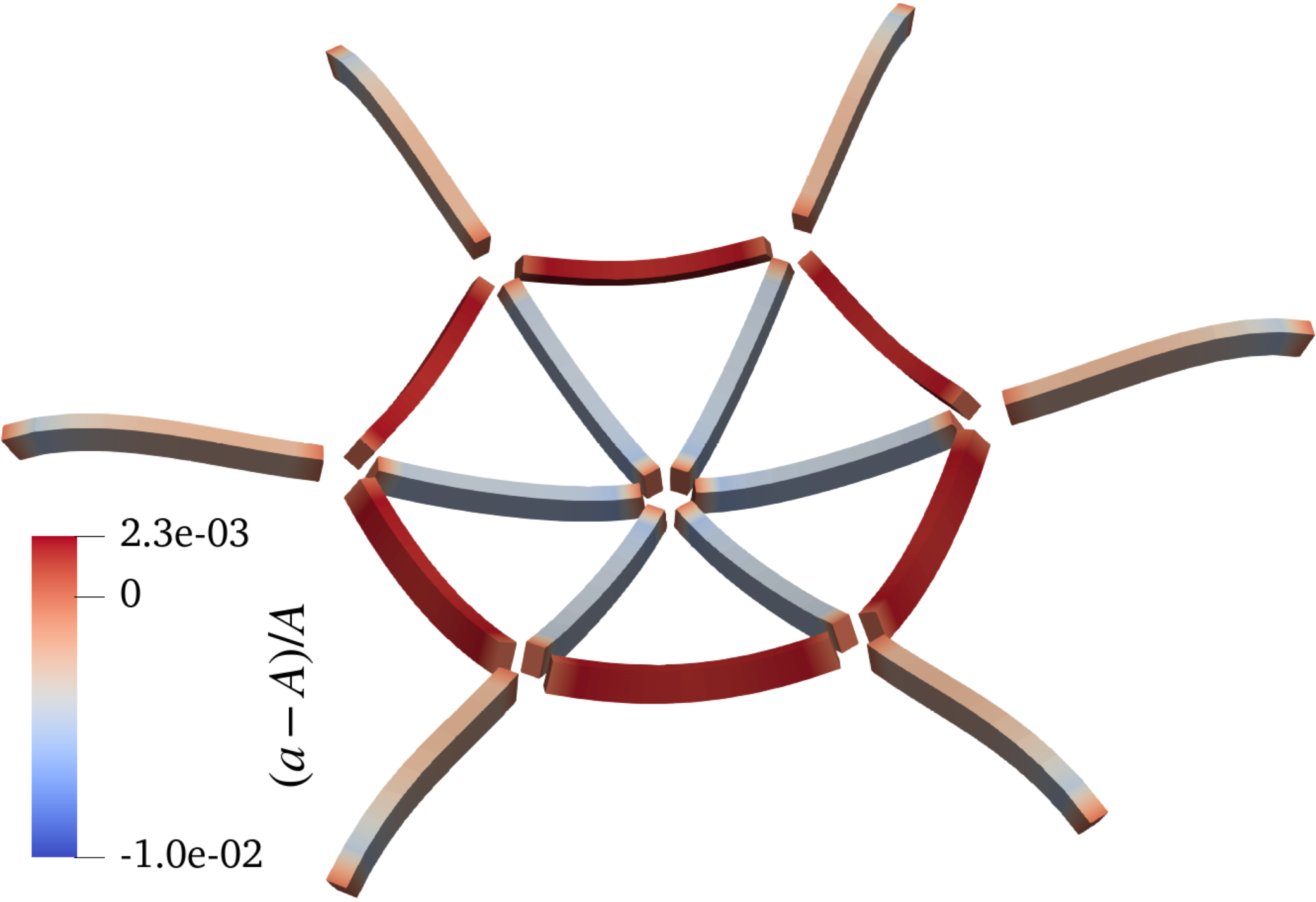}
		\caption{Brick (fully clamped)}
		\label{shallow_case1_deformed_compare_brick_full}			
	\end{subfigure}    
	\caption{Shallow dome (Case 1: displacement load only). Comparison of the final deformed configurations between the beam and brick solutions. The colors represent the relative change in cross-sectional area, where $a$ and $A$ denote the initial and current areas. In (c), the added volume patches are not visualized.}
	\label{shallow_case1_deformed_compare}	
\end{figure}

\begin{figure}[H]
	\centering
	\begin{subfigure}[b]{0.325\textwidth}\centering
		\includegraphics[width=\linewidth]{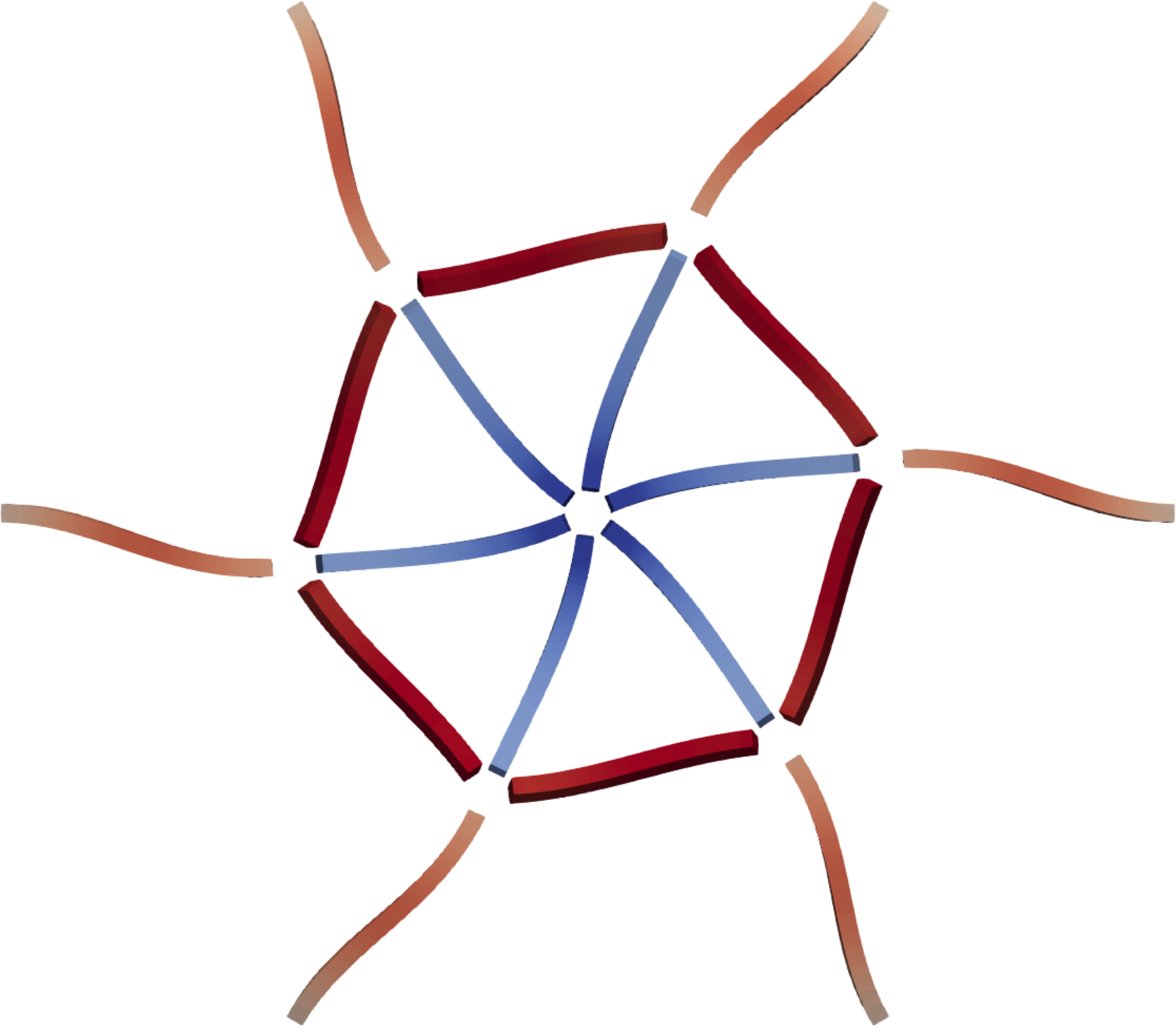}
		\caption{Beam}
		\label{shallow_case2_deformed_compare_beam}			
	\end{subfigure}        
	\begin{subfigure}[b]{0.325\textwidth}\centering
		\includegraphics[width=\linewidth]{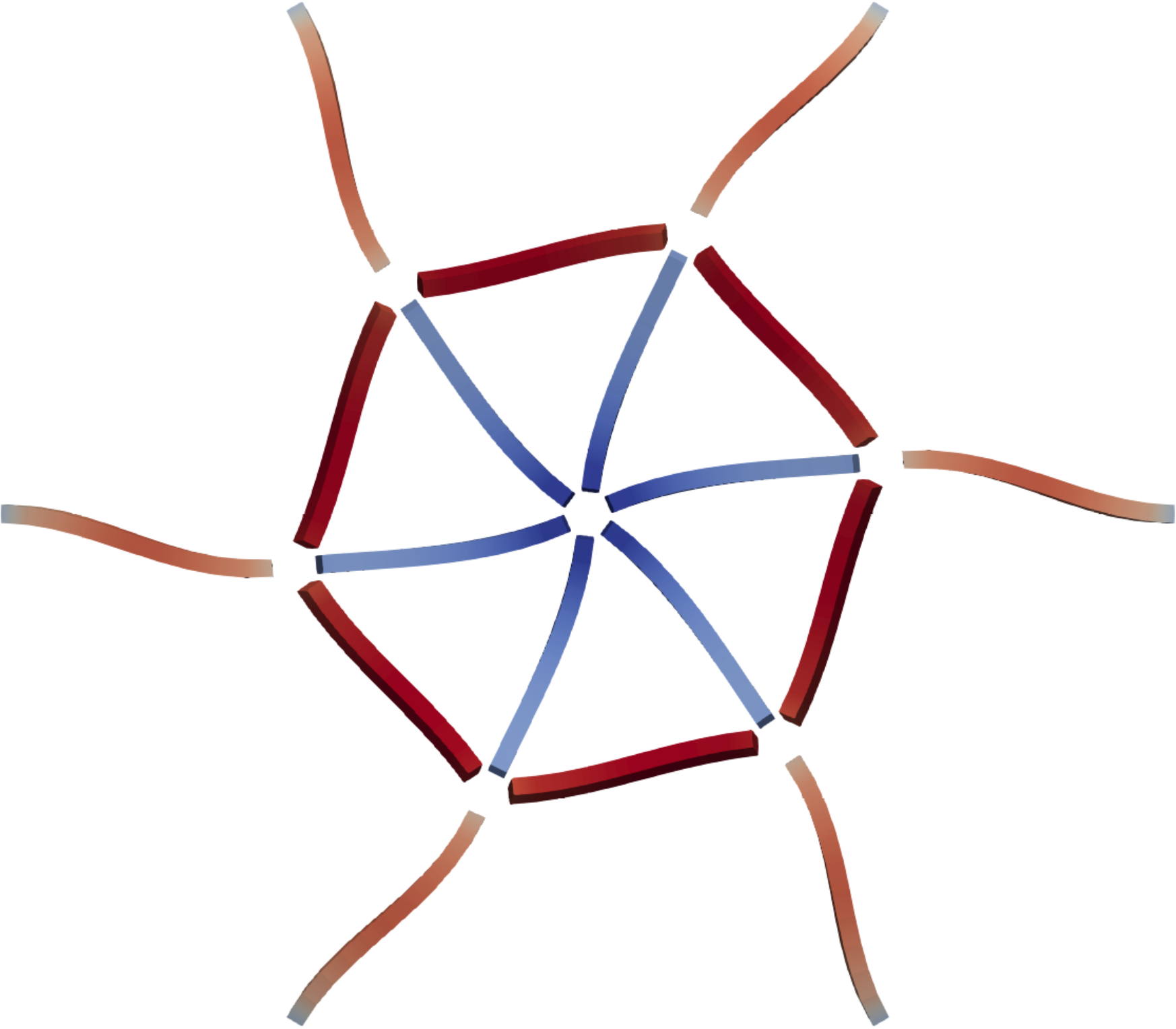}
		\caption{Brick}
		\label{shallow_case2_deformed_compare_brick}			
	\end{subfigure}
	\begin{subfigure}[b]{0.325\textwidth}\centering
		\includegraphics[width=\linewidth]{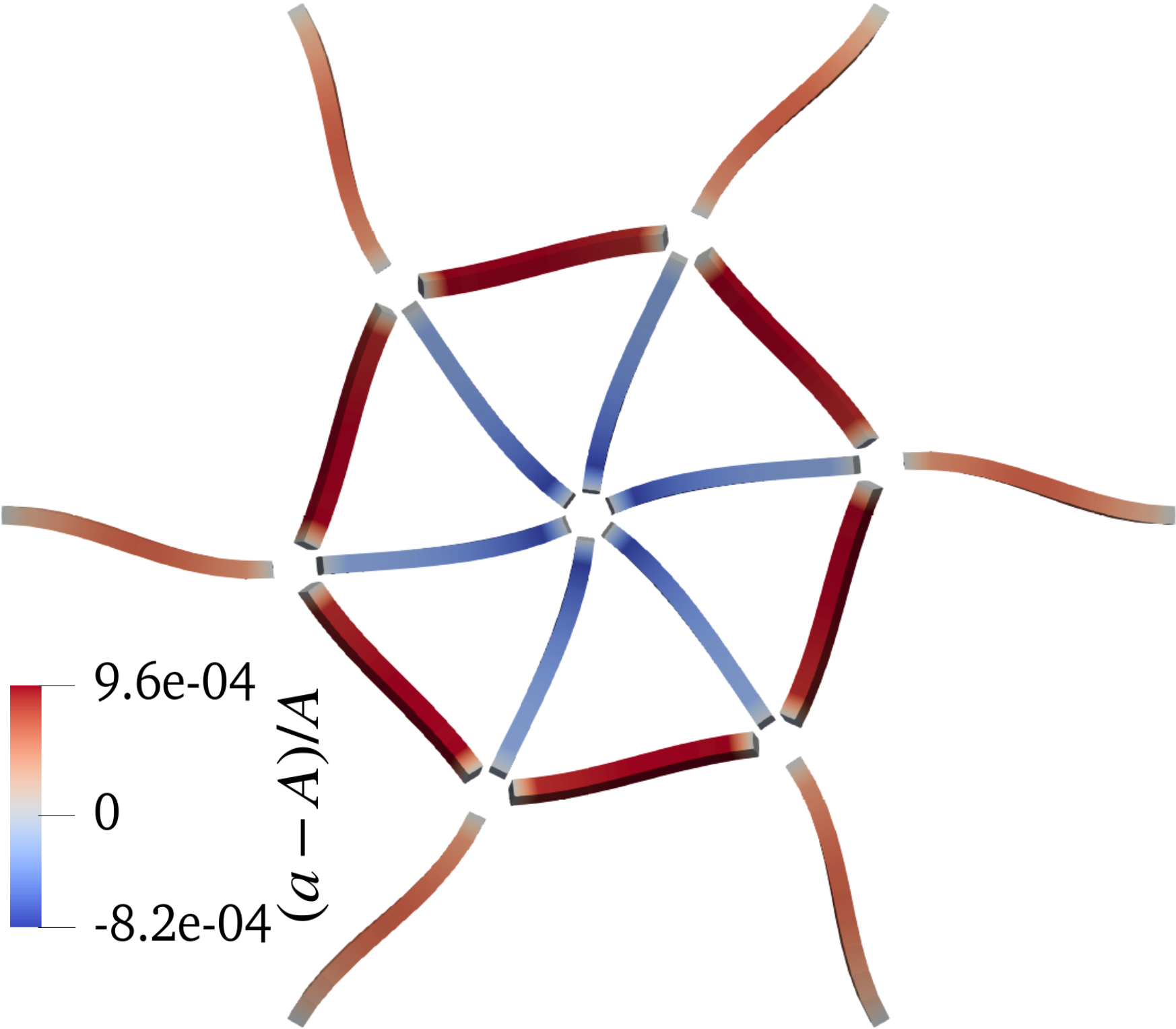}
		\caption{Brick (fully clamped)}
		\label{shallow_case2_deformed_compare_brick_full}			
	\end{subfigure}    
	\caption{Shallow dome (Case 2: added moment load). Comparison of the deformed configurations at the prescribed displacement $\bar w/2$ and moment $M/2$, between the beam and brick solutions (top view). The colors represent the relative change in the cross-sectional area. In (c), the added volume patches are not visualized.}
	\label{shallow_case2_deformed_compare}	
\end{figure}

\begin{figure}[H]
	\centering
	\begin{subfigure}[b]{0.475\textwidth}\centering
		\includegraphics[width=\linewidth]{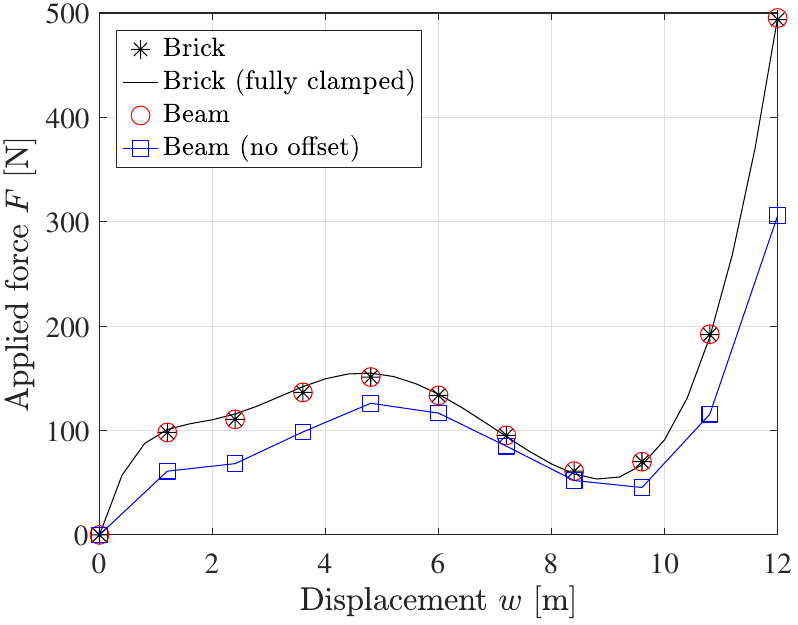}
		\caption{Displacement-load only (Case 1)}
	\end{subfigure}        
	%
	\begin{subfigure}[b]{0.475\textwidth}\centering
		\includegraphics[width=\linewidth]{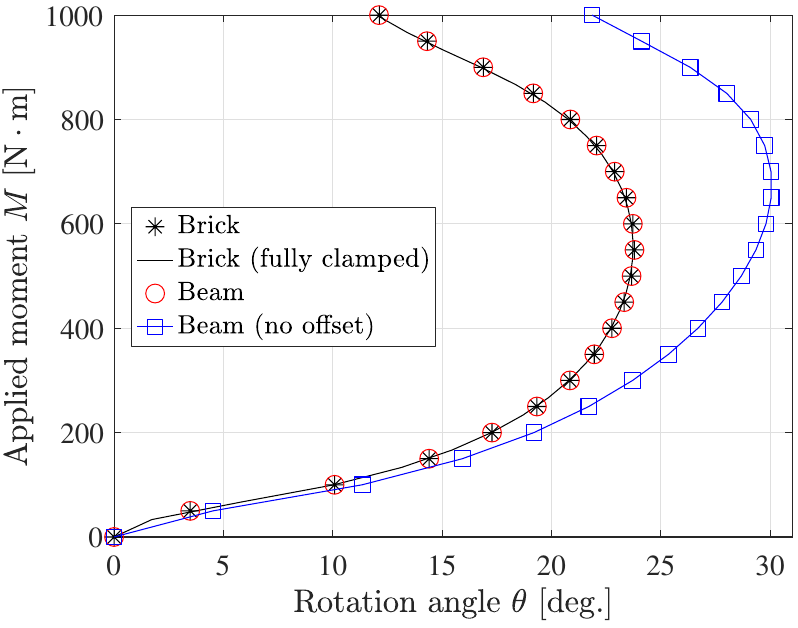}
		\caption{Additional moment-load (Case 2)}
	\end{subfigure}        
	\caption{Shallow dome. Comparison of the equilibrium paths from using beam and brick formulations. For the beam solutions, we have used $p=1$ and $n_\mathrm{el}=10$ in the axis, and $q_\mathrm{a}=1$, $n^\mathrm{a}_\mathrm{el}=8\times8$ in the cross-section. Table \ref{tab_shallow_dome_dof_info} compares the degrees-of-freedom.}
    \label{shallow_force_displacement_curve}	
\end{figure}

\begin{figure}[H]
	\centering
	\begin{subfigure}[b]{0.475\textwidth}\centering
		\includegraphics[width=\linewidth]{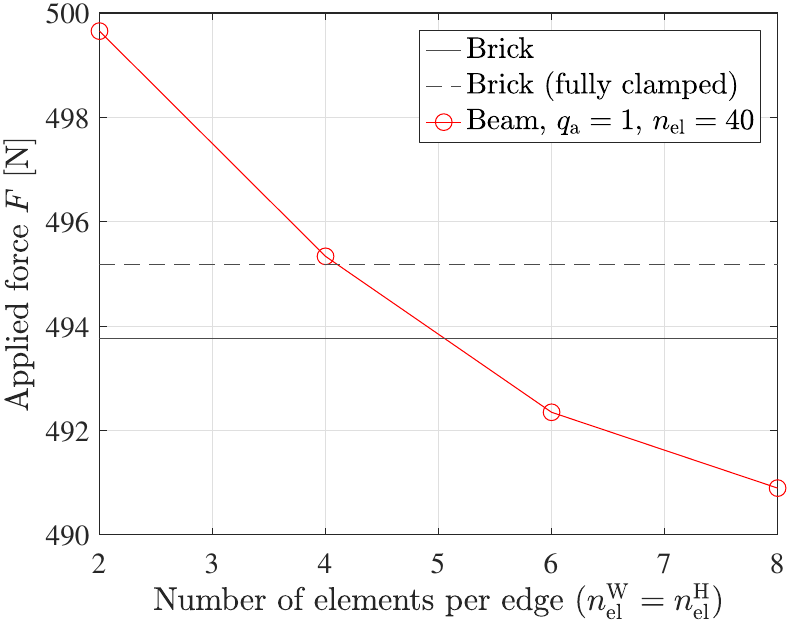}
		\caption{Refinement in the cross-section}
        \label{shallow_deformed_conv_reaction_force_wrt_cs_nelem_cs_case1}			
	\end{subfigure}        
	\begin{subfigure}[b]{0.475\textwidth}\centering
		\includegraphics[width=\linewidth]{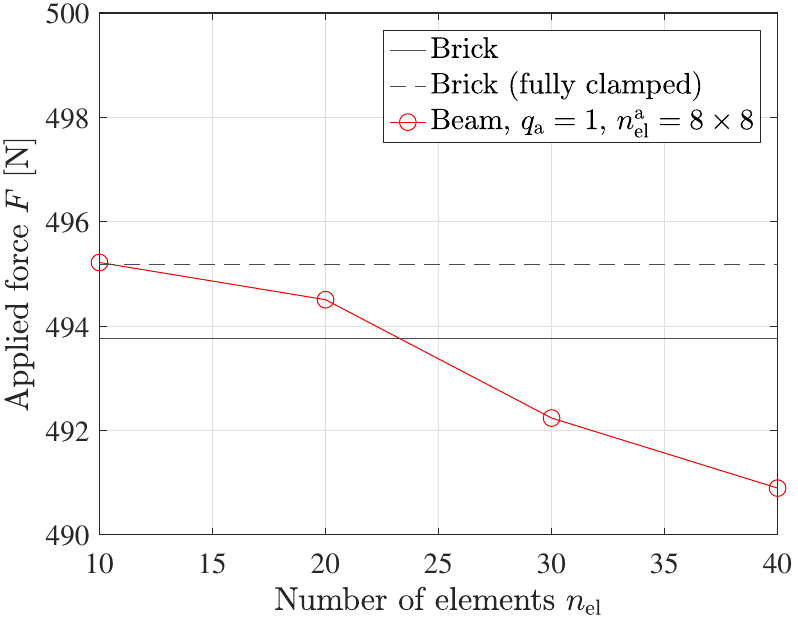}
		\caption{Refinement in the axis}
		\label{shallow_deformed_conv_reaction_force_wrt_nelem_ax_case1}			
	\end{subfigure}            
	\caption{Shallow dome (Case 1: displacement-load only). \textcolor{blue}{Convergence of the applied force $F$ at the final deformed configurations, as the number of elements increases. Table\,\ref{tab_shallow_dome_dof_info} compares the degrees-of-freedom.}}	
	\label{shallow_deformed_conv_wrt_nelem}
\end{figure}

\begin{figure}[H]
	\centering
	\begin{subfigure}[b]{0.475\textwidth}\centering
		\includegraphics[width=\linewidth]{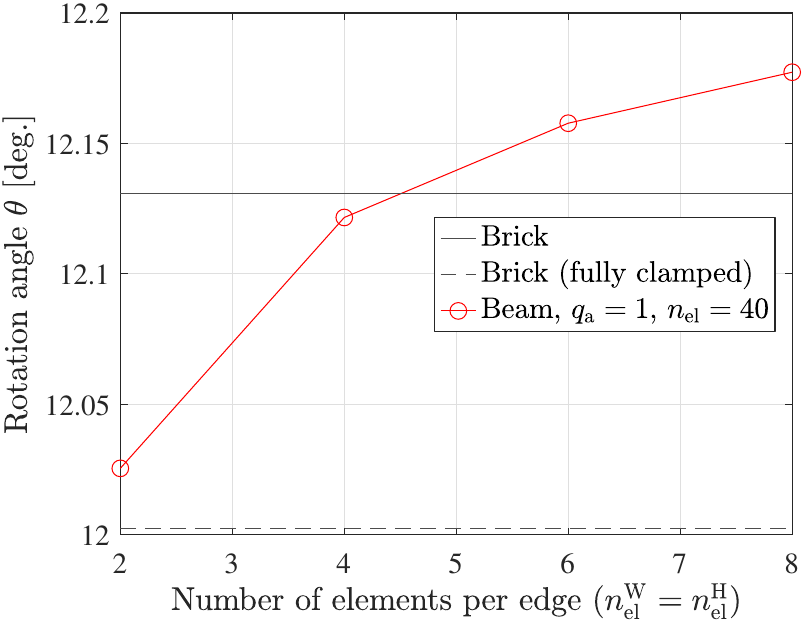}
		\caption{Refinement in the cross-section}
        \label{shallow_deformed_conv_reaction_force_wrt_cs_nelem_cs_case2}			
	\end{subfigure}        
	\begin{subfigure}[b]{0.475\textwidth}\centering
		\includegraphics[width=\linewidth]{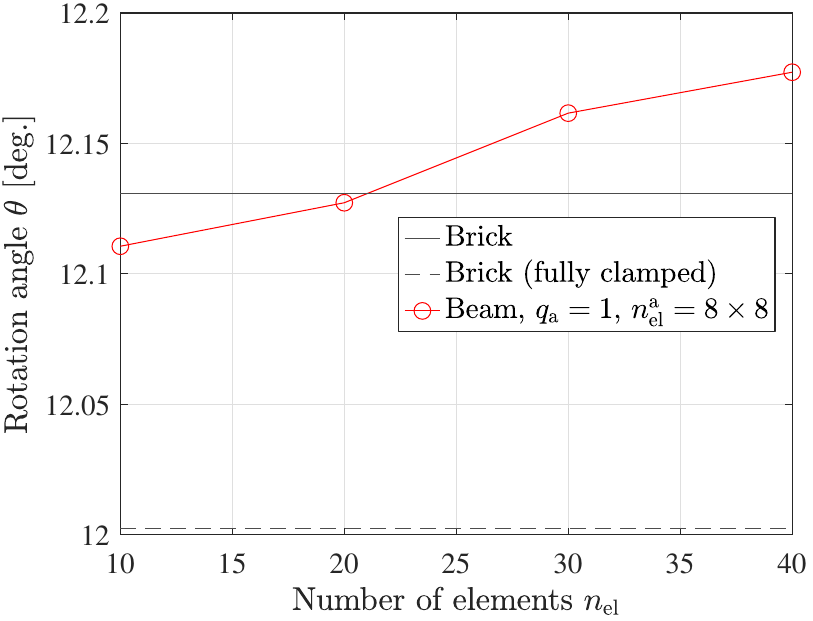}
		\caption{Refinement in the axis}
		\label{shallow_deformed_conv_reaction_force_wrt_nelem_ax_case2}			
	\end{subfigure}            
	\caption{Shallow dome (Case 2: added moment load). \textcolor{blue}{Convergence of the rotation angle $\theta$ at the final deformed configurations, as the number of elements increases. Table\,\ref{tab_shallow_dome_dof_info} compares the degrees-of-freedom.}}	
	\label{shallow_deformed_conv_wrt_nelem_case2}
\end{figure}

\begin{figure}[H]
	\centering
	\begin{subfigure}[b]{0.65\textwidth}\centering
		\includegraphics[width=\linewidth]{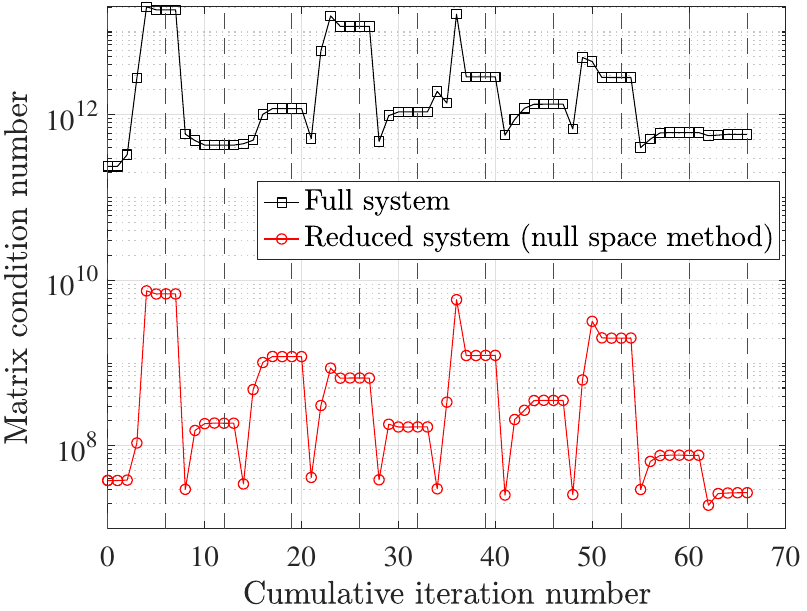}
		\label{shallow_condition number}			
	\end{subfigure}
	\caption{Shallow dome: Comparison of the condition number of the global tangent stiffness matrix in Case 1 (displacement-load only). Here, we have used $p=3$, $p_\mathrm{a}=1$, and $n_\mathrm{el}=40$, $n^\mathrm{a}_\mathrm{el}=2\times2$. \textcolor{blue}{Vertical dashed lines define load increments.}}
	\label{shallow_dome_comparison_condition_number}
\end{figure}

We consider a shallow frame structure consisting of initially straight beams, connected in a dome shape, as shown in Fig.\,\ref{shallow_dome_init_bdc_load_cond}. This example is a benchmark problem that exhibits structural instability (i.e., snap-through behavior); see \citet{battini2002co} and \citet{wackerfuss2009mixed}, from which we adopt the geometrical and material parameters. In conventional modeling of \textit{rigid joint} in such frame structures, it is typical to model the rotational and translational continuity by finite element assembly, in which we may encounter an incorrect description of the rigidity around the joint as well as unphysical overlap. Our objective here is to incorporate a more sophisticated and accurate representation of the joint's stiffness and configuration. To achieve this, we have trimmed the beams' initial axes by spheres at the junctions as illustrated in Fig.\,\ref{shallow_dome_init_bdc_load_cond}, and then coupled the beams using the present offset joint formulation. Here, all offsets are fixed, so that the joints are rigid, and the structure is also clamped at the bottom. \textcolor{blue}{In this example, we apply the same partially clamped joint and boundary conditions for the beam and brick formulations. Further, as mentioned in Remark\,\ref{rem_partial_clamped}, for a reference solution, we have also introduced the following additional treatments:
\begin{itemize}
    \item We fix all degrees-of-freedom of the bottom boundary surfaces.
    \item To constrain the cross-sectional warping at the joints, we have additional volume patches with increased stiffness, as shown in Fig.\,\ref{shallow_dome_init_jct_vol}. 
\end{itemize}
}
\noindent \textcolor{blue}{For the material model, we consider compressible Neo-Hookean-type hyperelasticity, with Young's modulus $E=20690\,\mathrm{Pa}$, and shear modulus $G=8830\,\mathrm{Pa}$. For the fully clamped brick model, we have increased Young's and shear moduli in the added-volume patches by a factor of $100$. Here, we consider two loading cases, as illustrated in Fig.\,\ref{shallow_dome_init_bdc_load_cond},
\begin{itemize}
	\item Case 1: A prescribed vertical displacement $\bar w$ on the top,
	\item Case 2: In addition to the load in Case 1, we further apply a moment $M$ to induce a rotationally symmetric deformation.
\end{itemize}
}
In Figs.\,\ref{shallow_case1_deformed_compare} and \ref{shallow_case2_deformed_compare}, we compare the final deformed configurations between the beam and brick element solutions for Cases 1 and 2, respectively. \textcolor{blue}{In Fig.\,\ref{shallow_case2_deformed_compare}, to clearly show the rotational symmetry, we plot the deformed configurations at half of the total load instead of at the full one.} The colors represent a relative change in the cross-sectional area. \textcolor{blue}{In case of fully clamped brick solutions in Figs.\,\ref{shallow_case1_deformed_compare_brick_full} and \ref{shallow_case2_deformed_compare_brick_full}, it is seen that the cross-sectional areas at the boundaries do not change, in contrast with the partially clamped beam and brick solutions. Note that the deformed configurations are plotted using only the global degrees-of-freedom. Despite this, an excellent agreement is seen between the beam and brick solutions, with the difference being relatively larger in the region near the boundaries, where the cross-sectional warping is more pronounced.}
\noindent \textcolor{blue}{In Fig.\,\ref{shallow_force_displacement_curve}, we compare the equilibrium paths between the beam and brick solutions, where an excellent agreement is also observed. Furthermore, the beam solution without offsets significantly underestimates the stiffness (blue curves), since it incorrectly represents the rigidity at the joint.}
\noindent \textcolor{blue}{In Figs.\,\ref{shallow_deformed_conv_wrt_nelem} and \ref{shallow_deformed_conv_wrt_nelem_case2}, we investigate the convergence of the beam solution, as we refine the mesh in the cross-section and axis. Since the fully clamped model has a larger stiffness, the applied force is larger, but the rotation angle (displacement) is smaller than that of the partially clamped model. It is seen that the beam solutions converge, but the converged values do not agree perfectly with the brick solutions under the same partially clamped conditions. We have relevant comments on this in \ref{observ_ex_z_torsion}, and further investigation into this remains future work. It is noted that the beam formulation requires a much smaller number of global degrees-of-freedom, as shown in Table \ref{tab_shallow_dome_dof_info}. For the convergence behavior of brick solutions, see Table \ref{tab_app_shallow_dome_brick_conv}.}
\begin{table}[H]
  \centering
  \footnotesize
  \caption{\textcolor{blue}{Shallow dome. Comparison of degrees-of-freedom. Here, the element counts are defined per patch.}}
    \begin{tabular}{lcccccccc}
    \toprule
    \multicolumn{1}{r}{} & \multicolumn{2}{c}{Degrees} &       & \multicolumn{2}{c}{Elements} &       & \multicolumn{2}{c}{Degrees-of-freedom} \\
\cmidrule{2-3}\cmidrule{5-6}\cmidrule{8-9}    \multicolumn{1}{r}{} & L     & \multicolumn{1}{c}{W,H} &       & L     & \multicolumn{1}{c}{W,H} &      & Global & Internal \\
    \midrule
    \multicolumn{1}{l}{Beam} & 3     & 1     &       & 40    & 8    &       & 6791      & 686880 \\
    \multicolumn{1}{l}{Brick} & 3     & 3     &       & 50    & 5     &       & 182993  & $-$ \\
    Brick (fully clamped) & 3     & 3     &       & 50    & 5     &       & {209549}  & $-$ \\
    \bottomrule
    \end{tabular}%
  \label{tab_shallow_dome_dof_info}%
\end{table}%
\noindent In Fig.\,\ref{shallow_dome_comparison_condition_number}, we show that the condition number of the global system matrix can be significantly improved by using the present null space method. 
\setcounter{remark}{0}
\section{Conclusions}
\label{conclusions_jct}
In this paper, we present a variational formulation that constrains or releases the relative rotation and translation of multiple connected beams with arbitrarily shaped cross-sections. Since this formulation requires no explicit interface to the rotational and translational degrees-of-freedom, it applies to a \textit{general} kinematic description of the cross-section, like a polynomial expansion-based one, as well as brick formulations. Here, we have implemented an \textit{implicit} interface to the rigid components of the cross-sectional motion by constraints. This formulation is \textit{unbiased}, that is, no master-slave relationship between connected beams is required. It also enables a straightforward calculation of the reaction forces and moments from the Lagrange multipliers. Further, we have eliminated those multipliers using a discrete null space method, which yields size-reduction and improved conditioning of the system matrix. Due to the localized constraints, the calculation of the fundamental bases is also localized, leading to a sparse system matrix and efficient post-processing. We have applied the present joint formulation to an extensible director-based beam formulation in nonlinear elastostatics and verified the results by comparison to brick solutions, where excellent agreement is observed.

An extension of the present joint formulation may consider the following: (i) In the present work, we have released the interface constraints for the relative translation only, using the variable offsets. One can further supplement the relative rotation by adding a released rotation angle. (ii) It is also straightforward to add friction to the released offset variable for a frictional prismatic joint. \textcolor{blue}{(iii) One can extend the dimensional reduction of the joint formulation by the first-order beam kinematics to incorporate the enhanced cross-sectional strains. Incorporating the enhanced kinematics may further improve the agreement between the beam and brick solutions.} These extensions remain future work. 

\section*{Acknowledgements}
M.-J. Choi was supported by the Deutsche Forschungsgemeinschaft (DFG, German Research Foundation) - Project number 523829370.

\appendix
\renewcommand{\thesection}{\Alph{section}}
\renewcommand{\theremark}{\thesection.\arabic{remark}}

\setcounter{remark}{0}
\section{\textcolor{blue}{Beam formulation}}
\label{app_beam_form}
\textcolor{blue}{In the first-order beam kinematics, the position of a material point in the current configuration can be parameterized by
\begin{align}
    \label{bkin_first_cur}
    {\boldsymbol{x}}_{\bar s} = {\boldsymbol{\varphi }}(\bar s) + {\zeta ^\alpha}{{\boldsymbol{d}}_\alpha}(\bar s),
\end{align}
using the position on the \textit{current axis}, $\boldsymbol{\varphi}\in\Bbb{R}^3$, and the \textit{current directors}, $\boldsymbol{d}_\alpha\in\Bbb{R}^3$. For a compact notation, we introduce the operator
\begin{align}
    \label{app_basis_Pi}
    {\boldsymbol{\Pi }} \coloneqq {\left[ {\begin{array}{*{20}{c}}
{{{\boldsymbol{1}}_3}}&{{\zeta ^1}{{\boldsymbol{1}}_3}}&{{\zeta ^2}{{\boldsymbol{1}}_3}}
\end{array}} \right]^{\rm{T}}},
\end{align}
such that
\begin{align}
    \label{x_sbar_oper_Pi}
    {\boldsymbol{x}} = {{\boldsymbol{\Pi }}^{\rm{T}}}{{\boldsymbol{y}}},
\end{align}
with the so-called \textit{configuration variable} $\boldsymbol{y}\coloneqq\left[\boldsymbol{\varphi}^\mathrm{T},\boldsymbol{d}_1^{\,\mathrm{T}},\boldsymbol{d}_2^{\,\mathrm{T}} \right]^{\mathrm{T}}\in\Bbb{R}^d$. Here $d=9$ denotes the number of kinematic degrees-of-freedom of the cross-section. For a subsequent finite element formulation for the present beam kinematics, we refer readers to \citet{choi2021isogeometric}. The assumed first-order kinematics in Eq.\,(\ref{bkin_first_cur}) may not be sufficient to represent the beam's correct stiffness, which makes it necessary to have additional degrees-of-freedom to account for cross-sectional deformations. This will be discussed below.}
\textcolor{blue}{
\subsection{Construction of warping basis functions}}
\label{app_construct_warp_basis}
\textcolor{blue}{
To enrich the cross-sectional warping modes, in Eq.\,(\ref{x_sbar_oper_Pi}), one can employ additional higher-order polynomial bases in Eq.\,(\ref{app_basis_Pi}) with unknown coefficients; see, e.g., \citet{moustacas2019enrichissement} and \citet{choi2022isogeometric} for relevant works. However, this approach increases the number of global DOFs and may also lead to an ill-conditioned tangent stiffness matrix due to different dimensions of those unknown coefficients. These can be resolved by introducing an element-wise static condensation of the unknown coefficients; see \citet{wackerfuss2009mixed}, where the enhancement is applied at the strain level using a so-called enhanced assumed strain (EAS) method, with properly designed coupled polynomial bases for the additional strains in the cross-section. Here, we have the total Green-Lagrange strain tensor, decomposed into \citep{buchter1994three,betsch1996four,bischoff1997shear}
\begin{align}
    \boldsymbol{E}^\mathrm{tot} = \underbrace{\boldsymbol{E}(\boldsymbol{y})}_\text{compatible} + \underbrace{\widetilde{\boldsymbol{E}}}_\text{enhanced},
\end{align}
where the detailed expressions of the compatible part can be found in \citet{choi2021isogeometric}. Here, we recall only the expression of the longitudinal strain component.
\begin{align}
    \label{app_recall_eq_longi_strn}
    E_{33}=\varepsilon + \zeta^\alpha\rho_\alpha + \zeta^\alpha\zeta^\beta\kappa_{\alpha\beta},
\end{align}
where $\varepsilon$ denotes the axial strain, $\rho_\alpha$ denotes the change of curvature due to bending, and $\kappa_{\alpha\beta}$ accounts for the torsion-induced (quadratic) longitudinal strain, the so-called Wagner effect \citep{Wagner1929}. For a more detailed expression of the relation between the beam strains and the configuration variable $\boldsymbol{y}$, we refer readers to \citet{choi2021isogeometric}. In \citet{wackerfuss2009mixed,wackerfuss2011nonlinear}, it turns out that having the quadratic terms in Eq.\,(\ref{app_recall_eq_longi_strn}) is sufficient to represent the beam's correct torsional stiffness. This observation leads to the omission of the enhancement of the longitudinal strain $E_{33}$ in the subsequent EAS approach. This approach is called \textit{E2-model}, which is more efficient than having the enhancement of ${\widetilde E}_{33}$ (so-called \textit{E1-model}), since the number of enhanced strain parameters is smaller. We employ this method in the present work. Here, the enhanced strain tensor can be represented by the parameters $\boldsymbol{\alpha}\in\Bbb{R}^{d_\mathrm{a}}$, as
\begin{equation}
	\label{phy_enh_dec}
	{\underline{\boldsymbol{\widetilde E}}} = {\boldsymbol{\Gamma }}({\zeta ^1},{\zeta ^2})\,{\boldsymbol{\alpha }}(s),
\end{equation}
where $\underline{(\bullet)}$ denotes the Voigt notation, and we have also introduced the following matrix of basis functions \citep{wackerfuss2009mixed, wackerfuss2011nonlinear}
\begin{align}
	\label{wg_basis_matrix_gamma}
	{\boldsymbol{\Gamma }}\!\!\:\left( {{\zeta ^1},{\zeta ^2}} \right) = \left[ {\begin{array}{*{20}{c}}
			{{{\bf{w}}_1}}&{}&{}&{}&{}\\
			{}&{{{\bf{w}}_1}}&{}&{}&{}\\
			{}&{}&\cancel{{{\bf{w}}_3}}&{}&{}\\
			{}&{}&{}&{{{\bf{w}}_2}}&{}\\
			{}&{}&{}&{}&{{\bf{w}}_{4,1}}\\
			{}&{}&{}&{}&{{{\bf{w}}_{4,2}}}
	\end{array}} \right]_{\textcolor{black}{6\times{d_\mathrm{a}}}}.
\end{align}
Here the blank entries represent zeros, and $(\bullet)_{,\alpha}$ denotes the derivative with respect to $\zeta^\alpha$ ($\alpha=1,2$). Here, ${\bf w}_i$ $(i=1,2,3,4)$ denote row arrays of the basis functions, where ${\bf{w}}_1$, ${\bf{w}}_2$, ${\bf{w}}_3$, and ${\bf{w}}_4$ enrich the transverse normal (${\widetilde E}_{11}$, ${\widetilde E}_{22}$), in-plane shear (${\widetilde E}_{12}$), longitudinal (${\widetilde E}_{33}$), and transverse shear (${\widetilde E}_{13}$, ${\widetilde E}_{23}$) strains, respectively. Note again that the enrichment of the longitudinal strain component ${\widetilde E}_{33}$ is dropped (i.e., the columns for ${\bf{w}}_3$ are eliminated), which improves the efficiency due to a reduced number of the enhanced strain parameters ($d_\mathrm{a}$). It is noted that the last two rows account for the coupling between ${\widetilde E}_{13}$ and ${\widetilde E}_{23}$ \citep{gruttmann2000theory}. Further, we have the following conditions for the orthogonality to the existing stress fields from the kinematics,
\begin{subequations}
    \label{app_ortho_for_w}
    \begin{align}
        \label{app_ortho_w1_w2}
        \int_\mathcal{A} {{{\bf{w}}_1}\,{\rm{d}}A}  = \int_\mathcal{A} {{{\bf{w}}_2}\,{\rm{d}}A}  ={\bf{0}},
    \end{align}
    \begin{align}
        \label{app_ortho_w4}    
        \int_\mathcal{A} {{{\bf{w}}_4}\,{\rm{d}}A}  = \int_\mathcal{A} {{\zeta ^\alpha }{{\bf{w}}_4}\,{\rm{d}}A}  = {\bf{0}},\,\,\alpha=1,2.
    \end{align}
\end{subequations}
\begin{remark}\small
Here, due to the additional three strain modes from the kinematics, which are illustrated in Fig.\,\ref{beam_cs_strn_mode_allow}, we also need to apply the orthogonality condition for ${\bf w}_1$ and ${\bf w}_2$ in Eq.\,(\ref{app_ortho_w1_w2}). This is an additional requirement compared with the descriptions in previous works \citep{wackerfuss2009mixed,wackerfuss2011nonlinear}, in which the directors are assumed to be orthonormal in the kinematics; see also \citet{choi2024objective} for the relevant discussions.
\end{remark}
}
\subsubsection{Corrected basis functions for orthogonality}
\label{app_warp_cor_local_shape}
\textcolor{blue}{
We introduce the so-called \textit{local concept} of enhancement to construct the bases ${\bf w}_i$ in Eq.\,(\ref{wg_basis_matrix_gamma}), using finite element basis functions instead of polynomial expansions (the so-called \textit{global concept}). The local concept turns out to be more effective for representing cross-sectional warping in open cross-sections; see \citet{wackerfuss2011nonlinear} for the relevant discussions. Here, we use the following basis functions 
\begin{subequations}
\label{app_row_arr_w_fe_shape}
\begin{align}
{{\bf{w}}_1} = {{\bf{w}}_2} = \left[ {{{\hat H}_1},{{\hat H}_2}, \cdots,{{\hat H}_{m^\mathrm{a}_\mathrm{cp}}}} \right] \in {\Bbb{R}^{1 \times m^\mathrm{a}_\mathrm{cp}}},
\end{align}
\begin{align}
{{\bf{w}}_4} = \left[ {{{\tilde H}_1},{{\tilde H}_2}, \cdots,{{\tilde H}_{m^\mathrm{a}_\mathrm{cp}}}} \right] \in {\Bbb{R}^{1 \times m^\mathrm{a}_\mathrm{cp}}},
\end{align}
\end{subequations}
constructed by using the \textit{corrected} finite element shape functions
\begin{subequations}  
\begin{align}
    {\hat H_I} &= {N^{q_\mathrm{a}}_I} + \alpha _I^{(0,0)},\label{app_warp_cor_w1_w2}\\
    {\tilde H_I} &= {N^{q_\mathrm{a}}_I} + \underbrace{\beta _I^{(0,0)} + \beta _I^{(1,0)}{\zeta ^1} + \beta _I^{(0,1)}{\zeta ^2}}_\text{correction},\label{app_warp_cor_w4}
\end{align}
\end{subequations}
for $I=1,\cdots,m^\mathrm{a}_\mathrm{cp}$, where we use the same number of shape functions, denoted by $m^\mathrm{a}_\mathrm{cp}$, for all ${\bf w}_i$ ($i=1,2,4$), which yields $d_\mathrm{a}=4\,m^\mathrm{a}_\mathrm{cp}$. Here, $N_I$ denotes the finite element shape functions, for which we utilize NURBS basis functions of degree $q_\mathrm{a}$ in the framework of IGA, which enables representing the exact geometry of the circular fillet shown in Fig.\,\ref{z_shape_with_fillet}. The constant correction coefficients $\alpha^{(0,0)}_I$, $\beta_I^{(0,0)}$, $\beta_I^{(1,0)}$, and $\beta_I^{(0,1)}\in\Bbb{R}$ can be uniquely determined from the orthogonality conditions in Eqs.\,(\ref{app_ortho_w1_w2}) and (\ref{app_ortho_w4}). This correction from the orthogonality conditions introduces redundancy in the enhanced strain parameters, which can be interpreted as rigid-body motions, and must be eliminated to ensure the uniqueness of the solution. Therefore, we eliminate one basis function (constant field) from ${\bf w}_1$ and ${\bf w}_2$ and three basis functions (one constant and two linear fields along $\zeta^1$ and $\zeta^2$) from ${\bf w}_4$; see \citet[Fig.\,3]{wackerfuss2011nonlinear} for an illustration. A detailed description of the subsequent finite element formulation of the EAS method can be found in \citet{choi2021isogeometric}. 
\begin{remark} 
\label{app_remark_eas_notation}\small
We define the following notations for the finite element discretization of the enhanced strain field in the cross-section: 
    \begin{itemize}
        \item $m^\mathrm{a}_\mathrm{cp}$ denotes the total number of control points in the cross-section, which determines the total number of enhanced strain parameters in the cross-section, $d_\mathrm{a}=4\,m^\mathrm{a}_\mathrm{cp}$.
        \item $n^\mathrm{a}_\mathrm{el}$ denotes the number of elements in the cross-section.
        \item $q_\mathrm{a}$ denotes the degree of basis functions for the cross-sectional discretization.
    \end{itemize}
\end{remark}
\noindent Further, for the discretization of the enhanced strain parameters in each element along the axis, Lagrange polynomials of degree $p_\mathrm{a}$, with inter-element discontinuity allowed. This enables the element-wise static condensation of the unknown coefficients of the enhanced strain field. That is, the number $d_\mathrm{a}$ does not increase the global decrees-of-freedom \citep{wackerfuss2009mixed}, but is associated with the number of internal variables. Further details on the static condensation and subsequent solution update process can be found in \citet{choi2021isogeometric}.}
\setcounter{remark}{0}
\section{Constraint formulation}
\subsection{Constraint Jacobian}
\label{app_cnst_jcb_general}
In the following, we present the detailed expressions of the constraint Jacobians, $\boldsymbol{G}_\mathrm{x}$, $\boldsymbol{G}_\mathrm{s}$, $\boldsymbol{G}_{\mathrm{\bar q}}$, and $\boldsymbol{G}_\Bbb{\bar q}$. We first recall
\begin{subequations}
    \label{app_cnst_fderiv}
    \begin{align}
        \label{app_cnst_jcb_Gx}
        {{\boldsymbol{G}}_{\rm{x}}} \coloneqq \dfrac{\partial \Bbb{c}}{\partial\boldsymbol{x}_{\bar s}}= \left[ {\renewcommand{\arraystretch}{1.5}\begin{array}{*{20}{c}}
    {{{\partial \boldsymbol{\phi}}}/{{\partial {{\boldsymbol{x}}_{{{\bar s}}}}}}}\\
    {{{\partial {\boldsymbol{\psi }}}}/{{\partial {{\boldsymbol{x}}_{{{\bar s}}}}}}}
    \end{array}} \right],
    \end{align}
    \begin{align}
        \label{app_cnst_jcb_Gs}    
        {{\boldsymbol{G}}_{\rm{s}}} \coloneqq \dfrac{\partial \Bbb{c}}{\partial\boldsymbol{s}} = \left[ {\renewcommand{\arraystretch}{1.5}\begin{array}{*{20}{c}}
    {{{\partial \boldsymbol{\phi}}}/{{\partial {s^1}}}}&{{{\partial \boldsymbol{\phi} }}/{{\partial {s^2}}}}&{{{\partial \boldsymbol{\phi}}}/{{\partial {s^3}}}}\\
    {{{\partial {\boldsymbol{\psi }}}}/{{\partial {s^1}}}}&{{{\partial {\boldsymbol{\psi }}}}/{{\partial {s^2}}}}&{{{\partial {\boldsymbol{\psi }}}}/{{\partial {s^3}}}}
    \end{array}} \right],
    \end{align}
    \begin{align}
        \label{app_cnst_jcb_Gqbar}        
        {{\boldsymbol{G}}_{{\rm{\bar q}}}} \coloneqq \dfrac{\partial \Bbb{c}}{\partial\boldsymbol{\bar q}}=\left[ {\renewcommand{\arraystretch}{1.5}\begin{array}{*{20}{c}}
    {\partial \boldsymbol{\phi}/\partial {\boldsymbol{\bar \varphi }}}&{\partial 
    \boldsymbol{\phi}/\partial {{{\boldsymbol{\bar t}}}_1}}&{\partial \boldsymbol{\phi}/\partial {{{\boldsymbol{\bar t}}}_2}}&{\partial \boldsymbol{\phi}/\partial {{{\boldsymbol{\bar t}}}_3}}\\
    {\partial {\boldsymbol{\psi }}/\partial {\boldsymbol{\bar \varphi }}}&{\partial {\boldsymbol{\psi }}/\partial {{{\boldsymbol{\bar t}}}_1}}&{\partial {\boldsymbol{\psi }}/\partial {{{\boldsymbol{\bar t}}}_2}}&{\partial {\boldsymbol{\psi }}/\partial {\boldsymbol{\bar t}_3}}
    \end{array}} \right],
    \end{align}
    and 
    \begin{align}
        \label{app_cnst_jcb_Gqqbar}            
        {{\boldsymbol{G}}_{{\Bbb{\bar q}}}} \coloneqq \dfrac{\partial \Bbb{c}}{\partial\Bbb{\bar q}} = \left[ {\renewcommand{\arraystretch}{1.5}\begin{array}{*{20}{c}}
    {{{\partial \boldsymbol{\phi}}}/{{\partial {\boldsymbol{\bar \varphi }}}}}&{{{\partial \boldsymbol{\phi}}}/{{\partial {\boldsymbol{\bar \theta }}}}}\\
    {{{\partial {\boldsymbol{\psi }}}}/{{\partial {\boldsymbol{\bar \varphi }}}}}&{{{\partial {\boldsymbol{\psi }}}}/{{\partial {\boldsymbol{\bar \theta }}}}}
    \end{array}} \right].
    \end{align}
\end{subequations}
Here, those first-order derivatives in Eq.\,(\ref{app_cnst_fderiv}) are given, as follows: First, in Eq.\,(\ref{app_cnst_jcb_Gx}) we have
\begin{subequations}
    \begin{align}
        \dfrac{{\partial \boldsymbol{\phi}}}{{\partial {{\boldsymbol{x}}_{\bar s}}}} = \left(\int_\mathcal{A} (\bullet)\,{{\rm{d}}A}\right) \,{\boldsymbol{1}},
    \end{align}
    \begin{align}
        \dfrac{{\partial \boldsymbol{\psi}}}{{\partial {{\boldsymbol{x}}_{\bar s}}}} = {\left[ {{\boldsymbol{t}}_\alpha } \right]_ \times }\left({\int_\mathcal{A} {{\zeta ^\alpha }(\bullet)\,{\rm{d}}A}}\right).
    \end{align}
\end{subequations}
Second, in Eq.\,(\ref{app_cnst_jcb_Gs}) we have
\begin{subequations}
    \begin{align}
        \label{app_deriv_phi_wrt_si}
        \dfrac{{\partial \boldsymbol{\phi}}}{{\partial {s^i}}} = A{{\boldsymbol{t}}_i},
    \end{align}
    and
    \begin{align}
        \label{app_deriv_psi_wrt_si}    
        \dfrac{{\partial \boldsymbol{\psi}}}{{\partial {s^i}}} = {I^\alpha }{{\boldsymbol{t}}_\alpha } \times {{\boldsymbol{t}}_i}.
    \end{align}
\end{subequations}
Third, in Eq.\,(\ref{app_cnst_jcb_Gqbar}) we have
\begin{subequations}
    \begin{align}
        \dfrac{{\partial \boldsymbol{\phi} }}{{\partial {\boldsymbol{\bar \varphi }}}} =  - A{\boldsymbol{1}},
    \end{align}
    \begin{align}
        \dfrac{{\partial \boldsymbol{\psi} }}{{\partial {\boldsymbol{\bar \varphi }}}} =  - {\left[ {{I^\alpha }{{{\boldsymbol{\bar t}}}_\alpha }} \right]_ \times },
    \end{align}    
\end{subequations}
and for $i=1,2,3$,
\begin{subequations}
\begin{align}
     \dfrac{{\partial \boldsymbol{\phi}}}{{\partial {{{\boldsymbol{\bar t}}}_i}}} =  {\dfrac{{\partial \boldsymbol{\phi} }}{{\partial {{\boldsymbol{t}}_i}}}{{\boldsymbol{R}}_0}},
\end{align}
\begin{align}
    \dfrac{{\partial \boldsymbol{\psi}}}{{\partial {{{\boldsymbol{\bar t}}}_i}}} =  {\dfrac{{\partial \boldsymbol{\psi} }}{{\partial {{\boldsymbol{t}}_i}}}{{\boldsymbol{R}}_0}},
\end{align}
\end{subequations}
with
\begin{subequations}
    \begin{align}
        {\dfrac{{\partial \boldsymbol{\phi}}}{{\partial {{\boldsymbol{t}}_i}}}}={\left( {A{s^i} - \delta _\alpha ^i{I^\alpha }} \right){\boldsymbol{1}}},
    \end{align}
    \begin{align}
        \dfrac{{\partial \boldsymbol{\psi}}}{{\partial {{\boldsymbol{t}}_i}}} ={\left[ {\delta _\alpha ^i\left( {{I^\alpha }{\boldsymbol{\bar \varphi }} - \int_\mathcal{A} {{\zeta ^\alpha }{{\boldsymbol{x}}_{\bar s}}\,{\rm{d}}A} } \right)} \right]_\times}.
    \end{align}
\end{subequations}
Fourth, in Eq.\,(\ref{app_cnst_jcb_Gqqbar}) we have
\begin{subequations}
    \begin{align}
        \dfrac{{\partial \boldsymbol{\phi} }}{{\partial {\boldsymbol{\bar \theta }}}} = -\dfrac{{\partial \boldsymbol{\phi}}}{{\partial {{{\boldsymbol{\bar t}}}_i}}}\left[{{\boldsymbol{\bar t}}_i}\right]_\times,
    \end{align}
    \begin{align}
        \dfrac{{\partial \boldsymbol{\psi} }}{{\partial {\boldsymbol{\bar \theta }}}} = -\dfrac{{\partial \boldsymbol{\psi} }}{{\partial {{{\boldsymbol{\bar t}}}_i}}}\left[{{\boldsymbol{\bar t}}_i}\right]_\times.
    \end{align}    
\end{subequations}
\subsection{Tangent operator}
\label{app_tangent_oper}
Here, we present the detailed expressions of the submatrices $\boldsymbol{k}_{\mathrm{x}\Bbb{\bar q}}$, $\boldsymbol{k}_{\mathrm{s}\Bbb{\bar q}}$, and $\boldsymbol{k}_{\Bbb{\bar q}\Bbb{\bar q}}$ of the tangent operator in Eq.\,(\ref{del_J_tangent_k}). We first recall
\begin{align}
    \label{app_k_xqbar}
    \boldsymbol{k}_{\mathrm{x}\Bbb{\bar q}} = \left[ {\begin{array}{*{20}{c}}
{\boldsymbol{0}_{3\times3}}&{{{\boldsymbol{k}}_{{\rm{x\bar \theta }}}}}
\end{array}} \right],
\end{align}
\begin{align}
    \label{app_k_sqbar}
    {{\boldsymbol{k}}_{{\rm{s}\Bbb{\bar q}}}} = \left[ {\begin{array}{*{20}{c}}
\boldsymbol{0}_{3\times3}&{{{\boldsymbol{k}}_{{\rm{s\bar \theta }}}}}
\end{array}} \right],
\end{align}
and
\begin{align}
    \label{app_k_qqbar}
    \boldsymbol{k}_{\Bbb{\bar q}\Bbb{\bar q}} =\left[ {\renewcommand{\arraystretch}{1.5}\begin{array}{*{20}{c}}
{\bf{0}}_{3\times{3}}&{{{\boldsymbol{k}}_{{\rm{\bar \varphi \bar \theta }}}}}\\
{{\boldsymbol{k}}_{{\rm{\bar \varphi \bar \theta }}}^{\rm{T}}}&{{{\boldsymbol{k}}_{{\rm{\bar \theta \bar \theta }}}}}
\end{array}} \right].
\end{align}
\noindent First, in Eq.\,(\ref{app_k_xqbar}) we have
\begin{align}
    \label{app_k_xtht}
    {{\boldsymbol{k}}_{{\rm{x\bar \theta }}}} = {\left( {\dfrac{{{\partial ^2}{\Bbb{c}}}}{{\partial {{\boldsymbol{x}}_{\bar s}}\partial {\boldsymbol{\bar \theta }}}}} \right)^{\!\rm{T}}}{\Bbb{f}} = - \left( {\int_\mathcal{A} {{\zeta ^\alpha }\left(\bullet\right){\rm{d}}A} }\right){\left[ {\boldsymbol{\mu }} \right]_ \times }{\left[ {{{\boldsymbol{t}}_\alpha }} \right]_ \times }{{\boldsymbol{R}}_0}.
\end{align}
\noindent Second, in Eq.\,(\ref{app_k_sqbar}) we have
\begin{align}
    {{\boldsymbol{k}}_{{\rm{s\bar \theta }}}} \coloneqq {\left( {\dfrac{{{\partial ^2}{\Bbb{c}}}}{{\partial {{\boldsymbol{s}}}\partial {\boldsymbol{\bar \theta }}}}} \right)^{\!\rm{T}}}{\Bbb{f}} = \left[ {\begin{array}{*{20}{c}}
{{{\boldsymbol{s}}_1}}&{{{\boldsymbol{s}}_2}}&{{{\boldsymbol{s}}_3}}
\end{array}} \right]^{\!\mathrm{T}}{{\boldsymbol{R}}_0},
\end{align}
with
\begin{align}
    {{\boldsymbol{s}}_i} &\coloneqq \dfrac{{\partial \boldsymbol{\phi}}}{{\partial {s^i}}} \times {\boldsymbol{\lambda }} + \dfrac{{\partial \boldsymbol{\psi}}}{{\partial {s^i}}} \times {\boldsymbol{\mu }} \nonumber\\
    &=A{{\boldsymbol{t}}_i}\times\boldsymbol{\lambda} + \left({I^\alpha }{{\boldsymbol{t}}_\alpha } \times {{\boldsymbol{t}}_i}\right)\times\boldsymbol{\mu},
\end{align}
from using Eqs.\,(\ref{app_deriv_phi_wrt_si}) and (\ref{app_deriv_psi_wrt_si}).
Third, in Eq.\,(\ref{app_k_qqbar}) we have 
\begin{align}
    \label{k_phi_tht}
    {{\boldsymbol{k}}_{{\rm{\bar \varphi \bar \theta }}}} \coloneqq {\left( {\dfrac{{{\partial ^2}{\Bbb{c}}}}{{\partial {\boldsymbol{\bar \varphi }}\partial {\boldsymbol{\bar \theta }}}}} \right)^{\!\rm{T}}}{\Bbb{f}} = \left[\boldsymbol{\mu}\right]_\times\left[I^\alpha\boldsymbol{\bar t}_\alpha\right]_\times,
\end{align}
and
\begin{align}
    \label{k_tht_tht_decompose}
    {{\boldsymbol{k}}_{{\rm{\bar \theta \bar \theta }}}} \coloneqq {\left( {\dfrac{{{\partial ^2}{\Bbb{c}}}}{{\partial {\boldsymbol{\bar \theta }}\partial {\boldsymbol{\bar \theta }}}}} \right)^{\!\rm{T}}}{\Bbb{f}} = {\boldsymbol{k}}_{{\rm{\bar \theta \bar \theta }}}^{{\rm{mat}}} + {\boldsymbol{k}}_{{\rm{\bar \theta \bar \theta }}}^{{\rm{geo}}},
\end{align}
\begin{subequations}
with the material part
\begin{align}
    {\boldsymbol{k}}_{{\rm{\bar \theta \bar \theta }}}^{{\rm{mat}}} = {\left[ {{I^\alpha }{{{\boldsymbol{\bar t}}}_\alpha }} \right]_ \times } {\left[{\boldsymbol{\mu }}\right]_\times\left[{s^i}{{{\boldsymbol{\bar t}}}_i}\right]_\times}  - {\left[ {{s^i}{{{\boldsymbol{\bar t}}}_i}} \right]_ \times }\left[{\boldsymbol{\mu }}\right]_\times\left[ {{I^\alpha }{{{\boldsymbol{\bar t}}}_\alpha }} \right]_\times,
\end{align}
and the geometric part
\begin{align}
    {\boldsymbol{k}}_{{\rm{\bar \theta \bar \theta }}}^{{\rm{geo}}}=\left[{\Bbb{\bar m}}_{\rm{t}}^i\right]_\times\left[{{{\boldsymbol{\bar t}}}_i}\right]_\times.
\end{align}
\end{subequations}
Note that $\boldsymbol{k}^\mathrm{mat}_{\bar\theta\bar\theta}$ is symmetric, and it can be shown that the unsymmetric part of ${\boldsymbol{k}}_{{\rm{\bar \theta \bar \theta }}}^{{\rm{geo}}}$ vanishes at the moment equilibrium. First, we have the unsymmetric part
\begin{align}
    \boldsymbol{k}^*_{\bar\theta\bar\theta}\coloneqq\dfrac{1}{2} \left({{\boldsymbol{k}}^\mathrm{geo}_{{\rm{\bar \theta \bar \theta }}}} - {{\boldsymbol{k}}^{\mathrm{geo}\,\mathrm{T}}_{{\rm{\bar \theta \bar \theta }}}}\right).
\end{align}
Then, we obtain
\begin{align}
    \delta {\boldsymbol{\bar \theta }} \cdot {{\boldsymbol{k}}^*_{{\rm{\bar \theta \bar \theta }}}}\Delta {\boldsymbol{\bar \theta }}&= \left( {\delta {\boldsymbol{\bar \theta }} \times \Delta {\boldsymbol{\bar \theta }}} \right) \cdot \left(\boldsymbol{\bar t}_i \times \Bbb{\bar m}^i_\mathrm{t}\right) =  \left( {\delta {\boldsymbol{\bar \theta }} \times \Delta {\boldsymbol{\bar \theta }}} \right) \cdot \Bbb{\bar m}.
\end{align}
This means that the unsymmetric part vanishes, that is $\boldsymbol{k}_{\bar\theta\bar\theta}$ becomes symmetric,  at the moment equilibrium, i.e., $\Bbb{\bar m}=\boldsymbol{0}$, if no external moment is applied.
\subsection{Expressions for the first-order beam kinematics}
\label{app_beam_new_oper}
Then, at $s=\bar s$, we have 
\begin{align}
    \label{x_sbar_oper_Pi_sbar}
    {{\boldsymbol{x}}_{\bar s}} = {{\boldsymbol{\Pi }}^{\rm{T}}}{{\boldsymbol{y}}_{\bar s}}.
\end{align} 
Using Eq.\,(\ref{x_sbar_oper_Pi}), the constraint equations in Eq.\,(\ref{cnst_cc_xs}) can be also rewritten, as 
\begin{align}
    \Bbb{\bar c}\equiv\Bbb{c}\left(\boldsymbol{y}_{\bar s},\boldsymbol{\bar q},\boldsymbol{s}\right)=\boldsymbol{0}.
\end{align}    
To correctly represent the stiffness, we need to further enhance the cross-sectional deformation modes. This will be discussed in the following section. For completeness, here, we rewrite $\delta \mathcal{J}$ in Eq.\,(\ref{fvar_J_cnst_jcb_rewrite_fc}) and its increment in Eq.\,(\ref{del_J_tangent_k}), by substituting Eq.\,(\ref{x_sbar_oper_Pi}) as
\begin{align}
    \label{fvar_J_cnst_jcb_rewrite}
	\delta \mathcal{\bar J} &={\delta {{\boldsymbol{y}}_{\bar s}} \cdot {{\boldsymbol{Q}}_{\rm{y}}}} + {\delta {\Bbb{\bar q}} \cdot \boldsymbol{Q}_{\Bbb{\bar q}}} + {\delta {\boldsymbol{s}} \cdot {{\boldsymbol{Q}}_{\rm{s}}}} + {\delta {\Bbb{f}} \cdot {\Bbb{\bar c}}},
\end{align}
where we have defined the generalized constraint force
\begin{align}
\boldsymbol{Q}_\mathrm{y}\coloneqq\boldsymbol{G}_\mathrm{y}^{\,\mathrm{T}}\Bbb{f}=\left\{ {\renewcommand{\arraystretch}{1.5}\begin{array}{*{20}{c}}
{A{\boldsymbol{\lambda }}}\\
{{\boldsymbol{\mu }} \times {I^{1\alpha }}{{\boldsymbol{t}}_\alpha }}\\
{{\boldsymbol{\mu }} \times {I^{2\alpha }}{{\boldsymbol{t}}_\alpha }}
\end{array}} \right\},
\end{align}
with the constraint Jacobian
\begin{align}
    \boldsymbol{G}_\mathrm{y}\coloneqq\dfrac{\partial \Bbb{\bar c}}{\partial \boldsymbol{y}}=\boldsymbol{G}_\mathrm{x}\boldsymbol{\Pi}^\mathrm{T}.
\end{align}
Further, in the same way, its increment can be also rewritten, as
\begin{align}
	\label{del_J_tangent_k_beam}
		\Delta \delta \mathcal{\bar J} =\left\{ {\renewcommand{\arraystretch}{1.5}\begin{array}{*{20}{c}}
				{\delta {\boldsymbol{y}}_\mathrm{\bar s}}\\
                {\delta {\boldsymbol{s}}}\\
                {\delta {\Bbb{\bar q}}}\\
                {\delta {\Bbb{f}}}
		\end{array}} \right\}.\left[ {\renewcommand{\arraystretch}{1.5}\begin{array}{*{20}{c}}
{{\bf{0}}_{d\times{d}}}&{{\bf{0}}_{d\times{3}}}&{{{\boldsymbol{k}}_{{\rm{y}\Bbb{\bar q}}}}}&{{{\boldsymbol{G}}^\mathrm{T}_{{\rm{y}}}}}\\
{{\bf{0}}_{3\times{d}}}&{{\bf{0}}_{3\times{3}}}&{{{\boldsymbol{k}}_{\rm{s}{\Bbb{\bar q}}}}}&{\boldsymbol{G}^\mathrm{T}_{\mathrm{s}}}\\
{{\boldsymbol{k}}_{{\rm{y}\Bbb{\bar q}}}^{\rm{T}}}&{{\boldsymbol{k}}_{\rm{s}{\Bbb{\bar q}}}^{\rm{T}}}&{{{\boldsymbol{k}}_{{\Bbb{\bar q}\Bbb{\bar q}}}}}&{\boldsymbol{G}^\mathrm{T}_{\Bbb{\bar q}}}\\
{{\boldsymbol{G}}_{\rm{y}}}&{{\boldsymbol{G}}_{{\rm{s}}}}&{{\boldsymbol{G}}_{{\Bbb{\bar q}}}}&{{{\bf{0}}_{m\times{m}}}}
\end{array}} \right]\left\{ {\renewcommand{\arraystretch}{1.5}\begin{array}{*{20}{c}}
				{\Delta {\boldsymbol{y}}_\mathrm{\bar s}}\\
                {\Delta {\boldsymbol{s}}}\\               
                {\Delta {\Bbb{\bar q}}}\\
				{\Delta {\Bbb{f}}}
		\end{array}} \right\},
\end{align}
with $\boldsymbol{k}_{\mathrm{y}\Bbb{\bar q}}\coloneqq\boldsymbol{\Pi}\boldsymbol{k}_{\mathrm{x}\Bbb{\bar q}}$. Here, more detailed expressions of $\boldsymbol{G}_\mathrm{y}$ and $\boldsymbol{k}_{\mathrm{y}\Bbb{\bar q}}$ can be found in Appendix \ref{app_beam_new_oper}.
The only new components in Eq.\,(\ref{del_J_tangent_k_beam}), compared with the tangent operator in Eq.\,(\ref{del_J_tangent_k}) for a brick formulation, are 
\begin{align}         \boldsymbol{G}_\mathrm{y}\coloneqq\boldsymbol{G}_\mathrm{x}\boldsymbol{\Pi}^\mathrm{T}=\left[ {\renewcommand{\arraystretch}{1.5}\begin{array}{*{20}{c}}
{A{\bf{1}}}\\
{{{\left[ {I^{1\alpha }}{{{\boldsymbol{t}}_\alpha }} \right]}_ \times }}\\
{{{\left[ {I^{2\alpha }}{{{\boldsymbol{t}}_\alpha }} \right]}_ \times }}
\end{array}} \right],
\end{align}
and
\begin{align}
    {{\boldsymbol{k}}_{{\rm{y\bar \theta }}}} \coloneqq  - \left[ {\renewcommand{\arraystretch}{1.5}\begin{array}{*{20}{c}}
    {{{\left[ {\boldsymbol{\mu }} \right]}_ \times }{{\left[ {{I^\alpha }{{\boldsymbol{t}}_\alpha }} \right]}_ \times }}\\
    {{{\left[ {\boldsymbol{\mu }} \right]}_ \times }{{\left[ {{I^{1\beta }}{{\boldsymbol{t}}_\beta }} \right]}_ \times }}\\
    {{{\left[ {\boldsymbol{\mu }} \right]}_ \times }{{\left[ {{I^{2\beta }}{{\boldsymbol{t}}_\beta }} \right]}_ \times }}
    \end{array}} \right]{{\boldsymbol{R}}_0}.
\end{align}
\setcounter{remark}{0}
\section{Numerical examples}
\subsection{Torsion of a straight beam with Z-section}	
\label{app_torsion_Z_shape}
\subsubsection{Analytical solution for torsional stiffness}
In isotropic linear elasticity, the torsional angle (radians) of a bar under a twisting moment $M_\mathrm{T}$ can be expressed by 
\begin{align}
    \label{app_z_shape_correct_K_asol}
    \theta = \dfrac{M_\mathrm{T}L}{KG},		
\end{align}
where $K$ denotes a torsional stiffness constant. For a Z-shaped cross-section, illustrated in Fig.\,\ref{z_shape_without_fillet}, we have an analytical solution of $K$ \citep{budynas2020roark}, as
\begin{align}
    \label{app_z_shape_correct_K}
    K = \frac{{{t^3}}}{3}\left( {2b + h} \right),
\end{align} 
where $t\ll{b,h}$ is assumed. This results in the linear curve in Fig.\,\ref{Z_torsion_mnt_case_1} (pink curve).
\subsubsection{Convergence test for brick solution}
\label{app_z_shape_conv_test_brick}
In Tables \ref{app_num_ex_z_shape_tab_conv_case1} and \ref{app_num_ex_z_shape_tab_conv_case2}, we show the convergence of the brick solutions for the applied moment ($M$) and the $X$-displacement at the end point B as the number of elements increases. \textcolor{blue}{For the brick solutions, we have used cubic NURBS basis functions in all directions: longitudinal (L), and transverse (W,H) ones}. It is seen that the solutions converge monotonically in all cases.
\begin{table}[htbp]
  \centering
  \footnotesize    
  \caption{Torsion of a Z-section beam (Case 1: no fillet). Convergence test for the brick solutions.}
    \begin{tabular}{cccrcc}
    \toprule
          & \multicolumn{2}{c}{Brick (partially clamped)} &       & \multicolumn{2}{c}{Brick (fully clamped)} \\
\cmidrule{2-3}\cmidrule{5-6}    $\left(n^\mathrm{L}_\mathrm{el},n^\mathrm{H}_\mathrm{el},n^\mathrm{W}_\mathrm{el}\right)$ & $M$ $[\mathrm{N}\cdot\mathrm{m}]$ & $u_\mathrm{B}$ $[\mathrm{m}]$ &       & $M$ $[\mathrm{N}\cdot\mathrm{m}]$ & $u_\mathrm{B}$ $[\mathrm{m}]$ \\
    \midrule
    $(30,3,3)$ & 1.5496E+06 & -5.2957E-02 &       & 1.9739E+06 & -5.0917E-02 \\
    $(40,4,4)$ & 1.5478E+06 & -5.2948E-02 &       & 1.9706E+06 & -5.0880E-02 \\
    $(50,5,5)$ & 1.5468E+06 & -5.2939E-02 &       & 1.9686E+06 & -5.0856E-02 \\
    \bottomrule
    \end{tabular}%
  \label{app_num_ex_z_shape_tab_conv_case1}%
\end{table}%

\begin{table}[htbp]
  \centering
  \footnotesize  
  \caption{Torsion of a Z-section beam (Case 2: added fillet). Convergence test for the brick solutions.}
    \begin{tabular}{cccrcc}
        \toprule
          & \multicolumn{2}{c}{Brick (partially clamped)} &       & \multicolumn{2}{c}{Brick (fully clamped)} \\
\cmidrule{2-3}\cmidrule{5-6}    $\left(n^\mathrm{L}_\mathrm{el},n^\mathrm{H}_\mathrm{el},n^\mathrm{W}_\mathrm{el}\right)$ & $M$ $[\mathrm{N}\cdot\mathrm{m}]$ & $u_\mathrm{B}$ $[\mathrm{m}]$ &       & $M$ $[\mathrm{N}\cdot\mathrm{m}]$ & $u_\mathrm{B}$ $[\mathrm{m}]$ \\
    \midrule
    $(30,3,3)$ & 1.7972E+06 & -5.2147E-02 &       & 2.2487E+06 & -5.0828E-02 \\
    $(40,4,4)$ & 1.7968E+06 & -5.2157E-02 &       & 2.2476E+06 & -5.0823E-02 \\
    $(50,5,5)$ & 1.7967E+06 & -5.2160E-02 &       & 2.2470E+06 & -5.0821E-02 \\
    \bottomrule
    \end{tabular}%
  \label{app_num_ex_z_shape_tab_conv_case2}%
\end{table}%

\subsection{Right-angle frame with a prismatic joint}
\label{app_num_ex_right_frame}
%
In Tables \ref{tab_frame_brick_conv_rigid} and \ref{tab_frame_brick_conv_prismatic}, we show the convergence of the brick solutions for the applied moment ($M$) and the $Y$-displacement at A ($v_\mathrm{A}$) as the number of elements increases. \textcolor{blue}{For the brick solutions, we have used cubic B-spline basis functions in all cases}. It is seen that these solutions converge monotonically in all cases.  
\begin{table}[H]
  \centering
  \footnotesize  
  \caption{Right-angle frame (Case 1: rigid joint). Convergence of the brick solutions.}
    \begin{tabular}{cccccc}
    \toprule
         & \multicolumn{2}{c}{Brick (partially clamped)} &       & \multicolumn{2}{c}{Brick (fully clamped)} \\
\cmidrule{2-3}\cmidrule{5-6}    $(n_\mathrm{el}^\mathrm{L},n_\mathrm{el}^\mathrm{W},n_\mathrm{el}^\mathrm{H})$ & $M$ $[\mathrm{N}\cdot\mathrm{m}]$  & $v_\mathrm{A}$ $[\mathrm{m}]$&       & $M$ $[\mathrm{N}\cdot\mathrm{m}]$  & $v_\mathrm{A}$ $[\mathrm{m}]$\\
    \midrule
    $(30,3,3)$ & 4.2059E+02 & 4.9154E-02 &       & 4.3643E+02 & 5.2022E-02 \\
    $(40,4,4)$ & 4.2043E+02 & 4.9128E-02 &       & 4.3627E+02 & 5.2003E-02 \\
    $(50,5,5)$ & 4.2037E+02 & 4.9118E-02 &       & 4.3620E+02 & 5.1993E-02 \\
    \bottomrule
    \end{tabular}%
  \label{tab_frame_brick_conv_rigid}%
\end{table}%

\begin{table}[H]
  \centering
  \footnotesize
  \caption{Right-angle frame (Case 2: prismatic joint). Convergence of the brick solutions.}
    \begin{tabular}{cccccc}
    \toprule
          & \multicolumn{2}{c}{Brick (partially clamped)} &       & \multicolumn{2}{c}{Brick (fully clamped)} \\
\cmidrule{2-3}\cmidrule{5-6}    $(n_\mathrm{el}^\mathrm{L},n_\mathrm{el}^\mathrm{W},n_\mathrm{el}^\mathrm{H})$ & $M$ $[\mathrm{N}\cdot\mathrm{m}]$   & $v_\mathrm{A}$ $[\mathrm{m}]$&       & $M$ $[\mathrm{N}\cdot\mathrm{m}]$  & $v_\mathrm{A}$ $[\mathrm{m}]$\\
    \midrule
    $(30,3,3)$ & 4.1812E+02 & 1.6289E-01 &       & 4.3395E+02 & 1.7884E-01 \\
    $(40,4,4)$ & 4.1797E+02 & 1.6281E-01 &       & 4.3377E+02 & 1.7876E-01 \\
    $(50,5,5)$ & 4.1791E+02 & 1.6277E-01 &       & 4.3370E+02 & 1.7873E-01 \\
    \bottomrule
    \end{tabular} %
  \label{tab_frame_brick_conv_prismatic} %
\end{table} %
\subsection{Shallow dome frame}
In Table \ref{tab_app_shallow_dome_brick_conv}, we show the convergence of the brick solutions for the applied force ($F$) and the rotation angle ($\theta$), as the number of elements increases, where we have used cubic B-spline basis functions in all cases. It is seen that the solutions converge monotonically in all cases. 
\begin{table}[H]
  \centering
  \footnotesize
  \caption{Shallow dome. Convergence of the brick solutions. \textcolor{blue}{Here, $F$ and $\theta$ are from Cases 1 and 2, respectively.}}
  \label{tab_app_shallow_dome_brick_conv}
    \begin{tabular}{cccccc}
    \toprule
          & \multicolumn{2}{c}{Brick} &       & \multicolumn{2}{c}{Brick (fully clamped)} \\
\cmidrule{2-3}\cmidrule{5-6}    $(n_\mathrm{el}^\mathrm{L},n_\mathrm{el}^\mathrm{W},n_\mathrm{el}^\mathrm{H})$ & $F$ $[\mathrm{N}]$ & $\theta$ $[\mathrm{deg}]$ &       & $F$ $[\mathrm{N}]$ & $\theta$ $[\mathrm{deg}]$\\
    \midrule
    $(30,3,3)$ & 4.9387E+02 & 1.2128E+01 &       & 4.9554E+02 & 1.1998E+01 \\
    $(40,4,4)$ & 4.9380E+02 & 1.2130E+01 &       & 4.9530E+02 & 1.2001E+01 \\
    $(50,5,5)$ & 4.9377E+02 & 1.2131E+01 &       & 4.9519E+02 & 1.2003E+01 \\
    \bottomrule
    \end{tabular}%
  \label{tab:addlabel}%
\end{table}%

\footnotesize
\bibliographystyle{elsarticle-num-names} 
\bibliography{cas-refs}

@book{hughes2003finite,
  title={The finite element method: linear static and dynamic finite element analysis},
  author={Hughes, Thomas JR},
  year={2003},
  publisher={Courier Corporation}
}

@book{Wagner1929,
  author    = {Wagner, Herbert},
  title     = {Verdrehen und Knicken von offenen Profilen (in German)},
  publisher = {Technische Hochschule Danzig},
  year      = {1929},
  address   = {Danzig},
  language  = {german}
}

@article{betsch1996four,
  title={A 4-node finite shell element for the implementation of general hyperelastic 3{D}-elasticity at finite strains},
  author={Betsch, P and Gruttmann, F and Stein, E},
  journal={Computer Methods in Applied Mechanics and Engineering},
  volume={130},
  number={1-2},
  pages={57--79},
  year={1996},
  publisher={Elsevier}
}

@article{bischoff1997shear,
  title={Shear deformable shell elements for large strains and rotations},
  author={Bischoff, M and Ramm, E},
  journal={International Journal for Numerical Methods in Engineering},
  volume={40},
  number={23},
  pages={4427--4449},
  year={1997},
  publisher={Wiley Online Library}
}

@article{buchter1994three,
  title={Three-dimensional extension of non-linear shell formulation based on the enhanced assumed strain concept},
  author={B{\"u}chter, Norbert and Ramm, Ekkehard and Roehl, Deane},
  journal={International Journal for Numerical Methods in Engineering},
  volume={37},
  number={15},
  pages={2551--2568},
  year={1994},
  publisher={Wiley Online Library}
}

@article{antman1966dynamical,
  title={Dynamical theory of hyperelastic rods},
  author={Antman, Stuart S and Warner, William H},
  journal={Archive for Rational Mechanics and Analysis},
  volume={23},
  number={2},
  pages={135--162},
  year={1966},
  publisher={Springer}
}

@book{angeles1982spatial,
  title={Spatial Kinematic Chains: Analysis—Synthesis—Optimization},
  author={Angeles, Jorge},
  year={1982},
  publisher={Springer Science \& Business Media}
}

@article{coda2009solid,
  title={A solid-like {FEM} for geometrically non-linear 3{D} frames},
  author={Coda, Humberto Breves},
  journal={Computer Methods in Applied Mechanics and Engineering},
  volume={198},
  number={47-48},
  pages={3712--3722},
  year={2009},
  publisher={Elsevier}
}

@article{greco2021non,
  title={A non-linear symmetric ${G}^1$-conforming {B}{\'e}zier finite element formulation for the analysis of {K}irchhoff beam assemblies},
  author={Greco, L and Scrofani, A and Cuomo, M},
  journal={Computer Methods in Applied Mechanics and Engineering},
  volume={387},
  pages={114176},
  year={2021},
  publisher={Elsevier}
}

@article{kiendl2010bending,
  title={The bending strip method for isogeometric analysis of {K}irchhoff--{L}ove shell structures comprised of multiple patches},
  author={Kiendl, Josef and Bazilevs, Yuri and Hsu, M-C and W{\"u}chner, Roland and Bletzinger, K-U},
  journal={Computer Methods in Applied Mechanics and Engineering},
  volume={199},
  number={37-40},
  pages={2403--2416},
  year={2010},
  publisher={Elsevier}
}

@article{schuss2019multi,
  title={Multi-patch isogeometric analysis for {K}irchhoff--{L}ove shell elements},
  author={Schu{\ss}, Stefan and Dittmann, Maik and Wohlmuth, Barbara and Klinkel, Sven and Hesch, Christian},
  journal={Computer Methods in Applied Mechanics and Engineering},
  volume={349},
  pages={91--116},
  year={2019},
  publisher={Elsevier}
}

@article{bae2000implementation,
  title={An implementation method for constrained flexible multibody dynamics using a virtual body and joint},
  author={Bae, DS and Han, JM and Choi, JH},
  journal={Multibody System Dynamics},
  volume={4},
  number={4},
  pages={297--315},
  year={2000},
  publisher={Springer}
}

@article{duong2017new,
  title={A new rotation-free isogeometric thin shell formulation and a corresponding continuity constraint for patch boundaries},
  author={Duong, Thang X and Roohbakhshan, Farshad and Sauer, Roger A},
  journal={Computer Methods in Applied Mechanics and Engineering},
  volume={316},
  pages={43--83},
  year={2017},
  publisher={Elsevier}
}

@article{gerstmayr2013review,
  title={Review on the absolute nodal coordinate formulation for large deformation analysis of multibody systems},
  author={Gerstmayr, Johannes and Sugiyama, Hiroyuki and Mikkola, Aki},
  journal={Journal of Computational and Nonlinear Dynamics},
  volume={8},
  number={3},
  pages={031016},
  year={2013},
  publisher={American Society of Mechanical Engineers}
}

@article{bauer2020weak,
  title={Weak coupling of nonlinear isogeometric spatial {B}ernoulli beams},
  author={Bauer, AM and W{\"u}chner, R and Bletzinger, K-U},
  journal={Computer Methods in Applied Mechanics and Engineering},
  volume={361},
  pages={112747},
  year={2020},
  publisher={Elsevier}
}

@article{munoz2003master,
  title={Master--slave approach for the modelling of joints with dependent degrees of freedom in flexible mechanisms},
  author={Mu{\~n}oz, J and Jeleni{\'c}, Gordan and Crisfield, MA},
  journal={Communications in Numerical Methods in Engineering},
  volume={19},
  number={9},
  pages={689--702},
  year={2003},
  publisher={Wiley Online Library}
}

@article{jelenic2001dynamic,
  title={Dynamic analysis of 3{D} beams with joints in presence of large rotations},
  author={Jeleni{\'c}, Gordan and Crisfield, MA},
  journal={Computer Methods in Applied Mechanics and Engineering},
  volume={190},
  number={32-33},
  pages={4195--4230},
  year={2001},
  publisher={Elsevier}
}

@article{bauchau2008review,
  title={Review of contemporary approaches for constraint enforcement in multibody systems},
  author={Bauchau, Olivier A and Laulusa, Andr{\'e}},
  journal={Journal of Computational and Nonlinear Dynamics},
  volume={3},
  number={1},
  pages={011005},
  year={2008}
}

@article{leyendecker2008discrete,
  title={The discrete null space method for the energy-consistent integration of constrained mechanical systems. {P}art {III}: Flexible multibody dynamics},
  author={Leyendecker, Sigrid and Betsch, Peter and Steinmann, Paul},
  journal={Multibody System Dynamics},
  volume={19},
  number={1},
  pages={45--72},
  year={2008},
  publisher={Springer}
}

@article{betsch2006discrete,
  title={The discrete null space method for the energy consistent integration of constrained mechanical systems. {P}art {II}: Multibody dynamics},
  author={Betsch, Peter and Leyendecker, Sigrid},
  journal={International Journal for Numerical Methods in Engineering},
  volume={67},
  number={4},
  pages={499--552},
  year={2006},
  publisher={Wiley Online Library}
}

@book{bertsekas2014constrained,
  title={Constrained optimization and Lagrange multiplier methods},
  author={Bertsekas, Dimitri P},
  year={2014},
  publisher={Academic press}
}

@article{di2024corotational,
  title={Corotational force-based beam finite element with rigid joint offsets for 3{D} framed structures},
  author={Di Re, Paolo and Addessi, Daniela and Gatta, Cristina and Parente, Luca and Sacco, Elio},
  journal={Computer Methods in Applied Mechanics and Engineering},
  volume={419},
  pages={116656},
  year={2024},
  publisher={Elsevier}
}

@article{hesch2009mortar,
  title={A mortar method for energy-momentum conserving schemes in frictionless dynamic contact problems},
  author={Hesch, Christian and Betsch, Peter},
  journal={International Journal for Numerical Methods in Engineering},
  volume={77},
  number={10},
  pages={1468--1500},
  year={2009},
  publisher={Wiley Online Library}
}

@article{greco2024objective,
  title={An objective minimal constraint formulation for the analysis of elastic articulated structures},
  author={Greco, L and Castello, D and Cuomo, M},
  journal={Computers \& Structures},
  volume={305},
  pages={107571},
  year={2024},
  publisher={Elsevier}
}

@article{hussein2009clamped,
  title={Clamped end conditions and cross section deformation in the finite element absolute nodal coordinate formulation},
  author={Hussein, Bassam A and Weed, David and Shabana, Ahmed A},
  journal={Multibody System Dynamics},
  volume={21},
  number={4},
  pages={375--393},
  year={2009},
  publisher={Springer}
}

@article{sugiyama2011spatial,
  title={Spatial joint constraints for the absolute nodal coordinate formulation using the non-generalized intermediate coordinates},
  author={Sugiyama, Hiroyuki and Yamashita, Hiroki},
  journal={Multibody System Dynamics},
  volume={26},
  number={1},
  pages={15--36},
  year={2011},
  publisher={Springer}
}

@article{hong2011modeling,
  title={A modeling of sliding joint on one-dimensional flexible medium},
  author={Hong, Difeng and Ren, Gexue},
  journal={Multibody System Dynamics},
  volume={26},
  number={1},
  pages={91--106},
  year={2011},
  publisher={Springer}
}

@article{munoz2008modelling,
  title={Modelling unilateral frictionless contact using the null-space method and cubic {B}-{S}pline interpolation},
  author={Mu{\~n}oz, JJ},
  journal={Computer Methods in Applied Mechanics and Engineering},
  volume={197},
  number={9-12},
  pages={979--993},
  year={2008},
  publisher={Elsevier}
}

@article{rees2018comparative,
  title={A comparative study of null-space factorizations for sparse symmetric saddle point systems},
  author={Rees, Tyrone and Scott, Jennifer},
  journal={Numerical Linear Algebra with Applications},
  volume={25},
  number={1},
  pages={e2103},
  year={2018},
  publisher={Wiley Online Library}
}

@article{wolfe1962reduced,
  title  = {The Reduced Gradient Method},
  author = {Wolfe, Philip},
  note   = {Unpublished manuscript, The RAND Corporation},
  year   = {1962}
}

@article{benzi2005numerical,
  title={Numerical solution of saddle point problems},
  author={Benzi, Michele and Golub, Gene H and Liesen, J{\"o}rg},
  journal={Acta Numerica},
  volume={14},
  pages={1--137},
  year={2005},
  publisher={Cambridge University Press}
}

@article{wehage1982generalized,
  title={Generalized coordinate partitioning for dimension reduction in analysis of constrained dynamic systems},
  author={Wehage, Roger A and Haug, Edward J},
  journal={Journal of Mechanical Design},
  volume={104},
  number={1},
  pages={247--255},
  year={1982}
}

@manual{onemkl_ref_2025,
  title        = {Intel® oneAPI Math Kernel Library (oneMKL) Developer Reference},
  author       = {{Intel Corporation}},
  year         = {2025},
  url          = {},
  note         = {}
}

@article{gerstmayr20063d,
  title={A 3{D} finite element method for flexible multibody systems},
  author={Gerstmayr, Johannes and Sch{\"o}berl, Joachim},
  journal={Multibody System Dynamics},
  volume={15},
  number={4},
  pages={305--320},
  year={2006},
  publisher={Springer}
}

@article{boungard2024master,
  title={Master--slave elimination scheme for arbitrary smooth nonlinear multi-point constraints},
  author={Boungard, Jonas and Wackerfu{\ss}, Jens},
  journal={Computational Mechanics},
  volume={74},
  number={5},
  pages={955--992},
  year={2024},
  publisher={Springer}
}

@article{sugiyama2003formulation,
  title={Formulation of three-dimensional joint constraints using the absolute nodal coordinates},
  author={Sugiyama, Hiroyuki and Escalona, Jos{\'e} L and Shabana, Ahmed A},
  journal={Nonlinear Dynamics},
  volume={31},
  number={2},
  pages={167--195},
  year={2003},
  publisher={Springer}
}

@article{carrera2010refined,
  title={Refined beam elements with arbitrary cross-section geometries},
  author={Carrera, Erasmo and Giunta, Gaetano and Nali, Pietro and Petrolo, Marco},
  journal={Computers \& Structures},
  volume={88},
  number={5-6},
  pages={283--293},
  year={2010},
  publisher={Elsevier}
}

@article{betsch2005discrete,
  title={The discrete null space method for the energy consistent integration of constrained mechanical systems: Part {I}: Holonomic constraints},
  author={Betsch, Peter},
  journal={Computer Methods in Applied Mechanics and Engineering},
  volume={194},
  number={50-52},
  pages={5159--5190},
  year={2005},
  publisher={Elsevier}
}

@article{leonetti2025locking,
  title={Locking mitigation in geometrically nonlinear solid beams through isogeometric analysis and a generalized constitutive approach},
  author={Leonetti, Leonardo and Patton, Alessia and Kiendl, Josef},
  journal={Computer Methods in Applied Mechanics and Engineering},
  volume={445},
  pages={118149},
  year={2025},
  publisher={Elsevier}
}

@article{argyris1982excursion,
  title={An excursion into large rotations},
  author={Argyris, John},
  journal={Computer Methods in Applied Mechanics and Engineering},
  volume={32},
  number={1-3},
  pages={85--155},
  year={1982},
  publisher={Elsevier}
}

@article{durville2012contact,
  title={Contact-friction modeling within elastic beam assemblies: an application to knot tightening},
  author={Durville, Damien},
  journal={Computational Mechanics},
  volume={49},
  number={6},
  pages={687--707},
  year={2012},
  publisher={Springer}
}

@article{tian2009dynamics,
  title={Dynamics of spatial flexible multibody systems with clearance and lubricated spherical joints},
  author={Tian, Qiang and Zhang, Yunqing and Chen, Liping and Flores, Paulo},
  journal={Computers \& Structures},
  volume={87},
  number={13-14},
  pages={913--929},
  year={2009},
  publisher={Elsevier}
}

@article{battini2002co,
  title={Co-rotational beam elements with warping effects in instability problems},
  author={Battini, Jean-Marc and Pacoste, Costin},
  journal={Computer Methods in Applied Mechanics and Engineering},
  volume={191},
  number={17-18},
  pages={1755--1789},
  year={2002},
  publisher={Elsevier}
}

@article{frischkorn2013solid,
  title={A solid-beam finite element and non-linear constitutive modelling},
  author={Frischkorn, J and Reese, S},
  journal={Computer Methods in Applied Mechanics and Engineering},
  volume={265},
  pages={195--212},
  year={2013},
  publisher={Elsevier}
}

@article{simo1985perturbed,
  title={A perturbed {L}agrangian formulation for the finite element solution of contact problems},
  author={Simo, Juan C and Wriggers, Peter and Taylor, Robert L},
  journal={Computer Methods in Applied Mechanics and Engineering},
  volume={50},
  number={2},
  pages={163--180},
  year={1985},
  publisher={Elsevier}
}

@article{gruttmann2000theory,
  title={Theory and numerics of three-dimensional beams with elastoplastic material behaviour},
  author={Gruttmann, F and Sauer, R and Wagner, W},
  journal={International Journal for Numerical Methods in Engineering},
  volume={48},
  number={12},
  pages={1675--1702},
  year={2000},
  publisher={Wiley Online Library}
}

@article{choi2022isogeometric,
  title={An isogeometric finite element formulation for frictionless contact of Cosserat rods with unconstrained directors},
  author={Choi, Myung-Jin and Klinkel, Sven and Sauer, Roger A},
  journal={Computational Mechanics},
  volume={70},
  number={6},
  pages={1107--1144},
  year={2022},
  publisher={Springer}
}

@article{rhim1998vectorial,
  title={A vectorial approach to computational modelling of beams undergoing finite rotations},
  author={Rhim, Jaewook and Lee, Sung W},
  journal={International Journal for Numerical Methods in Engineering},
  volume={41},
  number={3},
  pages={527--540},
  year={1998},
  publisher={Wiley Online Library}
}

@article{shafqat2024robust,
  title={A robust finite strain isogeometric solid-beam element},
  author={Shafqat, Abdullah and Weeger, Oliver and Xu, Bai-Xiang},
  journal={Computer Methods in Applied Mechanics and Engineering},
  volume={426},
  pages={116993},
  year={2024},
  publisher={Elsevier}
}

@article{duong2023variationally,
  title={A variationally consistent contact formulation based on a mixed interpolation point method and isogeometric discretization},
  author={Duong, Thang X and Leonetti, Leonardo and Kiendl, Josef},
  journal={Computer Methods in Applied Mechanics and Engineering},
  volume={417},
  pages={116361},
  year={2023},
  publisher={Elsevier}
}

@article{wackerfuss2009mixed,
  title={A mixed hybrid finite beam element with an interface to arbitrary three-dimensional material models},
  author={Wackerfu{\ss}, J and Gruttmann, F},
  journal={Computer Methods in Applied Mechanics and Engineering},
  volume={198},
  number={27-29},
  pages={2053--2066},
  year={2009},
  publisher={Elsevier}
}

@article{wackerfuss2011nonlinear,
  title={A nonlinear {H}u--{W}ashizu variational formulation and related finite-element implementation for spatial beams with arbitrary moderate thick cross-sections},
  author={Wackerfu{\ss}, J and Gruttmann, F},
  journal={Computer Methods in Applied Mechanics and Engineering},
  volume={200},
  number={17-20},
  pages={1671--1690},
  year={2011},
  publisher={Elsevier}
}

@Book{bonet2010nonlinear,
  Title                    = {Nonlinear continuum mechanics for finite element analysis},
  Author                   = {Bonet, Javier and Wood, Richard D},
  Publisher                = {Cambridge University Press},
  Year                     = {2008}
}

@article{hughes2005isogeometric,
  title={Isogeometric analysis: {CAD}, finite elements, {NURBS}, exact geometry and mesh refinement},
  author={Hughes, TJR and Cottrell, JA and Bazilevs, Y},
  journal={Computer Methods in Applied Mechanics and Engineering},
  volume={194},
  number={39-41},
  pages={4135--4195},
  year={2005},
  publisher={Elsevier}
}

@article{choi2021isogeometric,
  title={An isogeometric finite element formulation for geometrically exact {T}imoshenko beams with extensible directors},
  author={Choi, Myung-Jin and Sauer, Roger A. and Klinkel, S},
  journal={Computer Methods in Applied Mechanics and Engineering},
  volume={385},
  pages={113993},
  year={2021},
  publisher={Elsevier}
}

@article{simo1986three,
  title={A three-dimensional finite-strain rod model. {P}art {II}: {C}omputational aspects},
  author={Simo, Juan C and Vu-Quoc, Loc},
  journal={Computer Methods in Applied Mechanics and Engineering},
  volume={58},
  number={1},
  pages={79--116},
  year={1986},
  publisher={Elsevier}
}

@article{markovic2004micro,
  title={On micro--macro interface conditions for micro scale based {FEM} for inelastic behavior of heterogeneous materials},
  author={Markovic, Damijan and Ibrahimbegovic, Adnan},
  journal={Computer Methods in Applied Mechanics and Engineering},
  volume={193},
  number={48-51},
  pages={5503--5523},
  year={2004},
  publisher={Elsevier}
}

@article{choi2024objective,
  title={An objective isogeometric mixed finite element formulation for nonlinear elastodynamic beams with incompatible warping strains},
  author={Choi, Myung-Jin and Klinkel, Sven and Klarmann, Simon and Sauer, Roger A},
  journal={Multibody System Dynamics},
  pages={1--54},
  year={2024},
  publisher={Springer}
}

@article{jelenic1996non,
  title={Non-linear ‘master-slave’relationships for joints in 3-{D} beams with large rotations},
  author={Jeleni{\'c}, Gordan and Crisfield, MA},
  journal={Computer Methods in Applied Mechanics and Engineering},
  volume={135},
  number={3-4},
  pages={211--228},
  year={1996},
  publisher={Elsevier}
}

@article{choi2023selectively,
  title={A selectively reduced degree basis for efficient mixed nonlinear isogeometric beam formulations with extensible directors},
  author={Choi, Myung-Jin and Sauer, Roger A and Klinkel, Sven},
  journal={Computer Methods in Applied Mechanics and Engineering},
  volume={417},
  pages={116387},
  year={2023},
  publisher={Elsevier}
}

@phdthesis{klarmann2018geometrisch,
    author = {Klarmann, Simon},
    title = {Geometrisch und physikalisch nichtlineare Mehrskalenmodellierung r{\"a}umlicher Stabtragwerke (in German)},
    school = {Technical University of Darmstadt},
    year = {2018}
}

@book{budynas2020roark,
  title={Roark's formulas for stress and strain},
  author={Budynas, RG and Sadegh, AM},
  year={2020},
  publisher={McGraw-Hill Education}
}

@inproceedings{moustacas2019enrichissement,
  title={Enrichissement d’une cin{\'e}matique poutre Applications aux textiles en carbone (in {F}rench)},
  author={Moustacas, H and Durville, D and Wielhorski, Y},
  booktitle={14e Colloque National en Calcul des Structures (CSMA)},
  pages={1--8},
  year={2019}
}

@article{klarmann2020homogenization,
  title={Homogenization assumptions for coupled multiscale analysis of structural elements: beam kinematics},
  author={Klarmann, Simon and Gruttmann, Friedrich and Klinkel, Sven},
  journal={Computational Mechanics},
  volume={65},
  number={3},
  pages={635--661},
  year={2020},
  publisher={Springer}
}

\end{document}